\documentclass[12pt]{report}

\usepackage{fancyhdr}
\usepackage{epsfig}
\usepackage{amsmath,amssymb,amsfonts,textcomp,wasysym}
\usepackage{amscd}
\usepackage{dsfont}
\usepackage{amsthm}

\usepackage{sectsty}
\chapterfont{\large\sc\centering}
\chaptertitlefont{\centering}
\sectionfont{\centering}

\numberwithin{equation}{section}
\numberwithin{equation}{section}

\usepackage{xypic}

\theoremstyle{plain}

\theoremstyle{plain}

\newtheorem{theorem}{Theorem}[section]
\newtheorem{lemma}[theorem]{Lemma}
\newtheorem{proposition}[theorem]{Proposition}
\newtheorem{corollary}[theorem]{Corollary}
\newtheorem{definition}[theorem]{Definition}
\newtheorem{example}[theorem]{Example}
\newtheorem{note}{Note}[section]

\newenvironment{remark}[1][Remark]{\begin{trivlist}
\item[\hskip \labelsep {\bfseries #1}]}{\end{trivlist}}

\def \Z {\mathbb Z}
\def \N {\mathbb N}
\def \Q {\mathbb Q}
\def \eq {\equiv}
\def \zp {\Z_{(p)}}
\def \zn {\Z_{(n)}}
\def \qp {\Q_{(p)}}
\def \qn {\Q_{(n)}}
\def \r {2^{r-2}}
\def \w {2^{n-2}}

\def \v {2^{n-3}}

\def \n {\Z_{p^n} ^*}
\def \t {\Z_{2^n} ^*}
\def \p {p^{n - 1}}
\def \q {p^{n - 2}}
\def \a {\displaystyle\sum_{i = 1} ^ {\infty} a_i \, p^i}
\def \x {\displaystyle\sum_{i = 0} ^ {\infty} x_i \, p^i}
\def \b {\displaystyle\sum_{i = 1} ^ {\infty} b_i \, p^i}
\def \A {\sum_{i = 0} ^ {\infty} A_i \, x^i}
\def \B {\sum_{j = 1} ^ {\infty} B_j \, x^j}
\def \- {\, \in \,}
\def \h3 {\hspace{3em}}
\def \nr {p^{n - r}}
\def \nrone {p^{n - r - 1}}
\def \uone {1 + p \, \zp}
\def \dl {\displaystyle\left(}
\def \dr {\displaystyle\right)}
\def \iff {if and only if }

\begin{document}

{
\pagestyle{empty}
{\begin{center}
\vspace*{1cm}
{\bf \LARGE EXISTENCE AND DISTRIBUTION }\\
\vspace*{3mm}
{\bf \LARGE OF SOLUTIONS OF $a^x = b$ MODULO $p^n$}\\
%
\vspace*{1.cm}


{\bf Submitted in partial fulfillment of the requirements\\
of the degree of } \\
\vspace*{1.7cm}
\textbf{Master of Philosophy}\\
{\bf by}\\

{ \bf \large RUPALI RAVINDRA KHEDKAR \\
(Roll No.07409303)}\\
\vspace*{1cm}
{\bf Supervisor} \\ 
\vspace*{.5cm}
{\bf \large PROFESSOR RAVINDRA S. KULKARNI} \\
\vspace*{1.5cm}


\vspace*{1.5cm}
{\bf  DEPARTMENT OF MATHEMATICS}\\ 

{\bf  INDIAN INSTITUTE OF TECHNOLOGY BOMBAY\\ 21 JANUARY 2011 }
\end{center}
}

}

\pagenumbering{roman} 

\chapter*{Declaration}
  
\indent \textsc{I} declare that this written submission represents my ideas in my own words and where 
others' ideas or words have been included, I have adequately cited and referenced the original 
sources.  I also declare that I have adhered to all principles of academic honesty and integrity 
and   have   not   misrepresented   or   fabricated   or   falsified   any   idea/data/fact/source   in   my 
submission.  I understand that any violation of the above will be cause for disciplinary action 
by the Institute and can also evoke  penal action from the sources which have thus not been 
properly cited or from whom proper permission has not been taken when needed.

\vspace*{2cm}
\begin{flushright}
\textbf{-----------------------------------}\\
(Signature)
\vspace*{5mm}

\textbf{\underline{ RUPALI R. KHEDKAR. }}\\
 (Name of the student)
\vspace*{5mm}

\textbf{\underline{ 07409303 }}\\
(Roll No.)
\end{flushright}

Date:\, \textbf{\underline{$21^{st}$ Jan. 2011}}

\tableofcontents

\chapter*{Acknowledgement}
\indent \textsc{WORK} done in this thesis had started with some basic and fundamental questions in Number Theory and Group Theory.
Almost all the work done is an outcome of the discussions with my advisor Prof. Ravindra Kulkarni. 
His patience was really great, even when I was making silly mistakes or when I took a long time to answer some trivial questions. 
Even though sometimes I did not understand the matter, I really enjoyed each and every moment of those discussions.
He always gave me support and most important believed in my capabilities. 
I salute Prof. Kulkarni for his dedication towards mathematics and express heartfelt gratitude to him for everything he has done for me. \\
\indent I thank all my friends, Anuradha Ahuja, Saumya, Bakul, Shreedevi and Sayali for making my hostel life memorable and enjoyable.
Especially, I had very nice discussions with Anuradha not only on Mathematics but on Philosophy also.
I am very much thankful to Sayali and Mousumi for their kind help in Latex.\\
\indent It's beyond giving thanks to them and very difficult to express in words- I feel very fortunate that, I have such parents and grand parents who always encouraged 
me for higher studies and gave complete support all the time. I had long phone calls with the most lovable person, my younger sister, 
because of whom I used to forget all my tensions and get energy for doing further work. \\
\indent I salute the great mathematicians, Isaac Newton for his discovery of Binomial theorem and 
Carl Friedrich Gauss for his contribution to the primitive roots theorem.\\
\indent I thank National Board for Higher Mathematics for there financial support for 2 years.
I am very much thankful to the producers of LaTex, for creating such a nice and user friendly software.\\
\indent I dedicate this thesis to my parents, grand parents and my teachers. 

\chapter*{Abstract}

\indent  \textsc{INITIAL} objective of this dissertation is to study the existence of the solutions of the congruence $a^x \eq b^y \, (mod \, p^n)$
and distribution of solutions $(x, y)$ as $n$ varies in natural numbers, where $a$ and $b$ are integers coprime to prime $p.$
We observe that as $n \rightarrow \infty$, solutions take the form of $p$-adic integers. \\
\indent This motivates us, to study the existence of the solutions of equation $a^x = b$ in $p$-adic integers. 
 The relevant case is when $a$ and $b$ are units in $p$-adic integers. If the solution exists we try to find it out. We resolve the case of $a, b \- U_1$ completely.
 A necessary and sufficient condition for the existence of the solution of $a^x = b$ where $a, b$ are elements of $U_1,$ is `depth of $a$ is smaller than the depth of $b$'.
 In this case, if the solution exists then it is given by $\log b / \log a.$ In the other case, where $a$ and $b$ are $p$-adic units but not elements of $U_1$,
 we give the criteria for the existence of the solution of $a^x = b.$ Write $a$ and $b$ as product of Teichm\"{u}ller unit and an element of $U_1$.
 Suppose $a= a_1 a_2$ and $b = b_1 b_2$ where $a_1, b_1$ are Teichm\"{u}ller units and $a_2, b_2 \- U_1.$ Then the solution of $a^x = b$ exists if $b_1$ belongs
 to the group generated by $a_1$ and depth of $a_2$ is smaller than the depth of $b_2.$ 

\pagenumbering{arabic}

\chapter{INTRODUCTION} 

\section{Motivation}

\indent Suppose $a$, $b$ are integers and $p$ is a prime, consider the congruence $a^x \eq b^y \, (mod \, p^n) \hspace*{2em} (*)$ where $n$ is a natural number. 
We try to find when does the solution exist for $(*).$ Since this is a multiplicative congruence, solution in $(x, y)$ exists \iff the two sets 
$A = \{a, a^2, a^3, \cdots \} \, modulo\,  p^n$ and $B = \{ b, b^2, b^3, \cdots \} \, modulo\,  p^n$ intersects. In general if $(a, p) \neq 1$ or $(b, p) \neq 1$ then
 the sets $A$, $B$ are not groups under multiplication modulo $p^n.$ We will show that, the case $a$ and $b$ both are coprime to $p$ is the only interesting case to study 
the existence of solutions of $(*).$ In this case the sets $A \, (mod \, p^n)$ and $B \, (mod \, p^n)$ are actually subgroups of $(\n, *).$ Therefore, if $(x_o, y_o)$ satisfies
the congruence $(*)$ then $x_o \- \{1, 2, \cdots, o(a)\}$ and $y_o \- \{1, 2, \cdots, o(b)\}.$ If we fix $y = 1$ in $(*)$ and find the solution of $a^x \eq b \, (mod \, p^n) $ say
$(x_o, 1)$ then $a^{x_o y} \eq b^y \, (mod \, p^n)$ gives all the solutions of $(*).$ Hence we reduce to the case of studying solutions of $a^x \eq b \, (p^n)$ where $a$ and $b$
are coprime to $p.$ Let $x_n$ denote the solution of $a^x \eq b \, (p^n).$ To investigate the pattern of these solutions was a primary motivation of this dissertation.
 We will show that if $x_n$ exists for all $n$ then the sequence $\{x_n\}$ converges in $p$-adic integers.
That means the difference $x_{n + 1} - x_n$ is divisible by $p^{k_n}$ where $k_n$ goes to infinity as $n \rightarrow \infty.$ More precisely, we state the theorem \ref{25}
in section \ref{E}: \\
\indent {\it{Let $a$, $b$ be integers coprime to a prime $p$ and $a \neq 1.$ Let $x_n$ denote the smallest positive solution of the congruence $a^x \eq b \, (mod \, p^n)$.
If $x_n$ exists for all $n \geq 1$ then the sequence $\{x_n\}$ converges in $p$-adic integers. }}
\newpage
\section{Aim}

\indent The main objective of this dissertation is to study the existence and distribution of the solutions of equation $a^x = b$ modulo $p^n$,
where $a, b$ are units in $p$-adic integers, $a \neq 1$ and $n$ goes to infinity through natural numbers.
We divide the work in three parts:\\
 1. $a$ and $b$ are integers coprime to $p$, \\
 2. $a$ and $b$ are elements of $U_1,$\\
 3. $a, b$ are units in $p$-adic integers.

\indent We state here the main theorems in each of these cases: In the first case refer the theorem \ref{25}, stated as:\\
\noindent {\it{``Let $a$, $b$ be integers coprime to a prime $p$ and $a \neq 1.$ Let $x_n$ denote the smallest positive solution of the congruence $a^x \eq b \, (mod \, p^n)$.
If $x_n$ exists for all $n \geq 1$ then the sequence $\{x_n\}$ converges in $p$-adic integers.'' }}

\noindent In the second case, refer theorem \ref{27} from the section \ref{F} which gives the necessary and sufficient condition for the existence of $p$-adic solution,
when $a$ and $b$ are in $U_1.$\\
\noindent {\it{ ``Let $a, b $ be elements of $U_1$ and $a \neq 1$. Then $a^x = b$ has unique solution in $p$-adic integers if and only if depth of $a \leq$ depth of $b$.''}} 

\noindent For the third case, the main theorem \ref{28} is stated as below:\\
\noindent {\it{``Let $p$ be an odd prime. Let $a, b$ be the units in $\zp.$ Write $a = a_1 a_2$, $b = b_1 b_2$ where $a_1, b_1$ are Teichm\"{u}ller units and
 $a_2, b_2 \- U_1$ then $a^x = b$ has a solution in $\zp$ if (i) $b_1$ belong to the group generated by $a_1,$ and (ii) depth of $a_2 \leq$ depth of $b_2.$ \\
Let $p = 2$, Let $a, b$ be the units in $\Z_{(2)}.$ Write $a = a_1 a_2$, $b = b_1 b_2$ where $a_1, b_1$ are Teichm\"{u}ller units and $a_2, b_2 \- U_1$ then $a^x = b$ has a 
solution in $\Z_{(2)}$ if (i) $a_1 = b_1$ and (ii) depth of $a_2 \leq$ depth of $b_2.$ '' }}

The main examples for all the work in chapter 4 are $p = 2, (a = -3, b = 5 )$, for $p$ an odd prime, $(a = 1 - p, b = 1 + p).$ We discuss them in detail in section \ref{E},
example \ref{9} and section \ref{F}, example \ref{10}. These examples are also special pairs which we define at the end in section \ref{I}. 
\section{Preliminaries}

In this section we will highlight the main points from chapter 2 and 3. \\
\indent One of the most important contribution of Gauss in Disquisitiones Arithmaticae \cite{Gau} is his result on primitive roots. The result is stated as:
`A natural number $n$ has a primitive root \iff $n$ is 2 or 4 or of the form $p^m$ or $2 p^m$ for some odd prime $p.$' In modern group theoretic language,
theorem can be stated as: A multiplicative group $\Z_n ^*$ is cyclic \iff $n = 2, 4, p^m, 2 p^m$ for an odd prime $p.$ A constructible proof of the assertion
`every prime has a primitive root' is not known. If $r$ is a primitive root of both $p$ and $p^2$ then $r$ is a primitive root of $p^n$ for all $n.$ 
If $r$ is a primitive root of $p$ then $r(1 + p)$ and $r + p$ both are primitive roots of $p^2,$ hence they are also primitive root of $p^n$ for all $n \geq 1.$
In section \ref{G}, we discuss this in more detail. We in fact prove that, there exists an integer between $2$ and $p - 2$ which is a primitive root of $p^n$ for all $n \geq 1.$

In chapter 3, we define $n$-adic integers in two ways, one is using projective limit and the other is metric approach. In section \ref{D}, we will explicitly give how to construct 
the Teichm\"{u}ller units in $p$-adic integers. The structure of units in $p$-adic integers is important for our work in section \ref{H}. In the theory of quadratic forms,
the structure of $\qp ^* / {\qp ^*}^2$ is important. We will generalize this by considering the map $x \mapsto x^k$ on $\qp ^*$ and will study the structure of  $\qp ^* / {\qp ^*}^k.$

\chapter{THE GROUP $\Z_n ^*$}

\textsc{IN} this chapter we cover the basic material of group theory which we will need throughout the thesis. Section \ref{A} contains various cases in $a, b$ in the congruence 
$a^x \eq b^y \, (mod \, p^n)$. In section \ref{B} we will see some properties of the group $\Z_n ^*$, proof of structure theorem for $\Z_n ^*$ and Gauss's original proof of
`existence of primitive root of a prime.' In the last section of this chapter, we will prove Gauss's Theorem on Primitive Roots. We in fact show that there exists an integer
between 2 to $p - 1$ which is a primitive root of $p^n$ for all $n \geq 1$. 

\section{$a^x \eq b^y \, (mod \, p^n)$} \label{A}

\begin{definition} $\Z_n ^*$ is the set of all invertible elements of $\Z_n$, the ring of integers modulo $n$.
If we express elements of the ring $\Z_n$ as $0, 1, \cdots, n - 1$ then 
$\Z_n ^* = \{k \, \- \Z_n \, | \,  (k, n) = 1\}.$ \end{definition}

\begin{definition}: If $k$ is any element of $\Z_n ^*$ then $o(k) =$ order of $k$ is the smallest positive integer such that 
$k^{o(k)} \eq 1 \, (mod \, n)$.\end{definition}

\begin{proposition}
$\Z_n ^*$ is a group under multiplication. 
\end{proposition}
\begin{proof}
 If $k \- \Z_n ^*$ then $(k, n) = 1$. Therefore there exists integers $m, l$ such that $km + n l = 1$.
This implies $km \eq 1 \, (mod \, n)$. Also $(m, n) = 1$ because if $d = (m, n)$ then $d$ divides $m$ and $n$ implies $d$ divides $km + nl$.
So, $d \mid 1$. Therefore $d = 1$. Let $t \- \{1, 2, \cdots, n - 1\}$ be such that $t \eq m \, (mod \, n)$. Then, $t \- \Z_n ^*$ and $t$ is the inverse of $k$ in $\Z_n ^*$.
\end{proof}

 Let $a$ and $b$ be non-zero integers both not equal to 1. Let $p$ be a prime and $n \in \N.$ We look for solutions of above congruence. 
Two elementary cases are as follows:\\
\underline{Case 1. \, $p \mid a$, $p\mid b$}\\
$\Rightarrow a^x \eq 0 \, (mod \, p^n)$ for some $x$ and $ b^y \eq 0 \, (mod \, p^n)$ for some $y$\\
Hence $a^x \eq b^y (mod \, p^n)$ has non-zero solution in $(x, y)$.\\
\underline{Case 2. \, $p \mid a$ but $p \nmid b$}\\
$p \mid a \Rightarrow p \mid a^x \, ~~ \forall x > 0$ and $p \nmid b \Rightarrow p \nmid b^y \, ~~ \forall y > 0$\\
Therefore, for every $x, y > 0$, $p \nmid (a^x - b^y).$\\
Because, if $p \mid (a^x - b^y) \Rightarrow a^x - b^y = pk$ for some $k \in \Z$\\
~~~~~~~~~~~~~$ \Rightarrow a^x - pk = b^y$\\
~~~~~~~~~~~~~$ \Rightarrow p \mid b^y$.
~~~~~~~~~~~~ a contradiction.\\
Therefore, $a^x \eq b^y \, (mod \, p)$ has no non-zero solution in $(x, y)$.\\
\underline{Case 3. \, $p \nmid a$ and $p \nmid b$}\\
Therefore, $(a, p^n) = 1$ and $(b, p^n) = 1$. Consider, $a \, mod \,(p^n)$ and $b \, mod \, (p^n)$. Then $a$ and $b$ can be
considered as elements of $\Z_{p^n}^*$. Therefore, $a^{o(a)} \eq 1 \, (mod \, p^n)$ and $b^{o(b)} \eq 1 \, (mod \, p^n)$.\\
From this, $a^x \eq b^y \, (mod \, p^n)$ has non-zero solution $x = o(a)$, $y = o(b).$\\
We notice that, such solution will exists for every pair $(a, b)$ such that $p \nmid a$ and $p \nmid b.$ Hence we say that, 
$(x = o(a)$, $y = o(b))$ or $(x = o(a) k$, $y = o(b) l)$, $k, l \in \Z$ are \textbf{trivial} solutions of $a^x \eq b^y \, (mod \, p^n).$
(Here if $x$ is negative then $a^x$ is the inverse of $a^{-x}$ in the group $\n.$)\\

In other words, if $< a >$, $< b >$ denote the subgroups of $\n$ then we have, \textbf{$a^x \eq b^y \, (mod \, p^n)$ has non-trivial solution if and only if}
 \hbox{$<a> \cap <b> \, \neq \{1\}.$}\\

If $<a> \cap <b> \, = \, <c>$ then $c = a^r$ and $c = b^s$ for some $r, s$.
Therefore, $a^r \eq b^s \, (mod \, p^n) \Rightarrow (x = r, y = s)$ is non-zero solution of  $a^x \eq b^y \, (mod \, p^n).$
Actually, every other solution is $x = rk, y = sk ; \, \, k \in \Z.$\\

\underline{Summary of 3 cases:} We see that the first case is trivial, in second case 
there is no non-trivial solution and from the third case we can assume that, $(a, p) = (b, p) = 1$ and
$<a> \, = \, <b> \, \subseteq \, \Z_{p^n}^*.$\\

\framebox(225,25)[c]{$a^x \eq b^y \, (mod \, p^n); \, (a, p) = (b, p) = 1$}\\
From the above conditions on $a$ and $b$, we get $a$, $b \in \Z_{p^n}^*.$ We investigate the solutions $(x_n, y_n)$ of the congruence $a^{x_n} \eq b^{y_n} \, (mod \, p^n)$
 as $n \rightarrow \infty.$\\
\indent Let us further assume that, $<a> = <b>$ in $\Z_{p^n}^*.$ Then $o(a) = o(b).$
Therefore for every $x \, \in \{1, 2, \cdots, o(a)\}$, there exists unique $y \, \in \{1, 2, \cdots, o(b)\}$ such that 
$a^x \eq b^y \, (mod \, p^n).$ So this gives us all solutions of $a^x \eq b^y \, (mod \, p^n)$ and also the distribution of the solution as:\\
\noindent Fix $y = 1.$ Then there exists unique $x \, \in \{1, 2, \cdots, o(a)\}$ such that, 
\begin{eqnarray*}
 a^x &\eq& b \, (mod \, p^n)\\
 a^{2x} &\eq& b^2 \, (mod \, p^n)\\
&\vdots&\\
a^{x^2} &\eq& b^x \, (mod \, p^n)\\
&\vdots&
\end{eqnarray*}
Since, $<a> = <b>$, we must have that, $x$ is coprime to $o(a) = o(b)$. From above equations,
we get $\sigma : 1 \mapsto x \mapsto x^2 \mapsto \cdots$, $x^k \, (mod \, p^n)$ is a permutation of the set 
$\{1, 2, \cdots, o(a)=o(b)\}$. Hence $\sigma$ is an automorphism of $\Z_{o(a)}$. We shall see that for some special choices 
of $a$ and $b$ this permutation represents interesting patterns.  $Aut(\Z_{o(a)})$ is a group under composition.
In general, the group $Aut(\Z_n)$ is isomorphic to the multiplicative group $\Z_n ^*$. 

\section{Properties of $\Z_n ^*$} \label{B}
Recall the definition of $\Z_n ^*$ from the previous section. In this section we will prove the structure theorem for the group $\Z_n ^*$. 
This material can be found in any book on group theory, as a precise reference, we suggest \cite{DF}. 
We will prove the existence of the primitive root of prime by two ways: field theoretic approach and number theoretic approach. 
\begin{enumerate}
\item  $\Z_n ^*$ is an abelian group under multiplication with identity equal to 1.
\begin{definition} Euler totient function of $n$, $\phi(n)$ is defined as number of positive integers coprime to $n$ and less equal $n$. \end{definition}
\item order of $\Z_n ^*$ is $\phi (n)$
\item $\Z_n ^* \cong Aut \,(\Z_n , +).$\\
Proof: $k \longmapsto (1 \mapsto k)$ is an isomorphism.
\end{enumerate}

\begin{theorem}\label{5} Structure of the group $(\Z_n ^*, *)$:\\
if $n = p_1 ^{k_1} p_2 ^{k_2} \cdots p_r ^{k_r}$ is a prime factorisation of $n$ then,
$$ {\Z_n}^* \cong \Z_{{p_1}^{k_1}}^* \times \Z_{{p_2}^{k_2}}^* \times \dots \times \Z_{{p_r}^{k_r}}^*$$
\end{theorem}

\begin{proof}
 Since both sides are of same order, it is enough to prove that, there is an injective homomorphism for such a homomorphism would be surjective and would be an isomorphism.
Define a map $$f : {\Z_n}^* \longrightarrow \Z_{{p_1}^{k_1}}^* \times \dots \times \Z_{{p_r}^{k_r}}^*$$  
$$k \, (mod \, n) \mapsto (k \, (mod \, p_1 ^{k_1}), \cdots, k \, (mod \, p_r ^{k_r}))$$\\
We show that, $f$ is 1-1. Suppose $f(k) = f(l) = (a_1, \cdots, a_r)$. Then $k \eq l \, (mod \, p_i ^{k_i})$ for $i = 1, 2, \cdots, r$.
i.e. $p_i ^{k_i}$ divides $k - l$ for all $i$. Since $(p_i ^{k_i}, p_j ^{k_j}) = 1$ for $i \neq j$, we get product of  $p_i ^{k_i}$ divides $k - l$ i.e. 
$n$ divides $k - l$. Therefore $k \, (mod \, n) = l \, (mod \, n)$.\\
Therefore $f$ is 1-1. For any natural number $n$, we have $k \, (mod \, n) \times l \, (mod \, n) = k l \, (mod \, n)$. Therefore $f$ is also a homomorphism.
Hence $f$ is an isomorphism.
\end{proof}

\begin{theorem} \label{16}
 $\Z_{2^n}^* \cong \Z_2 \times \Z_{2^{n-2}}$ for $n \geq 2.$       
\end{theorem}

\begin{proof}
\begin{enumerate}
\item $(5, 2^n) = 1 \Rightarrow {5} \in \Z_{2^n} ^*.$ We prove (at the end) using binomial theorem that, order of $5$ is $2^{n-2}.$ 
\item $(2^n-1, 2^n) = 1 \Rightarrow {2^n-1} \in \Z_{2^n} ^*.$\\
$-1 \eq 2^n-1 \, (mod \, 2^n)$. $-1$ is of order 2. 
\item $-1 \notin \, <5>$\\
 Suppose, $-1 \in <5>$ then $-1 = 5^k$ for some $k$.\\
 $5 \eq 1 \, (mod \, 4) \Rightarrow 5^k \eq 1 \, (mod \, 4)$\\
 but $-1 \neq 1 \, (mod \,4)$. Therefore, $-1 \notin \, <5>$.
\item Let $H = <5>$, $K = < -1 > \, \subseteq \Z_{2^n} ^*$.  $H \cap K \, = \{1\}$.\\
Then $HK$ is a subgroup of $\Z_{2^n} ^*$ of order \\
$ \dfrac{|H| \, |K|}{|H \cap K|} = \dfrac{2^{n-2} \, 2}{1} = 2^{n - 1}$.\\
Also the order of $\Z_{2^n} ^*$ is $2^{n - 1}$. Therefore, $HK = \Z_{2^n} ^*$
\item $H \cap K = \{1\} \Rightarrow HK \cong H \times K.$\\
i.e. $\Z_{2^n} ^* \cong H \times K $ for $n \geq 2.$ \\
So $\t$ is an internal direct product of \hbox{$< 5 >$ and $< -1 >.$}
\end{enumerate}

It remains to prove that, order of ${5}$ in $\Z_{2^n} ^*$ is $2^{n-2}$, for $n \geq 3.$\\
For this we prove that, $ 5^{2^{n-2}} \eq 1 \, (mod \, 2^n)$ and $ 5^{2^{n-3}}  \neq 1 \, (mod \, 2^n) $.
\begin{eqnarray*}
    5^{\w} &=& (1 + 2^2)^{\w}\\
           &=& 1 + \w \, 2^2 + \dfrac{\w (\w -1)}{2!} (2^2)^2 + \cdots \\
           &\vspace{.2cm}&\\
           &+& \dfrac{\w (\w -1) \cdots (\w - k + 1)}{k!} (2^2)^k \\
           &\vspace{.2cm}&\\
           &+& \cdots + \w (2^2)^{(\w - 1)} + (2^2)^{(\w)}\\
           &\vspace{.2cm}&\\
           &=& 1+ 2^n + 2^n \dl{\dfrac{2^2(\w - 1)}{2!}}\dr + \cdots \\
           &\vspace{.2cm}&\\
           &+& 2^n \dl{\dfrac{2^{2k - 2} (\w - 1) \cdots (\w -k+1)}{k!}}\dr\\
           &\vspace{.2cm}&\\
           &+& \cdots + 2^n \, 2^{({2^{n-1} - 4})} + 2^n \, 2^{(2^{n-1} - n)}\\
           &\vspace{.2cm}&\\
           &\eq& 1 \, (mod \, 2^n) 
\end{eqnarray*}
Therefore order of 5 divides $\w$.

\begin{eqnarray*}
     5^{\v} &=& (1 + 2^2)^{\v}\\
            &\vspace{.2cm}&\\   
            &=& 1 + \v 2^2 + \dfrac{\v (\v -1)}{2!} (2^2)^2 + \cdots + \dfrac{\v (\v -1) \cdots (\v - k + 1)}{k!} (2^2)^k\\ 
            &\vspace{.2cm}&\\
            &+& \cdots + \v (2^2)^{(\v - 1)} + (2^2)^{(\v)}\\
            &\vspace{.2cm}&\\
            &\eq& 1 + 2^{n-1} \, (mod \, 2^n)\\
            &\vspace{.2cm}&\\
	    &\neq& 1 \, (mod \, 2^n)
\end{eqnarray*}
Therefore order of 5 is $\w.$
\end{proof}

\begin{definition}
 Let $n > 1$ be any natural number. We say that $r$ is a primitive root of $n$ if $\phi(n)$ is the smallest positive integer such that $$r^{\phi(n)} \eq 1 \, (mod \, n)$$
\end{definition}
In the group theoretic language, $r$ is a primitive root of $n$ is the same as the multiplicative group $\Z_n ^*$ is cyclic and $r$ is a generator of $\Z_n ^*$.
A major result of {\it{Gauss}} and {\it{Euler}} is the determination of $n$'s when the primitive root mod $n$ exists. There is an interesting issue about the existence
and construction of the primitive roots. 

\begin{theorem} \label{11}
 Every prime has a primitive root.
\end{theorem}

\begin{proof}
 A field theoretic argument for {\it{existence}} of primitive root mod $p$ is as follows:\\
We have to prove that the multiplicative group of the field $\Z_p$, namely $\Z_p ^*$ is cyclic.
We have already noted that order of $\Z_p ^*$ is $\phi(p) = p - 1$. So we prove that $\Z_p ^* \cong \Z_{p - 1}.$
Let us denote $\Z_p$ by $\mathbb{F}$ and $\Z_p ^*$ by $H$. As $H$ is of order $p - 1$, order of an element $t$ of $H$ 
divides $p - 1$. So $t$ is a ${p - 1}^{th}$ root of unity, a root of the polynomial $x^{p - 1} - 1.$
This polynomial has atmost $p - 1$ roots in any field. It follows that $H$ is the set of all ${p - 1}^{th}$ roots of unity in $\mathbb{F}$.\\
Now we prove that $H$ is cyclic. We use the structure theorem for abelian groups, which tells us that $H$ is isomorphic to a 
direct product of cyclic groups $$ H \cong \Z_{d_1} \times \cdots \times \Z_{d_k}$$ where $d_1 \, \mid \, d_2 \cdots \mid \, d_k$
and $p - 1 = d_1 d_2 \cdots d_k.$ The order of any element of $H$ divides $d_k$ because $d_k$ is the common multiple of all $d_i$'s.
So every element of $H$ is a root of the polynomial  $$x^{d_k} - 1.$$ This polynomial has atmost $d_k$ roots in $\mathbb{F}.$
But $H$ contains $p - 1$ elements and $p - 1 = d_1 d_2 \cdots d_k.$ The only possibility is that $d_k = p - 1$ and $k = 1$.
So $H$ is cyclic and isomorphic to $\Z_{p - 1}.$ 
\end{proof}

\begin{remark}
Gauss gave proof of existence of primitive root of $p$, without explicitly using field theory. We mention it below:
\end{remark}

\begin{proof}
 1 is the primitive root of 2.\\
Step 1: Let $p$ be any odd prime. Select a number $a$ which is coprime to $p$. Suppose order of $a \, (mod \, p)$ is $t$. Calculate the least residues of 
$a, a^2, a^3, \cdots, a^t$ modulo $p$. If $t = p - 1$ then $a$ is a primitive root of $p$.\\
Step 2: If $t < p - 1$ then choose another number $b$ such that $i) \, b$ is coprime to $p$, $ii) \, b$ is not congruent to any power of $a$, modulo $p$.   
Suppose order of $b \, (mod \, p)$ is $u$. Then $u$ cannot be equal to $t$ or $u$ cannot divide $t$. Because in any case, $b^t \eq 1 \, (mod \, p)$,
which cannot be true by condition $ii)$ on $b$. Now if $u = p - 1$ then $b$ is a primitive root of $p$. However if $u \neq p - 1$ but $u$ is a multiple of $t$
then we have obtained a number having order higher than $t$, hence we will be closer to our goal. 
We can find a number whose order modulo $p$ is equal to the least common multiple of $t, u$ say, $y$. Write $y = m n$ such that $(m, n) = 1$ and $m \mid t$, $n \mid u$.  
Suppose $a^{t/m} \eq A $ and $b^{u/n} \eq B $ modulo $p$. Then clearly, modulo $p$, order of $A$ is $m$ and order of $B$ is $n$.
Since $(m, n) = 1$, we must have order of $AB \, (mod \, p)$ is $mn = y$.\\
Step 3: If $y = p - 1$ then $AB$ is a primitive root of $p$. If $y \neq p - 1$ then find another number $c$ that does not appear in the least residues of
$AB, (AB)^2, \cdots, (AB)^y$  modulo $p$. Then either $c$ is a primitive root of $p$ or we get a number having order strictly greater than $y$. \\
Since the numbers we get by repeating this operation gives strictly increasing orders modulo $p$. We must finally get a number having maximum order.
This will be a primitive root. 
\end{proof}

\begin{remark}
 Above proof is ``non-constructive'' in the sense that there is no easy expression for the primitive root of $p$ in terms of $p$.
{\textit{Gauss}} has made the remark for the above theorem in his {\it{Disquisitiones Arithmeticae}} \cite{Gau} as: ``Finding primitive roots is reduced for the most part
to trial and error. {\textit{Euler}} admits that it is extremely difficult to pick out these numbers (primitive roots of prime $p$) and that their nature is one of the
deepest mysteries of numbers.''\\
We have taken above proof from \cite{Gau}, page no. 48-49.
\end{remark}

\begin{theorem} \label{12} $\Z_{p^n}^* \cong \Z_{p-1} \times \Z_{p^{n-1}}$ where $p$ is an odd prime and $n \geq 1.$ \end{theorem}

\begin{proof}
 \begin{enumerate}
\item For $n = 1$, from the above theorem, there exists a primitive root of $p$. 
\hbox{Therefore} $\Z_p ^*$ is isomorphic to $\Z_{p - 1}$.
\item For $n \geq 2$, we observe from the R.H.S. of the above isomorphism that, $(p-1, p^{n-1}) = 1 \Rightarrow \Z_{p-1} \times \Z_{p^{n-1}}$ is a cyclic group of order 
 $p^{n-1} (p-1).$ Since there is unique cyclic group of given order, it is enough to prove that, $\n$ is a cyclic group. 
\item We have already noted that, $\n$ is an abelian group of order $\phi(p^n) = p^{n-1} (p-1).$
\item Therefore, Sylow p-subgroup of $\n$ is of order $p^{n-1}.$\\
$(1+p, p) = 1 \, \Rightarrow (1+p, p^n) = 1$\\
$\Rightarrow 1 + p \in \n$. $1+p$ is of order $p^{n-1}$ in $\n.$ (Proof by binomial theorem, is given below.)
Hence $<1 + p>$ is a Sylow-p subgroup of $\n$. As $\n$ is an abelian group, this is the unique Sylow-p subgroup of $\n$
\item The map 
$$ \phi : \Z_{p^n} \longrightarrow \Z_p$$ 
$$ \bar{a} \, (mod \, p^n) \rightarrow \bar{a} \, (mod \, p)$$
is a surjective ring homomorphism. Because if $\bar{a} \- \Z_p$ then $a = r + p k$ for some $r \- \{0, 1, \cdots, p - 1\}$ and $k \- \Z$.
For $b = r + p^n k$, we get $\phi (\bar{b}) = \bar{a}$.\\
 Hence, this map restricted to group of units \\
$\n \xrightarrow{\phi} \Z_p ^*$ is also surjective. Therefore,
$$ \dfrac{\n}{ker \, \phi} \cong \Z_p ^*$$
$\Rightarrow $ order of $ker \, \phi = p^{n-1}.$
\item If $q \neq p$ is any prime dividing $\phi (p^n) = p^{n - 1} (p - 1)$ then $q \mid p-1$.
Therefore, the Sylow-$q$ subgroup of $\n$ can not be contained in the $ker \, \phi$. Hence, Sylow-$q$ subgroup of $\n$ is isomorphic to Sylow-$q$ subgroup of $\Z_p ^* \cong \Z_{p-1}.$
Therefore, Sylow-$q$ subgroup must be cyclic.
\item $\n$ is abelian $\Rightarrow$ it is product of it's Sylow subgroups. \\
Thus $\n$ is an internal direct product of it's $q$-Sylow subgroups, $q \neq p, \, q \,  \mid \, p -1$ and $< 1 + p>.$  
\end{enumerate}
It remains to prove that order of $1+p$ is $p^{n-1}.$
\begin{eqnarray*}
  (1+p)^{\p} &=& 1+ \p \, p + \dfrac{\p \,(\p - 1)}{2!} p^2 + \dfrac{\p \, (\p - 1) \, (\p - 2)}{3!} p^3\\
             &\vspace{.2cm}&\\
             &+&\cdots + \dfrac{\p \, (\p - 1) \cdots (\p - k + 1)}{k!} p^k\\
             &\vspace{.2cm}&\\
             &+&\cdots + \p \, p^{(\p - 1)} + p^{\p}\\
             &\vspace{.2cm}&\\
             &=& 1+ p^n + p^n \dl{\dfrac{p \, (\p - 1)}{2!}}\dr\\
             &\vspace{.2cm}&\\
             &+& p^n \, \dl{\dfrac{p^2 \, (\p - 1) \, (\p - 2)}{3!}}\dr\\
             &\vspace{.2cm}&\\
             &+&\cdots + p^n \, \dl{\dfrac{p^{k - 1} \, (\p - 1) \cdots (\p - k + 1)}{k!}}\dr\\
             &\vspace{.2cm}&\\
             &+&\cdots + p^n \, p^{(\p - 2)} + p^n \, p^{(\p - n)}\\
             &\eq& 1 \, (mod \, p^n)       
\end{eqnarray*}
Therefore, order of $1+p$ divides $\p.$\\
Now,
\begin{eqnarray*}
 (1+p)^{\q} &=& 1+ \q \, p + \dfrac{\q \,(\q - 1)}{2!} p^2 + \dfrac{\q \, (\q - 1) \, (\q - 2)}{3!} p^3\\
             &\vspace{.2cm}&\\
             &+&\cdots + \dfrac{\q \, (\q - 1) \cdots (\q - k + 1)}{k!} p^k\\
             &\vspace{.2cm}&\\
             &+&\cdots + \q \, p^{(\q - 1)} + p^{\q}\\
             &\vspace{.1cm}&\\
             &\eq& 1 + \p \, (mod \, p^n)\\
             &\vspace{.1cm}&\\
             &\neq& 1 \, (mod \, p^n)
\end{eqnarray*}
Therefore, order of $1 + p = \p. $
\end{proof}

\textbf{Observations related to above theorems:} \begin{enumerate}
\item If $p$ is an odd prime then $\n$ is cyclic. 
\item $\t$ is not cyclic.
\begin{enumerate}
\item Number of subgroups of order 2 is 3 if $n = 3.$
More explicitly, 
$$\Z_{2^3}^* = \Z_8 ^* = \{\bar{1}, \bar{3}, \bar{5}, \bar{7}\} \cong \Z_2 \times \Z_2 .$$
$$\{\bar{1}, \bar{3}\}, \, \{\bar{1}, \bar{5}\}, \, \{\bar{1}, \bar{7}\}$$ are 3 subgroups of order 2.
\item For $n \geq 4,$ number of subgroups of order $2^k$ is 2 if $1 \leq k \leq n-2.$ 
  $$\Z_{2^n}^* \cong \Z_2 \times \Z_{2^{n-2}}$$
If $\bar{x} \in \Z_{2^{n-2}}$ is of order $2^k, \, \, 2 \leq k \leq n-2 $ then $<(\bar{0}, \bar{x})>$ and $<(\bar{1}, \bar{x})>$
are the only subgroups of order $2^k$. We observe that, $<(\bar{0}, \bar{x})>$ is the unique subgroup of order $2^k$,
consisting of only squares.
\end{enumerate}
 \end{enumerate}

\section{Primitive roots of $p^n$} \label{G}

The aim of this section is to prove that, there exists a number between 2 and $p - 1$ which works as a primitive root of $p^n$ for all $n \geq 1.$
At the beginning we prove some essential results which will be required throughout the thesis.
We have given group theoretic proofs of some of these results in the previous section. Here we give elementary number theoretic proofs.
\begin{theorem} \textbf{Fermat's Theorem}
Let $p$ be a prime and $a$ be an integer such that $p \nmid a$ then $a^{p - 1} \eq 1 \, (mod \, p)$
\end{theorem}

\begin{proof}
In the group theoretic language, the statement of this theorem says that, if $a \- \Z_p ^*$ then $a^{|\Z_p ^*|} = 1.$ 
In any finite group $a^{|G|} = 1$ for all $a \- G$. A group theoretic proof of this theorem is clear. But we give here more elementary proof, using number theory arguments.\\
 Let $p = 2$. The condition $ 2 \nmid a$ implies that $a$ must be an odd integer. Then $a \eq 1 \, (mod \, 2)$.\\
Let $p$ be an odd prime. Then $0, 1, 2, \cdots, p- 1$ is the least residue system modulo $p$. Since $a$ is coprime to $p$,
we have $0, a, 2a, \cdots, a (p - 1)$ is another residue system modulo $p$. Therefore, for every $r \- \{1, 2, \cdots, p - 1\}$ there exists unique 
$s \- \{1, 2, \cdots, p - 1\}$ such that $ar \eq s \, (mod \, p)$. In this way, we get $p - 1$ number of congruences, one for each value of $r$.
Multiplying all the $p - 1$ congruences, we get $$a.2a.3a\cdots a(p - 1) \eq 1.2.3 \cdots (p - 1) \, (mod \, p)$$ Therefore we get, 
$ a^{p - 1} (p - 1)! \eq (p - 1)! $. Since $p$ does not divide $(p - 1)!$ \\ we must have, $a^{p - 1} \eq 1 \, (mod \, p)$.
\end{proof}

\begin{theorem} \textbf{Euler's Theorem}
 If $n \geq 1$ and $a$ is an integer such that $gcd (a, n) = 1$ then $a^{\phi(n)} \eq 1 \, (mod \, n)$
\end{theorem}
\begin{proof}
Special Case: Suppose $n$ is a power of a prime $p$ i.e. $n = p^k, \, k > 0$.\\
Proof by induction on $k$. If $k = 1$ then the assertion reduces to Fermat's theorem. \\
Induction hypothesis: For a fixed value of $k >0$, assume that,  $$a^{\phi(p^k)} \eq 1 \, (mod \, p^k)$$
We have to prove that,   $a^{\phi(p^{k + 1})} \eq 1 \, (mod \, p^{k + 1})$.\\
From the induction step, $a^{\phi(p^k)} = 1 + p^k u$ for some integer $u$.\\
Now, $\phi(p^{k + 1}) = p^{k + 1} - p^k = p \, (p^k - p^{k - 1}) = p \, \phi(p^k)$.\\
We use binomial theorem to prove the claim.
\begin{eqnarray*}
 a^{\phi(p^k)} &=& 1 + p^k u\\
 a^{p \, \phi(p^k)} &=& (1 + p^k u)^p \\
                 &=& 1 + p \, p^k u + \dfrac{p (p - 1)}{2} p^{2k} u^2 + \cdots + p^{pk} u^p\\
                 &=& 1 + p^{ k + 1} u + p^{k + 1} \dfrac{p^k (p - 1)}{2} u^2 + \cdots + p^{pk} u^p\\
                 &\eq& 1 \, (mod \, p^{k + 1})
\end{eqnarray*}

General case: Let $n = p_1 ^{k_1} p_2 ^{k_2} \cdots p_r ^{k_r}$ be the prime power factorisation of $n$.
Then, $\phi(n) = \phi(p_1 ^{k_1}) \, \phi(p_2 ^{k_2}) \cdots \phi(_r ^{k_r}).$\\
From the special case, $$a^{\phi(p_i ^{k_i})} \eq 1 \, (mod \, p_i ^{k_i}) \h3 i = 1, 2, \cdots, r$$
Therefore, $$a^{\frac{\phi(n)}{{\phi(p_i ^{k_i})}} {\phi(p_i ^{k_i})} } \eq 1 \, (mod \, p_i ^{k_i}) \h3 i = 1, 2, \cdots, r$$
$$ a^{\phi(n)} \eq 1 \, (mod \, p_i ^{k_i}) \h3 i = 1, 2, \cdots, r$$
Since all $p_i$'s are distinct, we get, $a^{\phi(n)} \eq 1 \, (mod \, n)$.
\end{proof}

Following result by Gauss gives all the natural numbers having primitive roots.
Proof of the theorem  \ref{15} is clear from the structure theorems \ref{5}, \ref{16} and \ref{12}, but we present here more elementary number theoretic proof.

\begin{theorem} \label{15}  \textbf{Gauss's Theorem on Primitive roots:}\\ 
A natural number $n$ has a primitive root if and only if $n = 2, 4, \,p^m, 2 p^m$ for $p$ odd prime and $m \geq 1$.
\end{theorem}
\begin{proof}  
If $n =  2, 4, p^m, 2 p^m$ for $p$ odd prime and $m \geq 1$ then $n$ has a primitive root is clear from the theorem \ref{5} and theorem \ref{12}.
Theorem \ref{11} gives the existence of primitive root of $p$. \\
Conversely, if $n$ is not 2 or 4 or of the form $p^m$ or $2 p^m$ then $n$ must have one of the following form. We prove that in each case $\Z_n ^*$ contains a non-cyclic group.
Therefore $\Z_n ^*$ cannot be cyclic and $n$ cannot have a primitive root. \\
\noindent \underline{Step 1:} For $n \geq 3$, the integer $2^n$ has no primitive root. \\
Because, for $n \geq 3$, we have $\t \cong \Z_2 \times \Z_{2^{n - 2}}$, so $\Z_2 \times \Z_2$ is a non-cyclic subgroup of $\t$. Therefore $\t$ is also not-cyclic.
Hence $2^n$ has no primitive root for $n \geq 3.$\\
\underline{Step 2:} The integer $n$ fails to have a primitive root if either \\
i) $n$ is divisible by two odd primes, or\\
ii) $n$ is of the form $2^m p^k$ for $m \geq 2$.\\ 
{\it{Proof.}} i) In this case also $\Z_2 \times \Z_2$ is a non-cyclic subgroup of $\Z_n ^*.$ Therefore $n$ has no primitive root. \\
ii) If $n = 2^m p^k$ for $m \geq 2.$ Then by the structure theorem \ref{5} we get, 
\begin{eqnarray*}
\Z_n ^* &\cong& \Z_{2^m} ^* \times \Z_{p^k} ^*\\
        &\cong& \Z_2 \times \Z_{2^{m - 2}} \times \Z_{p - 1} \times \Z_{p^{k - 1}}
\end{eqnarray*}
Since $p$ is an odd prime, $p - 1$ is even, and we get $\Z_2 \times \Z_{p - 1}$ is a non-cyclic subgroup of $\Z_n ^*$. Therefore $\Z_n ^*$ is not-cyclic.
Hence the proof is complete.
\end{proof}

 We extend the theorem \ref{15} partially as follows. We will be able to prove it at the end of this section.
\begin{theorem} \label{17}
Let $p$ be an odd prime. There exists $r \- \{2, 3, \cdots, p - 1\}$ which is a primitive root of $p^n$ for all $n \geq 1.$
\end{theorem}
The case $n = 1$ is already done. For $n \geq 2$ we need to prove some more results.

\begin{note} If $r$ is a primitive root of $p$ then $r$ belongs to non-square elements modulo $p$.
If $r \eq s^2 \, (mod \, p)$, then $$(s^2)^{\frac{p - 1}{2}} = s^{p - 1} \eq 1 \, (mod \, p)$$
Therefore, order of $r$ modulo $p$ is strictly smaller than $p - 1$. Hence, $r$ can not be square modulo $p$.
\end{note}

\begin{lemma}
 $p \eq 1 \, (mod \, 4)$ \iff $-1$ is a square $mod \,p.$ 
\end{lemma}
\begin{proof}
Suppose $p \eq 1 \, (mod \, 4)$. Let $r$ be a primitive root of $p$. Therefore,
\begin{eqnarray*}
 r^{p - 1} \eq 1 \, (mod \, p)\\
 {\dl{r^{\frac{p - 1}{2}}}\dr}^2 \eq 1 \, (mod \, p)\\
 r^{\frac{p - 1}{2}} \eq 1 \, or \, -1 \, (mod \, p)
\end{eqnarray*}
 But $r$ is a primitive root of $p$, therefore $r^{\frac{p - 1}{2}}$ must be congruent to $-1 \, (mod \, p)$.
 Since $\dfrac{p - 1}{2}$ is even, we can write, 
 \begin{eqnarray*}
  {\dl{r^{\frac{p - 1}{4}}}\dr}^2 \eq -1 \, (mod \, p)
 \end{eqnarray*}
 Conversely, Let $p$ be a prime and $x$ is such that
   \begin{eqnarray*}
  x^2 \eq -1 \, (mod \, p)\\
  (x^2)^{\frac{p - 1}{2}} \eq (-1)^{\frac{p - 1}{2}} \, (mod \, p)\\
  x^{p - 1} \eq (-1)^{\frac{p - 1}{2}} \, (mod \, p)
  \end{eqnarray*}
  By Fermat's theorem $x^{p - 1} \eq 1 \, (mod \, p)$. On the r.h.s.
  $(-1)^{\frac{p - 1}{2}} = 1$ \iff $\dfrac{p - 1}{2}$ is even i.e. $p \eq 1 \, (mod \, 4).$
 \end{proof}

 \begin{note}
  In the above proof, $r^{\frac{p - 1}{4}}$ is a square root of $-1 \, (mod \, p)$, which is given in terms of primitive root of $p$.
   It is possible to find the solution of $x^2 \eq -1 \, (mod \, p)$ only in terms of $p$, namely, $x = \pm \dl{\dfrac{p - 1}{2}}\dr !$
 But when $p$ is very large prime, the number $\dl{\dfrac{p - 1}{2}}\dr !$ is also very large, hence this solution is not useful. 
 \end{note}

\begin{lemma} Let $r$ be a primitive root of $p.$\\
 1. If $p \eq 1 \, (mod \, 4)$ then $-r \, (mod \, p)$ is a primitive root of $p.$\\
 2. If $p \eq 3 \, (mod \, 4)$ then $- r^2 \, (mod \, p)$ is a primitive root of $p.$
\end{lemma}
\begin{proof}
 1. $r$ is a primitive root of $p$, therefore from the note 1.3.1,  $r$ is non-square modulo $p$. Also $p \eq 1 \, (mod \, 4)$, therefore $-1$ is a square $mod \, p$.
 Hence $(-1) r$ is also a non-square $mod \, p$. Clearly, $(-r)^{p - 1} \eq 1 \, (mod \, p)$. We have to prove that $p - 1$ is the smallest integer
 satisfying $(-r)^{p - 1} \eq 1 \, (mod \, p)$. Suppose $m$ is the least positive integer such that $(-r)^m \eq 1 \, (mod \, p)$. Then $m$ divides $p - 1$.\\
 If $m$ is odd, 
\begin{eqnarray*}
(-r)^m &\eq& 1 \, (mod \, p)\\
-r^m   &\eq& 1 \, (mod \, p)\\
r^m    &\eq& -1 \, (mod \, p)\\
r^{2m} &\eq& 1 \, (mod \, p)
\end{eqnarray*}
But $r$ is a primitive root of $p$, implies that, $p - 1$ divides $2m$.  
\begin{eqnarray*}
 m \leq p - 1 \leq 2m \, &and& \, m \, | \, p - 1 , p - 1 \, | \, 2m\\
\Rightarrow m = p - 1 \, &or& \, m = \dfrac{p - 1}{2}
\end{eqnarray*}
This is contradiction, since  $p - 1$ is even and $p \eq 1 \, (mod \, 4) \Rightarrow \dfrac{p - 1}{2}$ is also even, but $m$ is an odd number.\\
Therefore $m$ must be even. Let $m = 2n$. Then, 
\begin{eqnarray*}
 (-r)^m &\eq& 1 \, (mod \, p)\\
 (-r)^{2n} &\eq& 1 \, (mod \, p)\\
r^{2n} &\eq& 1 \, (mod \, p)
\end{eqnarray*}
But $r$ is a primitive root of $p$, implies that, $p - 1$ divides $2n = m$.\\
We get, $m \, | \, p - 1$ and $p - 1 \, | \, m$ gives, $m = p - 1$.
Therefore, $-r$ is a primitive root of $p$.\\

2. Let $ p \eq 3 \, (mod \, 4)$ and $r$ is a primitive root of $p$. We have to prove that $-r^2 \, (mod \, p)$ is also a primitive root of $p$.
By Fermat's Theorem, $(-r^2)^{p - 1} \eq 1 \, (mod \, p)$. Let $m$ be the smallest positive integer such that,
\begin{eqnarray*}
(-r^2)^m &\eq&  1 \, (mod \, p) 
\end{eqnarray*}
Therefore, $m \, | \, p - 1$
\begin{eqnarray*}
(-1)^m r^{2m} &\eq& 1 \, (mod \, p)
\end{eqnarray*}
If $m$ is odd, then
\begin{eqnarray*}
 r^{2m} &\eq& -1 \, (mod \, p)\\
 r^{4m} &\eq&  1 \, (mod \, p)
\end{eqnarray*}
As $r$ is a primitive root of $p$, we get, $p - 1 \, | \, 4m $. Therefore, we get,
\begin{eqnarray*}
 m \, | \, p - 1 &and& p - 1 \, | \, 4m\\
\Rightarrow p - 1 = mk &and& 4m = (p - 1)l\\
\Rightarrow 4m = mkl\\
\Rightarrow  4 = kl\\
\Rightarrow \, either \, (k = 1, l = 4) \, or \,  (k = 4, l = 1) \, or \, (k = 2, l = 2)
\end{eqnarray*}
If $k = 1$ then $p - 1 = m$. But $m$ is odd and $p - 1$ is even. Hence contradiction. \\
If $l = 1$ then $p - 1 = 4m$. So, $p \eq 1 \, (mod \, 4)$. Hence contradiction. \\
If $k = l =2$ then $p - 1 = 2m $ so, $m = \dfrac{p - 1}{2}$. Since $p \eq 3 \, (mod \, 4)$, we get, $\dfrac{p - 1}{2}$ is an odd integer.
By the property of $m$, 
 \begin{eqnarray*}
  (-r^2)^m &\eq&  1 \, (mod \, p) \\
(-r^2)^{\frac{p - 1}{2}} &\eq&  1 \, (mod \, p) \\
 (-1)^{\frac{p - 1}{2}} r^{p - 1} &\eq&  1 \, (mod \, p) \\
- r^{p - 1} &\eq& 1 \, (mod \, p)\\
r^{p - 1} &\eq& -1 \, (mod \, p)
 \end{eqnarray*}
This is contradiction, because, $r^{p - 1} \eq 1 \, (mod \, p).$ Hence, $m$ must be an even number. Suppose $m = 2n$.
 \begin{eqnarray*}
  (-r^2)^m &\eq&  1 \, (mod \, p) \\
(-r^2)^{2n} &\eq&  1 \, (mod \, p) \\
r^{4n} &\eq&  1 \, (mod \, p) \\
r^{2m} &\eq&  1 \, (mod \, p) 
 \end{eqnarray*}
As $r$ is a primitive root of $p$, we get, $p - 1 \, | \, 2m $. We already have, $m$ divides $p - 1$.
 Therefore, we get,
 \begin{eqnarray*}
  m \leq p - 1 \leq 2m \, &and& \, m \, | \, p - 1 , p - 1 \, | \, 2m\\
\Rightarrow m = p - 1 \, &or& \, m = \dfrac{p - 1}{2}
 \end{eqnarray*}
But $m$ is even and $\dfrac{p - 1}{2}$ is odd. Hence $m$ must be equal to $p - 1$.
Therefore, $-r^2 \,(mod \, p)$ is a primitive root of $p$.
\end{proof}

\begin{lemma} \label{4}
If $r \- \{2, 3, \cdots, p - 1\}$ is a primitive root of $p$ then \\
1. for $p \eq 1 \, (mod \, 4)$ either $r$ or $- r \, (mod \, p)$ or both are primitive roots of $p^2$\\
2. for $p \eq 3 \, (mod \, 4)$ either $r$ or $- r^2 \, (mod \, p)$ or both are primitive roots of $p^2$.   
\end{lemma}
\begin{proof}
1. Let $p \eq 1 \, (mod \, 4)$. Suppose $r$ is a primitive root of $p$ but not of $p^2$.
Therefore $r^{p - 1} \eq 1 \, (mod \, p^2)$.
From the above theorem, $-r \, (mod \, p)$ is also a primitive root of $p$.\\
Claim: $-r \, (mod \, p)$ is a primitive root of $p^2$. i.e. $\phi (p^2) = p(p - 1)$ is the smallest positive integer such that, $(-r)^{p (p - 1)} \eq 1 \, (mod \, p^2)$.\\
For this, it is enough to show that, $(-r)^{p - 1} \neq 1 \, (mod \, p^2)$.\\
 Since $r^{p - 1} \eq 1 \, (mod \, p^2)$ we have $r^{p - 1} = 1 + p^2 k$ for some $k$.
\begin{eqnarray*}
 (-r)^{p - 1} &=& (p - r)^{p - 1}\\
              &=& (-r + p)^{p - 1}\\
              &=& (-r)^{p - 1} + (p - 1) (-r)^{p - 2} p + \dfrac{(p - 1)(p - 2)}{2} (-r)^{p - 3} p^2 + \cdots + p^{p - 1}\\
              &\eq& (-r)^{p - 1} + (p - 1) (-r)^{p - 2} p \, (mod \, p^2)\\
              &\eq& (-r)^{p - 1} - (-r)^{p - 2} \, (mod \, p^2)\\
              &\eq& r^{p - 1} + r^{p -2} \, (mod \, p^2)\\
              &\eq& 1 + r^{p - 2} \, (mod \, p^2)\\
              &\neq& 1 \, (mod \, p^2)
\end{eqnarray*}
Hence proved.\\

2. Let $p \eq 3 \, (mod \, 4)$. Suppose $r$ is a primitive root of $p$ but not of $p^2$.
Therefore $r^{p - 1} \eq 1 \, (mod \, p^2)$. From this we get $(-r^2)^{p - 1} \eq 1 \, (mod \, p^2)$
From the above theorem, $-r^2 \, (mod \, p)$ is also a primitive root of $p$.\\
Claim: $-r^2 \, (mod \, p)$ is a primitive root of $p^2$. i.e. $\phi (p^2) = p(p - 1)$ is the smallest positive integer such that, $(-r^2)^{p (p - 1)} \eq 1 \, (mod \, p^2)$.
Let $k$ be the smallest positive integer such that, $-r^2 + pk$ is a positive integer. Then $-r^2 + pk$ belongs to the set $\{2, 3, \cdots, p - 1\}$.
Since $-r^2 \eq -r^2 + pk \, (mod \, p),$ we prove that, $(-r^2 + pk)^{p - 1} \neq 1 \, (mod \, p^2)$.
\begin{eqnarray*}
 (-r^2 + pk)^{p - 1} &=& (-r^2)^{p -1} + (p - 1) (-r^2)^{p -2} pk + \dfrac{(p - 1)(p - 2)}{2} (-r^2)^{p - 3} (pk)^2 + \cdots + (pk)^{p - 1}\\
                     &\eq& (-r^2)^{p -1} + (p - 1) (-r^2)^{p -2} pk \, (mod \, p^2)\\
                     &\eq& (-r^2)^{p -1} - (-r^2)^{p -2} pk \, (mod \, p^2)\\
                     &\eq&  1 - (-r^2) \, pk \, (mod \, p^2)\\
                     &\neq& 1 \, (mod \, p^2) 
                     \end{eqnarray*}
Therefore, $- r^2 \, (mod \, p)$ is a primitive root of $p^2.$
\end{proof}

\begin{lemma} \label{14}
  If $r$ is a primitive root of both $p$ and $p^2$ then $r$ is a primitive root of $p^n$ for all $n \geq 3.$
\end{lemma}
\begin{proof}
$r$ is a primitive root of $p$, so $p - 1$ is the smallest positive integer such that $r^{p - 1} \eq 1 \, (mod \, p)$. Write $r^{p - 1} = 1 + pk$.
Since, $r$ is also a primitive root of $p^2$, so $p (p - 1)$ is the smallest positive integer such that $r^{p (p - 1)} \eq 1 \, (mod \, p^2)$.
Therefore, $r^{p -1} \neq 1 \, (mod \, p^2)$. This implies that, in the equation $r^{p - 1} = 1 + pk$, the integer $k$ is not divisible by $p$.
Now, $r^{p (p - 1)} \eq 1 \, (mod \, p^2)$, therefore $r^{p(p - 1)} = (r^{p - 1})^p = (1 + pk)^p$.
Let $m$ be the smallest positive integer such that, $r^m \eq 1 \, (mod \, p^3)$. As $r$ is coprime to $p$, by Fermat's theorem,
 $r^{\phi(p^3)} = r^{p^2 (p - 1)} \eq 1 \, (mod \, p^3)$.
Therefore, $m$ divides $p^2 (p - 1)$. Also, $r^m \eq 1 \, (mod \, p^3)$ implies that $r^m \eq 1 \, (mod \, p^2)$. But $r$ is a primitive root of $p^2$, so $p(p - 1)$ must divide $m$.
Therefore, either $m = p(p - 1)$ or $m = p^2 (p - 1)$. 
\begin{eqnarray*}
 r^{p (p - 1)} &=& (r^{p - 1})^p\\
               &=& (1 + pk)^p\\
               &=& 1 + p\, pk + \dfrac{p (p - 1) }{2} pk + \cdots + (pk)^p\\
               &=& 1 + p^2 k \,(mod \, p^3)
\end{eqnarray*}
Since $p$ does not divide $k$, $ r^{p (p - 1)} \neq 1 \,(mod \, p^3)$. Therefore $m = p^2 (p - 1)$.\\
Induction hypothesis: Suppose $r$ is a primitive root of $p, p^2, \cdots, p^n$.\\
$r$ is a primitive root of $p^{n - 1}$ implies that, $p^{n - 2} (p - 1)$ is the smallest number such that $r^{p^{n - 2} (p - 1)} \eq 1 \, (mod \, p^{n - 1})$.
Let $r^{p^{n - 2} (p - 1)} = 1 + p^{n - 1} k$. Since $r$ is also a primitive root of $p^n$, we must have, $p \nmid k$.
Now, $r$ is coprime to $p$, therefore $r^{\phi(p^{n + 1})} \eq 1 \, (mod \, p^{n + 1})$ i.e. $r^{p^n (p - 1)} \eq 1 \, (mod \, p^{n + 1})$. 
We have to prove that $p^n (p - 1)$ is the smallest number with this property. For this, it is enough to prove that $r^{p^{n - 1} (p - 1)} \neq 1 \, (mod \, p^{n + 1})$.
Now, $r^{p^{n - 2} (p - 1)} = 1 + p^{n - 1} k$ implies that $r^{p^{n - 1} (p - 1)} = (r^{p^{n - 2} (p - 1)})^p = (1 + p^{n - 1} k)^p \eq 1 + p^n k \,(mod \, p^{n + 1})$.
Since $p$ does not divide $k$, we get $1 + p^n k \neq 1 \, (p^{n + 1})$. Hence $p^n (p - 1)$ is the smallest positive integer such that $r^{p^n (p - 1)} \eq 1 \, (mod \, p^{n + 1})$
\end{proof}

\begin{proof}
of theorem \ref{17} is now clear from the lemma \ref{4} and lemma \ref{14}.
\end{proof}
Note that, from lemma \ref{4}, if $r$ is a primitive root of $p$ then either $r$ or $-r$ or $-r^2$ is primitive root of both $p$ and $p^2$ and hence of all $p^n$'s $n \geq 1$.
We already know from the Gauss's theorem that $r(1 + p)$ and $r + p$ both are primitive roots of $p^n$ for all $n$, but none of them modulo $p$ is primitive root for all $p^n$'s.
If $r$ is a primitive root of $p$ but not of $p^2$ then $r (1 + p)$ is a primitive root of $p^2.$
Because, $r$ is not a primitive root of $p^2$ implies that order of $r \, (mod \, p^2)$ is $p - 1$. By the proof of theorem \ref{12} 
order of $1 + p \, (mod \, p^2)$ is $p$. Since $p$ and $p - 1$ are coprime, order of $r (1 + p) \, (mod \, p^2)$ must be $p (p - 1)$.\\
Therefore by theorem \ref{14} we have, $r(1 + p)$ is a primitive root of $p^m$ for all $m \geq 1.$\\

\begin{example} 
\begin{enumerate}
\item (\textit{Gauss}) Let $p = 29$ and $a = 14$. Then, $28$ is the smallest positive integer satisfying 
$$ 14^{28} \eq 1 \, (mod \, 29)$$
$$ 14^{28} \eq 1 \, (mod \, 29^2)$$
Therefore, $14$ is a primitive root of $29$ but not of $29^2$. Since $29 \eq 1 \, (mod \, 4)$, we have $-14 (mod 29) = 15$ is a primitive root of $29^2$.
\item Let $p = 43$, $a = 19$ then $19$ is a primitive root of $43$ but not of $43^2.$ Since \hbox{$43 \eq 3\, (mod \, 4)$,} 
we calculate  $-19^2 \, (mod \, 43)$. It is equal to $26$. Then $26$ is a primitive root of both $43$ and $43^2$.
\end{enumerate}
\end{example}
\begin{example} \begin{enumerate}
\item For $p = 5$, primitive roots of $5$ are $2, 3.$  Squares modulo $5$ are $1, 4 = -1.$
\item For $p = 7$, primitive roots of $7$ are $3, 5.$ Squares modulo $7$ are $1, 2, 4.$
\end{enumerate}
\end{example}
\newpage
In the above table,
We have used `SAGE' to calculate primitive root of $p^n$ for all $n$ for first few primes.

\begin{table}
 \begin{tabular}[c]{|c|c||c|c|}
\hline
$p \eq 1\, (mod \, 4)$ & primitive root of $p^n$ & $p \eq 3\, (mod \, 4)$ & primitive root of $p^n$ \\
\hline
5 &2, 3& 7& 3, 5\\
\hline
13 & 2, 6 & 11& 2, 6, 7, 8\\
\hline
17 &3, 5, 6, 7 & 19 & 2,3, 10, 13, 14, 15\\
\hline
29 & 2, 3, 8, 10, 11, 15 & 23 & 5, 7, 10, 11, 14, 15, 17, 19, 20, 21\\
\hline
37 & 2, 5, 13, 15, 17, 19 & 31 & 3, 11, 12, 13, 17, 21, 22, 24\\
\hline
41 & 6, 7, 11, 12, 13, 15, 17, 19 & 43 & 3, 5, 12, 18, 20, 26, 28, 29, 30, 33, 34\\
\hline
 \end{tabular}
\end{table}

\chapter{ $n$-ADIC INTEGERS}

\textsc{THIS} chapter contains some information about $n$-adic \hbox{integers.} We will define $n$-adic \hbox{integers} by two ways: 1. projective limit approach and 2. metric approach,
that is using valuation on integers. In section \ref{C}, we will prove some algebraic properties of the ring of $n$-adic integers. In section \ref{D} we will give explicit
construction of Teichm\"{u}ller units in $p$-adic integers. We will also prove the structure theorem for the group of units of the ring $\zp$.\\

\textbf{ Two approaches to $n$-adic integers}\\
Given any natural number $x$ it can be written as $$ x_o + x_1 10 + x_2 10^2 + \cdots + x_k 10^k$$ where
$x_i \- \{0, 1, 2, 3, 4, 5, 6, 7, 8, 9\}$ and $k$ is a finite number.
Similarly, let $n > 1$ be a fixed natural number. An $n$-adic integer is a number written in base $n$.
i.e. write $x$ as $x_o + x_1 n + x_2 n^2 + \cdots + x_k n^k$, 
where $x_i \- \{0, 1, 2, \cdots, n - 1\}$. We say that, this is a $n$-adic expansion of $x$.
Difference between integers and $n$-adic integers is: the infinite sum $\displaystyle\sum_{i = 0} ^{\infty} x_i n^i$ does not make sense in integers 
but it does make sense in $n$-adic integers as defined below. There are two ways of defining $n$-adic integers. One is using projective limit of $\Z_{n^r}$ in
which we consider the structure of the ring $\Z_{n^r}$ simultaneously for all $r$. Other way is power series in $n$ with coefficients between 0 to $n - 1$, 
which is natural generalization of decimal expansion of natural numbers as explained above. 

\section{Projective limit approach to $n$-adic integers}

\indent In appendix we have given general definitions of a projective limits of sets, groups, rings, topological spaces. A reference for this is \cite{DF}.
Here we will consider only the projective limit of rings $\Z_{p^n}.$\\

\indent Let $n > 1$ be a fixed natural number having prime power factorisation $p_1 ^{k_1} p_2 ^{k_2} \cdots p_r ^{k_r}$. 
Then, we have following ring isomorphism
$$ \Z_n \cong \Z_{p_1 ^{k_1}} \times \Z_{p_2 ^{k_2}} \times \cdots \times \Z_{p_r ^{k_r}} \h3 (*)$$
Since primes appearing in prime power factorisation of $n$ and $n^t$ are same for any natural number $t$, we can write,
$$ \Z_{n^t} \cong \Z_{p_1 ^{k_1 t}} \times \Z_{p_2 ^{k_2 t}} \times \cdots \times \Z_{p_r ^{k_r t}} \h3 (**) $$

\textbf{projective limit approach:} Taking the projective limit on both sides of $(**)$ as $t \rightarrow \infty$, we get
$$ \Z_{(n)} \cong \Z_{(p_1 ^{k_1})} \times \Z_{(p_2 ^{k_2})} \times \cdots \times \Z_{(p_r ^{k_r})} $$
\begin{definition}: $\Z_{(p^k)}$ is projective limit of $\Z_{p^{kt}}$ as t goes to infinity.
$$ \Z_{(p^k)} = \{x = (\cdots, x_2, x_1) \- \prod_{t = 1} ^{\infty} \Z_{p^{kt}} \, \mid \, x_m \- \Z_{p^{km}}, \phi_m (x_m) = x_{m - 1} \}$$
where,  $$\phi_m : \Z_{p^{km}} \longrightarrow \Z_{p^{k(m - 1)}}$$
$$~~~~~~~~~~~~ a \, mod \, p^{km} \mapsto a \, mod \, p^{k(m - 1)}$$
\end{definition}
Thus $\Z_{(n)}$ can be described as follows: 
$$\zn = \{x = (\cdots, x_2, x_1) \, \in \, \prod_{t = 1} ^{\infty} \Z_{n^t}\, \mid \, x_i \, \in \, \Z_{n^i},\phi_i (x_i) = x_{i - 1} \, \, \forall i \geq 2 \}$$
Let $x = (\cdots, x_2, x_1) \, \in \, \zn$.
From the $(*)$ each $x_i$ corresponds to a tuple $(x_i ^1, x_i ^2, \cdots, x_i ^r) \- \Z_{p_1 ^{k_1 i}} \times \Z_{p_2 ^{k_2 i}} \times \cdots \times \Z_{p_r ^{k_r i}}$
and map $\phi_i = \phi_i ^1 \times \phi_i ^2 \times \cdots \times \phi_i ^r$
where $\phi_i ^j : \Z_{p_j ^{k_j i}} \rightarrow \Z_{p_j ^{k_j (i - 1)}}$ sending $a \, (mod \, p_j ^{k_j i})$ to $a \, (mod \, p_j ^{k_j (i - 1)})$.
$\zn$ is a ring with `coordinate by coordinate' addition and multiplication. This is possible because $\phi_n$ is a ring homomorphism. 
$\zn$ is a commutative ring with identity $(\cdots, 1, 1).$\\

By definition of the projective limit,
for $x = (\cdots, x_2, x_1) \- \zn$ if $x_m = 0$ for some $m$, then $\phi_m$'s are all ring homomorphisms implies that
 $x_{m - 1}, x_{m - 2}, \cdots, x_1$ all are equal to zero.
Therefore, $x$ looks like  $(\cdots, x_{m + 1}, 0, 0, \cdots, 0).$

\section{Algebraic properties of ring $\zn$} \label{C}

We recall some standard definitions in commutative algebra. Let $R$ denote the commutative ring with unity.

\begin{definition}
$R$ is called an \textbf{integral domain} if for $x, y \- R$, we have $x y = 0$ then either $x = 0$ or $y = 0$.
\end{definition}

\begin{definition}
 $R$ is \textbf{principal ideal domain} (PID) if every ideal of $R$ is generated by a single element. 
\end{definition}

\begin{definition}
 A commutative ring having unique maximal ideal is called as \textbf{local ring}.
\end{definition}

\begin{definition}
  A commutative ring having finitely many maximal ideals is called as \textbf{semilocal ring}.
\end{definition}

Now we prove some algebraic properties of the ring of $n$-adic integers.

\begin{theorem} \label{6}
Two rings $\zp$ and $\Z_{(p^k)}$ are isomorphic. 
\end{theorem}

\begin{proof}
 $\zp = \varprojlim \Z_{p^n}$ and $\Z_{(p^k)} = \varprojlim \Z_{p^{kn}}$. The directed set for $\zp$ is $\N$ and for $\Z_{(p^{kn})}$ is $k \N = \{k, 2k, \cdots\}$.
Then $k \N \subset \N$ and for every $n \- \N$, there exist $km \- k \N$ such that $km \geq n$. Thus the projective system $(\Z_{p^{kn}}, k \N)$ is cofinal in 
the projective system $(\Z_{p^n}, \N)$. Therefore, there projective limits are isomorphic. That is, $\zp \cong \Z_{p^k}.$
\end{proof}

\begin{theorem} Let $p$ be a prime number. 
Then every non-zero element $x$ of $\zp$ can be uniquely written as $p^n u$ for some $n \- \N$ and $p$ not dividing $u.$ \end{theorem} 

\begin{proof}Let $x = (\cdots, x_3, x_2, x_1)$ be non-zero element of $\zp$. Let $n$ be the smallest number such that $x_n \neq 0$. 
Then $x$ is of the form $(\cdots, x_{n + 1}, x_n, 0, \cdots, 0)$. By the definition of $\zp$ we get $x_n \eq 0 \, (mod \, p^{n - 1}).$
That is $p^{n - 1}$ divides $x_n$. Since $x_n \- \Z_{p^n}$ and $x_n \neq 0$ we get $x_n = p^{n - 1} y_n$ where $p \nmid y_n$. 
For $m \geq n$, by the definition of $\zp$, $\, x_m \eq x_n \, (mod \, p^{n - 1}).$ Since $p^{n - 1}$ divides $x_n$, we have $p^{n - 1}$ also divides $x_m$, for all $m \geq n.$
Therefore $p^{n - 1}$ divides $x$. Also $p \nmid y_n$ implies $p^n \nmid x$. We can write $x = p^{n - 1} u$,  where $u = (\cdots, y_{n + 1}, y_n)$. 
\end{proof}

\begin{theorem}
 Let $u \- \zp$. Then $u$ is a unit if and only if $p$ does not divide $u$.
\end{theorem}

\begin{proof}
Let $u = (\cdots, u_3, u_2, u_1)$ be a unit in $\zp$. Let $v = (\cdots, v_3, v_2, v_1)$ be the inverse of $u$ in $\zp$. 
By the multiplication in $\zp$, we get $u_1 v_1 \eq 1 \, (mod \, p)$. Therefore $p \nmid u_1$, hence $p$ does not divide $u$.\\
Conversely, Suppose $p$ does not divide $u$. Let $n$ be the smallest number such that $p \nmid u_n$. Consider the map $\phi_n : \Z_{p^n} \rightarrow \Z_p$
given by $a \, (mod \, p^n) \mapsto a \, (mod \, p)$. Then  $p \nmid u_n$ implies $\phi_n (u_n)$ is non-zero. But  $\phi_n (u_n)$ is $u_1$ by definition of $\zp$.
Let $v_1 = u_1 ^{-1}$ in $\Z_p$. This infact shows that the smallest number $n$ such that $p \nmid u_n$ is 1. Since $u_1$ is non-zero in $\Z_p$, and $u_2 \eq u_1 \, (mod \,p)$
 we get $u_2$ is also not divisible by $p$ in $\Z_{p^2}$. So $u_2$ is invertible in $\Z_{p^2}$. Similarly, we get $u_i$ is invertible in $\Z_{p^i}$ for all $i$.
 Let $v_i = u_i ^{-1}$ in $\Z_{p^i}$. Then $v = (\cdots, v_3, v_2, v_1)$ is the inverse of $u$ in $\zp$.
\end{proof}

\begin{theorem}
  If $n = p_1 ^{k_1} p_2 ^{k_2} \cdots p_r ^{k_r}$ then group of units
$\zn ^* \cong \Z_{(p_1 ^{k_1})} ^* \times \Z_{(p_2 ^{k_2})} ^* \times \cdots \times \Z_{(p_r ^{k_r})} ^*$
\end{theorem}

\begin{proof}
 Let $x = (x_1, x_2, \cdots, x_r)$ be an element of $\zn$ and $x_i \- \Z_{(p_i ^{k_i})}$. Then $x$ is invertible in $\zn$ \iff each $x_i$ is invertible in $\Z_{(p_i ^{k_i})}$.
Therefore $\zn ^* \cong \Z_{(p_1 ^{k_1})} ^* \times \Z_{(p_2 ^{k_2})} ^* \times \cdots \times \Z_{(p_r ^{k_r})} ^*$.
\end{proof}

\begin{theorem} \label{19}
$\zn$ is an integral domain if and only if $n$ is power of a prime number.
\end{theorem}

\begin{proof} 
Suppose $n$ is not power of a prime.\\
Let $n = p_1 ^{k_1} p_2 ^{k_2} \cdots p_r ^{k_r}$ and $r \geq 2$. Then $\zn \cong \Z_{(p_1 ^{k_1})} \times \cdots \times \Z_{(p_r ^{k_r})}$.
Let $A = (a, 0, 0,\cdots, 0)$ and $B = (0, b, 0, \cdots, 0)$
are in $\Z_{(p_1 ^{k_1})} \times \cdots \times \Z_{(p_r ^{k_r})}$ where $a \neq 0$ and $b \neq 0$. But $AB = 0$. Therefore $\zn$ is not an integral domain.\\
Conversely, suppose $n = p^r$ for some prime $p$ and $r \- \N$. Then by theorem \ref{6}, we get $\zn$ is isomorphic to $\zp$. Therefore it is enough to prove that $\zp$
is an integral domain. Let $x$ and $y$ be non-zero elements of $\zp$. Then $x = p^n u$ and $y = p^m v$ for some $n, m \- \N $ and $p$ does not divide $u, v$.
Therefore $u, v$ are units in $\zp$. So $uv \neq 0$. Now $x y = p^{m + n} u v$, which is non-zero because, $uv$ is non-zero.
\end{proof}

\begin{theorem}
 $\zp$ is a \textbf{PID}. In fact, every non-zero ideal of $\zp$ is generated by $p^n$ for some $n \geq 0.$
\end{theorem}
 
\begin{proof} Let $I$ be any non-zero ideal of $\zp$. If $I = \zp$ then it is generated by $ 1 = (\cdots, 1, 1)$.
Suppose, $I \neq \zp$. Let $n$ be the largest natural number such that $p^n$ divides every element of $I$.\\
\underline{Claim:} $I$ is generated by $p^n$ in $\zp.$\\
By the choice of $n$, there exists an element $a = p^n u \, \in \, I$ such that $p$ does not divide $u$. Therefore, $u$ is unit in $\zp$.
$I$ is an ideal implies $a u^{-1} \, \in I$ i.e. $p^n$ is in $I$. Therefore, ideal generated by $p^n$ is contained in $I$.\\
Conversely, every element of $I$ is of the type $p^n x$ for $x \, \in \zp$. Therefore, $I \subseteq p^n \zp.$ 
\end{proof}

\begin{theorem}
 $\zp$ is a \textbf{local ring}.
\end{theorem}
 
\begin{proof} If $n > m$ then $p^n \zp \subset p^m \zp$,
because $p^m \, \in \, p^m \zp$ and $p^{n - m} \, \in \, \zp$ implies that $p^m p^{n - m} = p^n \, \in \, p^m \zp$. We have,
$$ \cdots \subset p^n \zp \subset p^{n - 1} \zp \subset \cdots \subset p^2 \zp \subset p \zp \subset \zp$$
Therefore, ideal generated by $p$ is the unique maximal ideal of $\zp.$ 
\end{proof}

\begin{theorem}
$\zn$ is a \textbf{semilocal ring.} 
\end{theorem}

\begin{proof}
 Let $n = p_1 ^{k_1} p_2 ^{k_2} \cdots p_r ^{k_r}$. Then, 
\begin{eqnarray*}
\zn &\cong& \Z_{(p_1 ^{k_1})} \times \cdots \times \Z_{(p_r ^{k_r})} \\
    &\cong& \Z_{(p_1)} \times \cdots \times \Z_{(p_r)}
\end{eqnarray*}
Since, ideal generated by $p_i$ is the unique maximal ideal of $\Z_{(p_i)}$, all the maximal ideals of $\zn$ are given by 
$$ M_i = \Z_{(p_1)} \times \cdots \times p_i \, \Z_{(p_i)} \times \cdots \times \Z_{(p_r)} \hspace*{4mm}  i = 1, 2, \cdots, r.$$
\end{proof}

\begin{theorem}
$\zp/{p \zp}$ is isomorphic to $\mathbb{F}_p$ a field of order $p$. 
\end{theorem}

\begin{proof} Since ideal generated by $p$ is the maximal ideal of $\zp$, the quotient $\zp/{p \zp}$ must be a field. We show that this field is of order exactly $p$.
So that, it is isomorphic to $\mathbb{F}_p$. We can write,\\
$$\dfrac{\zp}{p \zp} = \{ p \zp, \, 1 + p \zp, \cdots, (p - 1) + p \zp \}$$\\
Every element of $\zp$ is of the form $x = (\cdots, x_2, x_1)$ where $x_1 \- \{0, 1, \cdots, p - 1\}$.
Then $x$ belong to the coset $x_1 \, + \, p \zp$. Hence, we have listed all the cosets.
Therefore Order of $\zp/{p \zp}$ is equal to $p$.
Therefore, $\zp/{p \zp} \cong \mathbb{F}_p.$ 
\end{proof}

\begin{theorem}
For $p$ and $q$ distinct primes, $\Z_{(p)}$ and $\Z_{(q)}$ are not isomorphic as rings.
\end{theorem}

\begin{proof} If they were isomorphic then there unique maximal ideals would also be isomorphic and hence the quotient fields $\mathbb{F}_p$ and $\mathbb{F}_q$
will become isomorphic. This is not true because,  $\mathbb{F}_p$ is of order $p$ and $\mathbb{F}_q$ is of order $q.$
\end{proof}

%
%

\textbf{\underline{The field of $p$-adic numbers, $\qp$}}\\
From the theorem \ref{19} it is clear that, quotient field of $\zn$ is defined only when $n$ is a power of prime.
\begin{definition} The field of $p$-adic numbers $\qp$ is defined as the quotient field of $\zp$. \end{definition}
If $n = p^r$ for some prime $p$ and $r \geq 1$ then by theorem \ref{6}, we get $\Z_{(p^r)} \cong \zp.$ Therefore $\Q_{(p^r)} \cong \qp.$
Every element of $\qp$ can be written as $p^n x$ where $n \- \Z$ and $x \- \zp.$
$$\qp \cong \Z \times \zp.$$
If $n = p_1 ^{k_1} p_2 ^{k_2} \cdots p_r ^{k_r}$, then define $\qn$ as $\Q_{(p_1)} \times \Q_{(p_2)} \times \cdots \times \Q_{(p_r)}$.
Clearly, $\qn$ is a field if and only if $n$ is power of a prime number. $\hspace*{3em} \square$\\

Let $x = (\cdots, x_2, x_1) \, \in \, \zn$. Then 
\begin{eqnarray*}
x_1 \, \in \, \Z_n &\Rightarrow& 0 \leq x_1 \leq n - 1
\end{eqnarray*}
\noindent Let \, $a_o = x_1$
\begin{eqnarray*}
x_2 \- \Z_{n^2}, \, \phi_2 (x_2) = x_1 &\Rightarrow& x_2 \eq x_1 \, (mod \, n)\\
                                       &\Rightarrow& x_2 = a_o + a_1 n \, for \, some \, 0 \leq a_1 \leq n - 1\\
x_3 \- \Z_{n^3}, \, \phi_3 (x_3) = x_2 &\Rightarrow& x_3 \eq x_2 \, (mod \, n^2) \\
                                       &\Rightarrow& x_3 = x_2 + a_2 n^2 \, for \, some\,  0 \leq a_2 \leq n - 1\\
                                       &\Rightarrow& x_3 = a_o + a_1 n + a_2 n^2
\end{eqnarray*}
\noindent In \, general,
\begin{eqnarray*}
 x_i \- \Z_{n^i}, \, \phi_i (x_i) = x_{i - 1} &\Rightarrow& x_i \eq x_{i - 1} \, (mod \, n^{i - 1})\\
                                             &\Rightarrow& x_i = x_{i - 1} + a_{i - 1} n^{i - 1} \, for \,some \, 0 \leq a_{i - 1} \leq n - 1\\
                                             &\Rightarrow& x_i = a_o + a_1 n + a_2 n^2 + \cdots + a_{n - 1} n^{i - 1}
\end{eqnarray*}

These equations actually tells us that,  $n$-adic integers can be defined in another way-
 power series in $n$ with coefficients in $\{0, 1, \cdots, n - 1\}.$

\section{Metric approach to $n$-adic integers}

\underline{\textbf{$n$-adic valuation and $n$-adic metric:}}\\
Fix a prime $p$. Define \textbf{$p$-adic valuation $v_p$} as,
\begin{eqnarray*}
 v_p : \Z &\longrightarrow& \{0\} \cup \N \cup \{\infty \}\\
v_p (x) &=& m  \hspace{4em}  x = p^m x_o, \, p \nmid x_o\\
        &=& \infty \hspace{4em}  x = 0
\end{eqnarray*}

\begin{example} $v_2 (32) = 5$, $v_7 (-98) = 2$ \end{example}

Using $p$-adic valuation we define a function $\v_n$ on $\Z$ as follows:\\
Let $x$ be any non-zero integer. If $p_1, p_2, \cdots, p_r$ are all the primes appearing in prime factorisation of $n$, then define
 $$ v_n (x) = max\, \{v_{p_1} (x), \, v_{p_2} (x), \cdots, v_{p_r} (x)\}$$

\begin{example} $$v_6 (32) = max \, \{v_2 (32), v_3 (32)\} = max \, \{5, 0\} = 5$$
 $$v_{36} (18) = max \, \{v_2 (36), v_3 (36)\} = max \, \{2, 2\} = 2$$ \end{example}
\begin{remark}
Note that, $v_n$ is not a valuation on $\Z.$ I would like to thank Prof. Shripad Garge, IIT Bombay, for pointing out this fact. 
\end{remark}


Now, define \textbf{$n$-adic metric} from the function $v_n$.
If $x$ and $y$ are integers then define $$ d (x , y) = e^{-v_n (x - y)}.$$

\begin{proposition}
 $(\Z, d)$ is a metric space.
\end{proposition}

\begin{proof}
We check that, $d$ satisfies the metric properties:
Let $x$, $y$ and $z$ be any elements of $\Z$.
\begin{enumerate}
 \item  $d (x, y) = e^{-v_n (x - y)} \geq 0$
\item $d (x, y) = 0$ if and only if $ e^{-v_n (x - y)} = 0$ if and only if $ v_n (x - y) = \infty $ if and only if $ x - y = 0$
\item $v_n (x - y) = v_n (y - x)$\\
$~~ \Rightarrow e^{-v_n (x - y)} = e^{-v_n (y - x)}$
\indent i.e. $d (x, y) = d (y, x)$
\item Triangle inequality: $d (x, y) \leq d (x, z) + d (z, y)$ \\
i.e. $e^{-v_n (x - y)} \leq e^{-v_n (x - z)} + e^{-v_n (z - y)}$\\
\underline{Proof:} \begin{eqnarray*}
x - y &=& (x - z) + (z - y)\\
\Rightarrow v_n (x - y) &\geq& min \, \{v_n (x - z), \, v_n (z - y) \}
\end{eqnarray*}
Suppose, minimum is $v_n (x - z).$ Then,
\begin{eqnarray*}
 v_n (x - y) &\geq& v_n (x - z)\\
 - v_n (x - y) &\leq& - v_n (x - z)\\
 e^{-v_n (x - y)} &\leq& e^{-v_n (x - z)}\\
 e^{-v_n (x - y)} &\leq& e^{-v_n (x - z)} + e^{-v_n (z - y)}\\
i.e. \, d (x, y) &\leq& d (x, z) + d (z, y)
 \end{eqnarray*}
\end{enumerate}
\end{proof}

$d$ actually satisfies stronger property than the triangle inequality, viz. the \textbf{ultrametric} property:

\begin{proposition} Ultrametric property: For all $x, y, z$ in $\Z$,
$$ d (x, y) \leq \, max \, \{d (x, z), d(y, z)\}$$
~~~~~~~~~~~~~~~~~~~~~~~~~i.e. $e^{-v_n (x - y)} \leq \, max \, \{e^{-v_n (x - z)} , e^{-v_n (z - y)}\}$
\end{proposition}

\begin{proof}
Suppose $max \, \{e^{-v_n (x - z)} , e^{-v_n (z - y)}\} = e^{-v_n (x - z)}$. Therefore,
\begin{eqnarray*}
 e^{-v_n (z - y)} &\leq& e^{-v_n (x - z)}\\
-v_n (z - y) &\leq& -v_n (x - z)\\
v_n (z - y) &\geq& v_n (x - z)\\
min \, \{v_n (x - z), \, v_n (z - y)\} &=& v_n (x - z)\\
min \, \{v_n (x - z), \, v_n (z - y)\} &\leq& v_n ((x - z) + (z - y)) = v_n (x - y)\\
v_n (x - z) &\leq& v_n (x - y)\\
 e^{-v_n (x - z)} &\geq&  e^{-v_n (x - y)}
\end{eqnarray*}
i.e. $e^{-v_n (x - y)} \leq \, max \, \{e^{-v_n (x - z)} , e^{-v_n (z - y)}\} $
\end{proof}

Define the set of \textbf{$n$-adic integers $\zn$} as completion of $\N$ with respect to the \\ 
$n$-adic metric.\\
Let $x$ be any positive integer. Then $x$ can be uniquely written as $$ x_o + x_1 n + x_2 n^2 + \cdots + x_k n^k $$ where $x_i \- \{0, 1, 2, \cdots, n - 1\}$.
That is, $x$ can be written as polynomial in $n$ with coefficient in $\{0, 1, 2, \cdots, n - 1\}$. Clearly, power series in $n$, 
$$ x_o + x_1 n + x_2 n^2 + \cdots + x_k n^k  + \cdots \h3 (*)$$  do not make sense in integers, but we will see that it makes sense as an $n$-adic integer.\\
Consider a sequence of integers obtained from  $(*)$, where $x_i \- \{0, 1, 2, \cdots, n - 1\}$.
\begin{eqnarray*}
X_o &=& x_o\\
X_1 &=& x_o + x_1 n\\
X_2 &=& x_o + x_1 n + x_2 n^2\\
\vdots\\
X_k &=& x_o + x_1 n + x_2 n^2 + \cdots + x_k n^k\\
\vdots
\end{eqnarray*}
and let $X = x_o + x_1 n + x_2 n^2 + \cdots + x_k n^k + \cdots$ be just an expression. We will give meaning to this $X$. Let $k, m \- \N$ and $k > m$. Then 
\begin{eqnarray*}
X_k - X_m &=& x_{m + 1} \, n^{m + 1} + \cdots + x_k \, n^k\\
          &=& n^{m + 1} \, (x_{m + 1} + x_{m + 2} \, n + \cdots + x_k \,  n^{k - (m + 1)})
\end{eqnarray*}
WLOG we can assume that, $x_{m + 1} \neq 0$. Therefore, 
\begin{eqnarray*}
v_n (X_k - X_m) &=& m + 1\\
d (X_k, X_m) &=& e^{- v_n (X_k - X_m)} = e^{- (m + 1)}
\end{eqnarray*}
$e^{- (m + 1)}$ can be made as small as possible for all $m \geq m_o$ for some $m_o$.\\
i.e. we have proved that, $\{X_k\}$ is a Cauchy sequence in $\N$ with $n$-adic metric.


\underline{Define:} $\lim_{k \rightarrow \infty} X_k = X = x_o + x_1 n + x_2 n^2 + \cdots$\\
Therefore $\zn$ can be described as follows:
$$ \zn = \{X = x_o + x_1 n + x_2 n^2 + \cdots \, \mid \, x_i \- \{0, 1, 2, \cdots, n - 1\} \, \}$$

$X$ is a positive integer if and only if all but finitely many $x_i$'s are zero.
%
%
%

$\zn$ can be given a ring structure by defining addition and multiplication in a similar (not same) way as
we define addition and multiplication of power series.\\

\textbf{Addition:} \\
Let $X = x_o + x_1 n + x_2 n^2 + \cdots$ and $Y = y_o + y_1 n + y_2 n^2 + \cdots$ be elements of $\zn$ 
then $(x_o + y_o) + (x_1 + y_1) n + (x_2 + y_2) n^2 + \cdots$ may not be in $\zn$, because, 
$x_i$ and $y_i$ lies between 0 and $n - 1$ implies that $x_i + y_i$ lies between 0 and $2 (n - 1)$.\\

Write $X = x_o + x_1 n + x_2 n^2 + \cdots$ as $\cdots x_2 \, x_1 \, x_o$\\
and $Y = y_o + y_1 n + y_2 n^2 + \cdots$ as $\cdots y_2 \, y_1 \, y_o$\\
(This is exactly similar to decimal expansion of integer.\\
e.g. in the decimal system $1 + 2*10 + 3* 10^2 + 4* 10^3$ corresponds to the integer 4321)\\
Like the addition of integers here also, we will do the same columnwise addition of digits from right to left, with `remainder-carry' rule.\\
Formulation of these rules in general terms will contain only complicated notations but will not add anything more to the knowledge. Hence we explain the addition
and multiplication rules by following examples:

\begin{example} Consider $\Z_{(7)} = \{\, \displaystyle\sum_{i = 0} ^{\infty} x_i 7^i \, \mid \, x_i \- \{0, 1, 2, 3, 4, 5, 6\} \,\}$\\
Let $X = \cdots 2355 \longleftrightarrow 5 + 5 * 7 + 3 * 7^2 + 2 * 7^3 + \cdots$\\
and $Y = \cdots 1646 \longleftrightarrow 6 + 4 * 7 + 6 * 7^2 + 1 * 7^3 + \cdots$ be elements of $\Z_{(7)}$.\\
Then $X + Y$ can be calculated as follows: 
\begin{table}[h]
  \begin{tabular}{cccccc}
\underline{In $\Z_{(7)}$}\\
 & & 1 & 1 & 1 &  \\
\hline
 & $\cdots$ & 2 & 3 & 5 & 5\\ 
+& $\cdots$& 1 & 6 & 4 & 6\\ 
\hline
& $\cdots$ & 4 & 3 & 3 & 4
 
\end{tabular}

%

\end{table}
\end{example}
\begin{remark} If $X, Y \- \zn$ are actually integers then their addition and multiplication in $\zn$ and in $\Z$ is same with the understanding that,\\
denote $X = x_o + x_1 n + \cdots + x_k n^k $ as $ x_k \cdots x_1 \, x_o$\\
and $Y = y_o + y_1 n + \cdots + y_l n^l$ as $ y_l \cdots y_1 \, y_o$\\
Suppose, $X + Y = Z = z_m \cdots z_1 \, z_o \- \zn$\\
then $(x_o \, + \, x_1 n \, + \, \cdots \, + \, x_k n^k) \,+ \, (y_o \, + \, y_1 n \, + \, \cdots \, + \, y_l n^l)$ is same as \\
$z_o + z_1 n + \cdots + z_m n^m$.
\end{remark}
\begin{example} In $\Z_{(7)}$, we proved $2355 + 1646 = 4334$\\
Now, \begin{eqnarray*}
 2355 &\longleftrightarrow & 5 + 5 * 7 + 3 * 7^2 + 2 * 7^3 = 873\\
 1646 &\longleftrightarrow& 6 + 4 * 7 + 6 * 7^2 + 1 * 7^3 = 671\\
 4334 &\longleftrightarrow& 4 + 3 * 7 + 3 * 7^2 + 4 * 7^3 = 1544
\end{eqnarray*}
and verify that, $873 + 671 = 1544$.
\end{example}
\begin{example} multiplication: 
In $\Z_{(7)}$ we find $35 \times 64.$
\begin{table}[h]
  \begin{tabular}{ccccc}
\underline{In $\Z_{(7)}$}\\
  &  &  & 3 & 5\\ 
 $\times$ &  &  & 6 & 4\\ 
\hline
 &  & 2 & 0 & 6\\
+ & 3 & 1 & 2 & 0\\
\hline
 & 3 & 3 & 2 & 6 
\end{tabular} 
%
\end{table}

Now, \begin{eqnarray*}
 35 &\longleftrightarrow & 5 + 3 * 7 = 26\\
 64 &\longleftrightarrow& 4 + 6 * 7 = 46\\
 3326 &\longleftrightarrow& 6 + 2 * 7 + 3 * 7^2 + 3 * 7^3 = 1196
\end{eqnarray*}
and verify that, $26 \times 46 = 1196$.
\end{example}

\textbf{Power series representation of some special numbers}\\
We have seen that, the $n$-adic power series expansion of any natural number is finite i.e. only finitely many $n^i$ have non-zero coefficient.
We now represent negative integers as elements of $n$-adic integers $\zn$.

\begin{enumerate}
\item The infinite series $1 + p + p^2 + \cdots  \textit{converges to} \, \dfrac{1}{1 - p} $ in $\zp$
\begin{proof} Let $a_n = 1 + p + p^2 + \cdots + p^n = \dfrac{p^{n + 1} - 1}{p - 1}$.\\
We use the fact that, valuation of $p^n$ goes to infinity as $n \rightarrow \infty$, implies that $p^n$ converges to 0 in $\zp$.
\begin{eqnarray*}
 1 + p + p^2 + \cdots  &=& \lim_{n \rightarrow \infty} \, a_n \\
                       &=& \lim_{n \rightarrow \infty} {\dfrac{p^{n + 1} - 1}{p - 1}}\\
                       &=& \lim_{n \rightarrow \infty} \dfrac{p^{n + 1}}{p - 1} - \lim_{n \rightarrow \infty} \dfrac{1}{p - 1}\\
                       &=& 0 - \dfrac{1}{p - 1}\\
                       &=& \dfrac{1}{1 - p}
\end{eqnarray*}
\end{proof}

\item (From the above fact) \, $-1 = (p - 1) + (p - 1)p + (p - 1)p^2 + \cdots $
\begin{proof}
 \begin{eqnarray*}
-1 &=& \dfrac{p - 1}{1 - p}\\
   &=& (p - 1) \dfrac{1}{1 - p}\\
   &=& (p - 1) \, (1 + p + p^2 + \cdots)\\
   &=& (p - 1) + (p - 1)p + (p - 1)p^2 + \cdots
            \end{eqnarray*}
\end{proof}

\item If $n$ is a natural number then, $-n = \displaystyle\sum_{i = 0} ^{\infty} x_i p^i$,
where all but finitely many $x_i$'s are equal to $p - 1$.
\begin{proof} Case 1: If $n = p^k$ then 
\begin{eqnarray*}
-n &=& p^k \, ((p - 1) + (p - 1)p + (p - 1)p^2 + \cdots)\\
   &=& (p - 1) p^k + (p - 1)p^{k + 1} + (p - 1)p^{k +2} + \cdots
\end{eqnarray*}
Case 2: Suppose $n$ is not equal to a power of $p$. \\
Let $k$ be the smallest integer such that, $p^k > n$. Then the $p$-adic expansion of $p^k - n$ will contain powers of $p$ upto $ k - 1$. 
\begin{eqnarray*}
 p^k - n &=& x_o + x_1 \, p + x_2 \, p^2 + \cdots + x_{k - 1} \, p^{k - 1}\\
     -n  &=& x_o + x_1 \, p + x_2 \, p^2 + \cdots + x_{k - 1} \, p^{k - 1} - p^k\\
         &=& x_o + x_1 \, p + x_2 \, p^2 + \cdots + x_{k - 1} \, p^{k - 1} + p^k \, ((p - 1) + (p - 1)\, p + \cdots)\\
         &=& x_o + x_1 \, p + x_2 \, p^2 + \cdots + x_{k - 1} \, p^{k - 1} + (p - 1) \, p^k + (p - 1) \, p^{k + 1} + \cdots 
\end{eqnarray*}
\end{proof}

\item If $n$ is coprime to $p$ then we find $p$-adic expansion of $1/n$.\\
Answer: $(n, p) = 1$, therefore there exists integers $a$ and $b$ such that $an + bp = 1.$\\
Let $a = \displaystyle\sum_{i = 0} ^s a_i p^i$ and $b = \displaystyle\sum_{i = 0} ^t b_i p^i$ be $p$-adic expansions of $a$ and $b$ respectively.
\begin{eqnarray*}
 an + bp = 1 &\Rightarrow&  an = 1- bp\\
             &\Rightarrow&  \dfrac{1}{n} = \dfrac{a}{1 - bp}\\
             &\Rightarrow& \dfrac{1}{n} = (\displaystyle\sum_{i = 0} ^s a_i p^i) * (\displaystyle\sum_{i = 0} ^{\infty} (bp)^i)\\
             &\Rightarrow& \dfrac{1}{n} = (\displaystyle\sum_{i = 0} ^s a_i p^i) * (\displaystyle\sum_{i = 0} ^{\infty} ((\displaystyle\sum_{j = 0} ^t b_j p^j)p)^i).
\end{eqnarray*}
\end{enumerate}

%
%
%
%
%
%
%
%
%
%
%
%
%
%
%
%
\newpage
\section{Teichm\"{u}ller Units} \label{D}

In this section we discuss about units in $p$-adic integers which are of finite order. We will give the structure of the group of units $U$ of $p$-adic integers.

\begin{definition} 
Finite order elements of $\zp ^*$ are called as \textbf{Teichm\"{u}ller units}.
\end{definition}

 Let $V$ denote the set of all Teichm\"{u}ller units. $V$ forms a subgroup of $\zp ^*,$ group of units in $p$-adic integers.
If $a = a_o + a_1 p + a_2 p^2 + \cdots$ is an element of $V$ of order $k$ then $a \, (mod \, p^n)$ is of order $k$ in $\n$, for every $n \geq 1.$
It is not at all clear from the definition that how many elements belong to the group $V.$ 

\begin{theorem}  \label{8}
 For an odd prime $p$, $\zp ^* \cong \Z_{p - 1} \times \zp$\\
For $p =2$, $\Z_{(2)} ^* \cong \Z_2 \times \Z_{(2)}$.
\end{theorem}

\begin{proof} Let $p$ be an odd prime. Then,
$$ \n \cong \, \Z_{p - 1} \times \Z_{p^{n - 1}}.$$
Taking the projective limit on both sides, as $n \rightarrow \infty$,
$$ U = \zp ^* \cong \, \Z_{p - 1} \times \zp$$
Let $p = 2$. Then for $n \geq 2$,
$$ \t \cong \, \Z_2 \times \Z_{\w}$$
taking the projective limit on both sides, as $n \rightarrow \infty$,
$$ \Z_{(2)} ^* \cong \, \Z_2 \times \Z_{(2)}$$
\end{proof}

Each element of $(\zp, \, +)$ is of infinite order. Above theorem proves that, $V$ is finite and as a group it is isomorphic to the cyclic group of order $p - 1$ if $p$ is an odd
 prime and $V$ is isomorphic to the cyclic group of order 2, when $p = 2.$
We now give the actual construction of Teichm\"{u}ller units.\\

\textbf{Construction of Teichm\"{u}ller units:}\\
Let $a = a_o + a_1 p + a_2 p^2 + \cdots $ where $a_o, a_1, a_2, \cdots \, \in \{0, 1, 2, \cdots, p - 1\}$ \\
Recall: $a \, \in \, \Z_{(p)}$ is a unit if and only if $a_o \neq 0.$
%
%
%
\begin{theorem} 
If $p = 2$ then $V = \{1, 1 + 2 + 2^2 + \cdots + 2^n + \cdots \}$.  
\end{theorem}

\begin{proof}
 From the theorem \ref{8}, group of Teichm\"{u}ller units in $\Z_{(2)}$ is of order 2. Therefore it is enough to prove that, $1 + 2 + 2^2 + \cdots$ is of order 2 in $\Z_{(2)} ^*.$
But $1 + 2 + 2^2 + \cdots + 2^n + \cdots = -1$ in 2-adic integers and $(-1)^2 = 1$ in $\Z_{(2)}.$
\end{proof}

\begin{theorem} \label{18}
Let $p$ be an odd prime. Teichm\"{u}ller units are in one to one correspondence with the non-zero constant term.
i.e. if $a = a_o + a_1 p + a_2 p^2 + \cdots $ is a Teichm\"{u}ller unit then $a$
is uniquely determined by the constant term $a_o \neq 0$ and conversely given $a_o \- \{1, 2, \cdots, p - 1\}$ there exists unique Teichm\"{u}ller unit having constant term $a_o.$
\end{theorem}

\begin{proof} 
\textbf{Construction of Teichm\"{u}ller units from it's constant term:}\\
Let $a_o \, \neq 0$ be given and $a_o \, \in \, \{1, 2, \cdots, p - 1\}$. Let $o(a_o) = k.$ Then $k \mid {p - 1}.$
We want to find $a_i$'s between 0 to $ p - 1$ such that order of $a = a_o + a_1 p + \cdots + a_n p^n + \cdots$ is $k$. Therefore,
$$a^k= (a_o + a_1 p + a_2 p^2 + \cdots )^k = 1$$
\noindent Using $a_o$ we can find $a_1$ ,\\
Using $a_o, a_1$ we can find $a_2$,\\
Using $a_o, a_1, a_2$ we can find $a_3$,\\
and so on.\\
$a_o \, \in \{1, 2, \cdots, p - 1\}$ is of order $k$, modulo $p$. Therefore $a_o ^k \eq 1 \, (mod \, p).$\\
Then there exists unique $a_1 \, \in \, \{0, 1, 2, \cdots, p - 1\}$ such that  $(a_o + a_1 p)^k \eq 1 \, (mod \, p^2).$
Following is the construction of $a_1$:
\begin{eqnarray*} 
 (a_o + a_1 p)^k &=& {a_o}^k + k {a_o}^{k - 1} a_1 p + \frac{k (k - 1)}{2} {a_o}^{k - 2} a_1 ^2 p^2  \\
                      &\vspace{1mm}&\\
                       &+& \cdots + {\binom k r } {a_o}^{k - r} {a_1}^r p^r + \cdots + {a_1}^k p^k \\
		       &\vspace{1mm}&\\
                       &\eq&  {a_o}^k + k {a_o}^{k - 1} a_1 p \, (mod \, p^2) \\
(a_o + a_1 p)^k &\eq& 1 \, (mod \, p^2) \\
\Leftrightarrow {a_o}^k + k {a_o}^{k - 1} a_1 p  &\eq& 1 \, (mod \, p^2)  \\
\Leftrightarrow k {a_o}^{k - 1} a_1 p  &\eq& 1 - a_o ^k \, (mod \, p^2)  \\
\Leftrightarrow k a_o ^{k - 1} a_1 &\eq& \dfrac{1 - a_o ^k}{p} \, (mod \, p) \h3 (*)
\end{eqnarray*}
Since, $a_o ^k \eq 1 \, (mod \, p)$, we get $a_o ^{k - 1} = {a_o}^{-1} \, (mod \, p)$. Also $k \mid \,{p - 1}$ implies $(k, p) = 1.$ 
Therefore $k^{-1}$ exists, mod $p.$
Therefore 
$(*)$ can be written as 
\begin{eqnarray*}
k {a_o}^{-1} a_1 &\eq&  \dfrac{1 - a_o ^k}{p} \, (mod \, p)\\
 a_1  &\eq& k^{-1} a_o \dfrac{1 - a_o ^k}{p}  \, (mod \, p)\\
      &\eq& k^{-1} a_o \dfrac{a_o ^k - 1}{p}  (p - 1) \, (mod \, p)
\end{eqnarray*}
$a_o$ is known, hence right hand side can be calculated and modulo $p$ 
operation gives that, $a_1 \, \in \, \{0, 1, \cdots, p - 1\}$.
Uniqueness of $a_1$ follows from the fact that $(*)$ is linear in $a_1$.

\begin{example}
Take $p = 5$. Given $a_o$ we find $a_1$ using above formula:
 \begin{table}[h]
 \begin{tabular}{|c|c|c|c|c|c|}
 \hline
  $a_o$ & $k$ & $k^{-1}$ & $\dfrac{a_o ^k - 1}{p}$ & $A_1 = k^{-1} a_o \dfrac{a_o ^k - 1}{p}  (p - 1)$ & $a_1 = A_1 \, (mod \, p)$ \\ 
  \hline
  \hline
  1 & 1 & 1 & 0 & 0 & 0\\
  \hline
  2 & 4 & 4 & 3 &  96 & 1\\
  \hline
  3 & 4 & 4 & 16 & 768 & 3\\
  \hline
  4 & 2 & 3 & 3 & 144 & 4\\
  \hline
 \end{tabular}
\end{table}
\end{example}
\newpage
\underline{To find $a_2$ using $a_o$ and $a_1$}\\
Given $a_o ^k \eq 1 \, (mod \, p)$ and \\
$(a_o + a_1 p)^k \eq 1 \, (mod \, p^2).$\\
Find $a_2 \, \in \, \{0, 1, 2, \cdots, p - 1\}$ such that $(a_o + a_1 p + a_2 p^2)^k \, \eq \, 1 \, (mod \, p^3)$. 
\begin{eqnarray*}
 (a_o + a_1 p + a_2 p^2)^k &=& ((a_o + a_1 p) + a_2 p^2)^k\\
                                 &\vspace{0.2cm}& \\
                                 &=& (a_o + a_1 p)^k + k \, (a_o + a_1 p)^{k - 1} \, a_2 \, p^2\\
                                 &\vspace{2mm}&\\
                                 &+& \frac{ k (k - 1)}{2} \, (a_o + a_1 p)^{k - 2} \, a_2 ^2 \, p^4 + \cdots + (a_2 p^2)^k\\
                                 &\vspace{0.2cm}&\\
                                 &\eq& (a_o + a_1 p)^k + k \, (a_o + a_1 p)^{k - 1} \, a_2 \, p^2 \, (mod \, p^3)\\
				 &\vspace{0.2cm}&\\
 (a_o + a_1 p + a_2 p^2)^k &\eq& 1 \, (mod \, p^3)\\
                                 &\vspace{0.2cm}&\\
\Leftrightarrow (a_o + a_1 p)^k + k \, (a_o + a_1 p)^{k - 1} \, a_2 \, p^2   &\eq& 1 \, (mod \, p^3)\\
                                 &\vspace{0.2cm}&\\
\Leftrightarrow  k \, (a_o + a_1 p)^{k - 1} \, a_2 \, p^2 &\eq& 1 - (a_o + a_1 p)^k \, (mod \, p^3)\\
                                 &\vspace{0.2cm}&\\
\Leftrightarrow k \, (a_o + a_1 p)^{k - 1} \, a_2  &\eq& \frac{1 - (a_o + a_1 p)^k }{p^2} \, (mod \, p) \h3 (**)
\end{eqnarray*}
Since $a_o + a_1 p \eq a_o \, (mod \, p)$\\
$\Rightarrow (a_o + a_1 p)^{k - 1} \eq a_o ^{k - 1} = a_o ^{-1} \, (mod \, p)$
Therefore $(**)$ can be written as 
\begin{eqnarray*}
k \, a_o ^{-1} \, a_2  &\eq& \dfrac{1 - (a_o + a_1 p)^k }{p^2} \, (mod \, p)\\
                                      &\vspace{0.2cm}&\\
a_2  &\eq& k^{-1} \, a_o \, \dfrac{1 - (a_o + a_1 p)^k }{p^2} \, (mod \, p)\\
     &\vspace{0.2cm}&\\
     &\eq& k^{-1} \, a_o \, \dfrac{(a_o + a_1 p)^k - 1 }{p^2} (p - 1) \, (mod \, p)  
\end{eqnarray*}
 $a_o$ and $a_1$ are known, hence right hand side can be calculated and modulo $p$ 
 operation gives that, $a_2 \, \in \, \{0, 1, \cdots, p - 1\}$.
 Uniqueness of $a_2$ follows from the fact that $(**)$ is linear in $a_2$.
\begin{example}
For $p = 5$, from the earlier table using values of $a_o$ and $a_1$ we give here values of $a_2$, using the formula
$ k^{-1} \, a_o \, \dfrac{(a_o + a_1 p)^k - 1 }{p^2} (p - 1) \, (mod \, p). $
\begin{table}[h]
\begin{tabular}{|c|c|c|c|c|c|c|}
 \hline
  $a_o$ & $a_1$ & $k$ & $k^{-1}$ & $\frac{(a_o + a_1 p)^k - 1 }{p^2}$ &  $A_2 = k^{-1} a_o \frac{(a_o + a_1 p)^k - 1 }{p^2} (p - 1)$ & $a_2 = A_2 \, (mod \, p)$ \\ 
  \hline
  \hline
  1 & 0 & 1 & 1 & 0 & 0 & 0 \\
  \hline
  2 & 1 & 4& 4 & 96 & 3072 & 2 \\
  \hline
  3 & 3 & 4 & 4 & 4199 & 201552 & 2 \\
  \hline
  4 & 4 & 2 & 3 & 23 & 1104 & 4\\
  \hline
 \end{tabular}
\end{table}
 \end{example}
To find $a_n$, we use induction on $n$. \\
\noindent suppose $a_o, a_1, a_2, \cdots, a_{n - 1}$ are known. 
Let $ x_n = a_o + a_1 p + a_2 p^2 + \cdots + a_{n - 1} p^{n - 1}$ be such that $x_n ^k \eq 1 \, (mod \, p^n)$. 
Then we prove that, there exists unique $a_n \, \in \, \{0, 1, 2, \cdots, p - 1\}$ such that
$(a_o + a_1 p + a_2 p^2 + \cdots + a_n p^n)^k \eq 1 \, (mod \, p^{n + 1})$.\\
Let $x_{n + 1} = a_o + a_1 p + a_2 p^2 + \cdots + a_n p^n$.
\begin{eqnarray*}
 x_{n + 1} ^k &=& (x_n + a_n p^n)^k\\
                     &=& x_n ^k + k x_n ^{k - 1} a_n p^n + \dfrac{k (k - 1)}{2} x_n ^{k - 2} a_n ^2 p^{2n} + \cdots + (a_n p^n)^k\\
		     &\eq& x_n ^k +  k x_n ^{k - 1} a_n p^n \, (mod \, p^{n + 1})\\
x_{n + 1} ^k &\eq& 1 \, (mod \, p^{n + 1})\\
\Leftrightarrow x_n ^k +  k x_n ^{k - 1} a_n p^n  &\eq& 1 \, (mod \, p^{n + 1})\\
\Leftrightarrow  k x_n ^{k - 1} a_n p^n  &\eq& 1 - x_n ^k \, (mod \, p^{n + 1})\\
\Leftrightarrow k x_n ^{k - 1} a_n &\eq& \dfrac{1 - x_n ^k}{p^n}  \, (mod \, p) \h3 (***)
\end{eqnarray*}
Since $x_n \eq a_o \, (mod \, p)$ we get $x_n ^{k - 1} \eq a_o ^{k - 1} = a_o ^{-1} \, (mod \, p)$\\
Therefore $(***)$ can be written as 
\begin{eqnarray*}
k \, a_o ^{-1} a_n &\eq& \dfrac{1 - x_n ^k}{p^n} \, (mod \, p)\\
               a_n &\eq& k^{-1} a_o \, \dfrac{1 - x_n ^k}{p^n} \, (mod \, p)\\
                   &\eq& k^{-1} a_o \, \dfrac{x_n ^k - 1}{p^n} (p - 1) \, (mod \, p)            
\end{eqnarray*}
 $a_o, a_1, \cdots, a_{n - 1}$ are known, hence right hand side can be calculated and modulo $p$ 
 operation gives that, $a_n \, \in \, \{0, 1, \cdots, p - 1\}$.\\

\underline{ Uniqueness of $a_n$ follows from the fact that $(***)$ is linear in $a_n$.}
Therefore, Teichm\"{u}ller units are fixed once their constant term is known.\\
\end{proof}

\begin{example} Using above procedure, suppose $a = a_o + a_1 5 + a_2 5^2 + \cdots$ is the Teichm\"{u}ller unit in the ring of 5-adic integers,
 $\Z_{(5)}$ having constant term $2$.
With the help of `SAGE' we have calculated first few $a_i$'s as follows:
\begin{table}[h]
 \begin{tabular}{|c|c|c|c|c|c|c|c|c|c|c|}
\hline
$a_o$ & $a_1$ & $a_2$ & $a_3$ & $a_4$ & $a_5$ & $a_6$ & $a_7$ & $a_8$ & $a_9$ & $a_{10}$ \\
\hline
2 & 1 & 2 & 1 & 3 & 4 & 2 & 3 & 0 & 3 & 2\\
\hline
\hline
$a_{11}$ & $a_{12}$ & $a_{13}$ & $a_{14}$ & $a_{15}$ & $a_{16}$ & $a_{17}$ & $a_{18}$ & $a_{19}$ & $a_{20}$ & $a_{21}$ \\
\hline
2 & 0 & 4 & 1 & 3 & 2 & 4 & 0 & 4 & 3 & 4 \\
\hline
 \end{tabular}
\end{table}
\end{example}
%

Recall that, we have given the structure of the group of units of $\zp$ in theorem \ref{8}. But now we express $U = \zp ^*$ as an internal direct product of its subgroups.

\begin{theorem} \label{7}
Let $U$ be the group of units of $\zp$, and $V$ be the group of Teichm\"{u}ller units of $\zp.$
Then i) for $p$ an odd prime, $U$ is an internal direct product of $V$ and $U_1.$\\
ii) for $p = 2$, $U = U_1 $ and $U_1$ is an internal direct product of $V$ and $U_2.$
\end{theorem}

\begin{proof}
i) Let $p$ be an odd prime. We have a map from $U$ to $\Z_p ^*$ which maps $x = x_o + x_1 p + x_2 p^2 + \cdots$ to $x_o$. Since $x$ is a unit in $\zp$, 
constant term $x_o$ is non-zero.
Therefore the map is well-defined. Since $\Z_p ^*$ is a group under multiplication with 1 as identity, we get $U_1$ is the kernel of above map. 
Therefore we get an exact sequence $$U_1 \rightarrow U \rightarrow \Z_p ^* \cong \Z_{p - 1}.$$ Since we also have a map from $\Z_p ^*$ to $U$ 
given by $x_o$ goes to the Teichm\"{u}ller unit having constant term $x_o.$ So the composition map from $\Z_p ^* \rightarrow U \rightarrow \Z_p ^*$ is
identity. Therefore above exact sequence splits. Hence $U \cong U_1 \times \Z_p ^*$. As $\Z_p ^*$ is isomorphic to $V$ the group of Teichm\"{u}ller
units, $U$ is an internal direct product of $U_1$ and $V$.\\
ii) Let $p = 2$. Then $x = x_o + x_1 p + x_2 p^2 + \cdots$ is a unit \iff $x_o = 1.$ Therefore, $U = U_1.$ As in the earlier case, define a map from $U_1$ to $\Z_2$ given by
 \hbox{$1 + 2 (x_o + 2 x_1 + \cdots) \longmapsto x_o.$} This is a group homomorphism having kernel
$U_2 = \{1 + 2 x \- U_1 \, | \, x_o = 0\}.$ Therefore we get an exact sequence $$U_2 \rightarrow U_1 \rightarrow \Z_2.$$ Since we also have a map from $\Z_2$ to $U_1$ 
given by $0 \mapsto 1$ and $1 \mapsto \sum_{i = 0} ^{\infty} 2^i.$ So the composition map from $\Z_2 \rightarrow U_1 \rightarrow \Z_2$ is
identity. Therefore above exact sequence splits. Hence $U_1 \cong \Z_2 \times U_2$. As $\Z_2$ is isomorphic to $V$ the group of Teichm\"{u}ller
units, $U_1$ is an internal direct product of $U_2$ and $V$.
\end{proof}

\begin{remark}
 From theorem \ref{8} and theorem \ref{7}, for an odd prime $p$ we see that, $\uone$ must be isomorphic to $\zp$. 
For $p = 2$, $U_2$ is isomorphic to $\Z_{(2)}.$
\end{remark}

\begin{note}
 In the appendix we will prove that $\exp: p\, \zp \rightarrow U_1$ and $\log: U_1 \rightarrow p \, \zp$ maps are group homomorphisms and are inverses of each other. 
In fact, we will prove this in more general set up namely, formal power series ring over rationals $\Q[[x]].$ We note here some important results from the appendix
which we will require for our work. 
\end{note}

\begin{definition} \label{21} Let $p$ be an odd prime.
 \begin{eqnarray*}
 \exp: p \, \zp &\longrightarrow& \uone\\
        px &\mapsto& 1 + px + \dfrac{(px)^2}{2!} +  \cdots + \dfrac{(px)^n}{n!} + \cdots
\end{eqnarray*}
\end{definition}
 Series on the R.H.S. above is convergent in $p$-adic integers iff valuation of the $n^{th}$ term goes to infinity as $n \rightarrow \infty$.
 $$ v_p ((px)^n) \geq n $$ 
We have de Polignac's formula which gives maximum power of a given prime dividing $n!$. Proof of de Polignac's formula can be found in \cite{NZM} page 182.
\begin{eqnarray*}
v_p (n!) &=& \sum_{i = 1} ^{\infty} \displaystyle\left\lfloor \dfrac{n}{p^i}\displaystyle\right\rfloor \\
         &\leq& \sum_{i = 1} ^{\infty} \dfrac{n}{p^i} \\
         &=& \dfrac{n}{p - 1}
\end{eqnarray*}
Therefore, $v_p \dl\dfrac{(px)^n}{n!}\dr = v_p (px)^n - v_p (n!) \geq n - \dfrac{n}{p - 1}$. As $n \rightarrow \infty$, $\dfrac{n}{p - 1} \rightarrow \infty$.
Therefore $\exp$ map is well defined.

\begin{definition} \label{22} Let $p$ be an odd prime.
\begin{eqnarray*}
 \log: \uone &\longrightarrow& p \, \zp\\
 1 + px &\mapsto& px - \dfrac{(px)^2}{2} +  \cdots + (-1)^{n + 1} \dfrac{(px)^n}{n} - \cdots
\end{eqnarray*}
\end{definition}
 Again the series on the R.H.S. is convergent \iff valuation of $\dfrac{(px)^n}{n}$ goes to infinity as $n$ goes to infinity.
As $v_p (n) < n$, and for the subsequence $n_k = p^k$, we have $v_p (p^k) = k$. So $v_p \dl {\dfrac{(px)^{p^k}}{p^k}} \dr \geq p^k - k$.
As $p^k - k$ goes to infinity when $k \rightarrow \infty$, we get $v_p \dl {\dfrac{(px)^{p^k}}{p^k}} \dr \rightarrow \infty$ when $k \rightarrow \infty$.
Therefore $\log$ map is well defined.

\begin{theorem}
$\log \circ \exp = I_{p \, \zp}.$ 
\end{theorem}

\begin{theorem} \label{20}
Let $p$ be an odd prime. The maps $\exp: p\, \zp \rightarrow U_1$ and $\log: U_1 \rightarrow p \, \zp$ are group isomorphisms.
\end{theorem}

Proofs of above two theorems are discussed in more general context of $\Q[[x]],$ in appendix B.
Appendix can be referred only when $p$ is an odd prime. We consider the case $p = 2$ separately.

\textbf{$\exp$ and $\log$ map for 2-adic integers:}\\
In contrast with the case of odd prime, $\exp$ map is not defined on $U_1 = 1 + 2 \, \Z_{(2)}.$ But it is defined on $U_2 = 1 + 2^2 \, \Z_{(2)}$. 

\begin{definition}
 \begin{eqnarray*}
  \exp : 2^2 \, \Z_{(2)} &\rightarrow& 1 + 2^2 \, \Z_{(2)}\\
          2^2 x &\mapsto& 1 + 2^2 x + \dfrac{(2^2 x)^2}{2!} + \cdots + \dfrac{(2^2 x)^n}{n!} + \cdots
 \end{eqnarray*}
\end{definition}

\begin{definition} \label{29}
 \begin{eqnarray*}
  \log : 1 + 2^2 \, \Z_{(2)} &\rightarrow&  2^2 \, \Z_{(2)}\\
         1 + 2^2 x &\mapsto& 2^2 x - \dfrac{(2^2 x)^2}{2} + \cdots + (-1)^{n + 1} \dfrac{(2^2 x)^n}{n} - \cdots
 \end{eqnarray*}
\end{definition}

Calculations as we have done in definitions \ref{21} and \ref{22}, will also work here, and we will get above two series on R.H.S. are convergent in $\zp.$
Proofs of $\exp$ and $\log$ are group isomorphisms is as same in case of $p$ an odd prime. Hence refer appendix B.

\newpage
\section{Structure of $\qp ^* / {\qp ^*}^k$}
We denote  the set $1 + p^k \, \zp$ by $U_k$, $k \geq 1$.
In the previous section, we have proved that, the $\log$ map gives group isomorphism from $U_1$ to $p \, \zp$.  
In this section we will prove that, multiplicative group $U_k$ is isomorphic to the additive group $p^k \, \zp$ by the $\log$ map.
We shall also discuss about the map $\phi_k : U_1 \rightarrow U_{k + 1}$ given by $x \mapsto x^k$.

\begin{theorem}
The $\log$ map is an isomorphism from the group $(U_k, *)$ to the group $(p^k \, \zp, +)$. 
\end{theorem}

\begin{proof}
 Consider the restriction of $\log$ map on $U_k = 1 + p^k \, \zp.$ Let $1 + p^k x$ be any element of $U_k$. Then 
$$\log (1 + p^k x) = p^k x - \dfrac{(p^k x)^2}{2} + \dfrac{(p^k x)^3}{3} - \cdots + (-1)^{n + 1} \dfrac{(p^k x)^n}{n} - \cdots$$
Clearly, R.H.S. of the above equation belongs to $p^k \, \zp$. Therefore, $\log (U_k) \subseteq p^k \, \zp$.\\
For the reverse way, take the restriction of the $\exp$ map on $p^k \, \zp$. Let $p^k x$ be any element of $p^k \, \zp$. Then,
$$\exp (p^k x) = 1 + p^k x + \dfrac{(p^k x)^2}{2!} + \cdots + \dfrac{(p^k x)^n}{n!} + \cdots$$
Then, R.H.S. of the above equation belongs to $1 + p^k \, \zp$. Therefore, $\exp (p^k \, \zp) \subseteq 1 + p^k \, \zp$.
Since, $\log: \uone \rightarrow p \, \zp$ and $\exp: p \, \zp \rightarrow \uone$ are one-one maps, there restrictions are also one-one. 
Hence, $\log_{\mid_{U_k}}: U_k \rightarrow p^k \, \zp$ and $\exp_{\mid_{p^k \, \zp}}: p^k \, \zp \rightarrow U_k$ 
are bijections. These maps are also group homomorphisms follows from $\exp$ and $\log$ are homomorphisms.   
\end{proof}

\begin{corollary}
\begin{eqnarray*}
\phi_k : U_1 &\longrightarrow& U_k\\
 \alpha &\mapsto& {\alpha}^{p^{k - 1}}
\end{eqnarray*}
is an isomorphism of groups.
\end{corollary}
\begin{proof}
If $\alpha = 1 + px$ then using binomial theorem, $(1 + px)^{p^{k - 1}} = 1 + p^{k - 1} px + \cdots + (px)^{p^{k - 1}}$. This element belongs to $U_k$. Therefore, the map $\phi_k$
is well defined.\\
From the above theorem, $(U_1, *) \cong (p \, \zp, +)$ and $(U_k, *) \cong (p^k \, \zp, +)$. Therefore, we have following commutative diagram:
$$\xymatrix{
U_1 \ar[d] \, \ar[r] & U_k\ar[d]\\
p \, \zp \ar[r] & p^k \, \zp }$$
$\phi_k$ induces the map $\psi_k : p \, \zp \rightarrow p^k \, \zp$ given by $px \mapsto p^k x$.  
If $\psi_k (x) = p^k x = 0$ then $x = 0$. Therefore, $\psi_k$ is one-one. Hence, $\phi_k$ is also one-one.
Clearly, $\psi_k$ is an onto map. Therefore in the above commutative diagram all four groups are isomorphic to each other.
\end{proof}

\indent Now we will discuss the maps $x \mapsto x^k$ on $\qp ^*$ for $k \geq 1.$ The case $k = 2$ is important in number theory,
in connection with the theory of quadratic forms. \\
\indent Let $\qp$ be a field of $p$-adic numbers. Recall that, every non-zero element of $\qp$ can be uniquely written as $p^n u$ where $n \- \Z$ and $u \- U$ is a unit in $\zp.$
From the theorem \ref{7}, we know the structure of $U,$ the group of units of $\zp.$ From the theorem \ref{20},  we get (by the $\log$ map) $U_1 = 1 + p \, \zp \cong p \,\zp$ 
But $(p \, \zp +)$ is isomorphic to $(\zp, +)$ by the map $px \mapsto x.$ For $p =2$, $U_1 \cong \Z_2 \times U_2$ and by the $\log$ map $(U_2, * ) \cong (2^2 \Z_{(2)}, +).$ 
Also $(2^2 \Z_{(2)}, +) \cong (\Z_{(2)}, +)$ (by the map $2^2 x \mapsto x$.) Therefore, structure of  $\qp ^* = \qp / \{0\}$ can be described as follows:
\begin{eqnarray*}
      \qp ^* &\cong& \Z \times U\\
             &\cong& \Z \times \Z_{p - 1} \times \zp  \hspace*{2em} p = \,odd \, prime.\\
 \Q_{(2)} ^* &\cong& \Z \times \Z_2 \times \Z_{(2)}  \hspace*{3em} p = 2.
\end{eqnarray*}

Let $k$ be any natural number. Define a map 
\begin{eqnarray*}
 \phi_k : \qp ^* &\rightarrow& \qp ^*\\
               x &\mapsto& x^k
\end{eqnarray*}
We find the structure of $\qp ^* / {\qp ^*}^k$. For this we will use the structure of $\qp ^*$ stated above. We first study $\phi_k$ on $\Z$, $\Z_{p - 1}$, $\zp.$
Under the map $\phi_k$ we will find the image, kernel, cokernel for the maps on $\Z$, $\Z_{p - 1}$, $\zp.$ Calculations involved in it are trivial hence we omit them.

\begin{enumerate}
\item \begin{eqnarray*}
 \phi_k ^1 : \Z &\rightarrow& \Z \\
               x &\mapsto& kx
\end{eqnarray*}
Then, \begin{enumerate}
\item  $Im \, \phi_k ^1 = k \Z $  
\item $ker \, \phi_k ^1 = \{0\}$ 
\item $\Z / {Im \, \phi_k ^1} = \Z_k$
\item $\Z / {ker \, \phi_k ^1} = \Z$ 
\end{enumerate}

\item  For $p$ is an odd prime.
 \begin{eqnarray*}
 \phi_k ^2 : \Z_{p - 1} &\rightarrow& \Z_{p - 1} \\
                  x &\mapsto& kx
\end{eqnarray*}                 
Then, \begin{enumerate}
\item $Im \, \phi_k ^2 =  \Z_{d_1} $ where $d_1 = \dfrac{p - 1}{(p - 1, k)}$
\item $ker \, \phi_k ^2 = \Z_{d_2}$ where $d_2 = (p - 1, k)$ 
\item $\Z_{p - 1} / {ker \, \phi_k ^2} = \Z_{d_1}$ 
\item $\Z_{p - 1} / {Im \, \phi_k ^2} = \Z_{d_2}$ 
\end{enumerate}

For $p = 2$
 \begin{eqnarray*}
 \phi_k ^2 : \Z_2 &\rightarrow& \Z_2 \\
                  x &\mapsto& kx
\end{eqnarray*}                 
Then, \begin{enumerate}
\item $Im \, \phi_k ^2 = \{0\}$ if $k$ is even.\\
      $Im \, \phi_k ^2 = \Z_2$ if $k$ is odd.
\item $ker \, \phi_k ^2 = \Z_2$ if $k$ is even.\\
      $ker \, \phi_k ^2 = \{0\}$ if $k$ is odd.
\item $\Z_2 / {ker \, \phi_k ^2} = \{0\}$ if $k$ is even.\\
      $\Z_2 / {ker \, \phi_k ^2} = \Z_2$ if $k$ is odd.
\item $\Z_2 / {Im \, \phi_k ^2} = \Z_2$  if $k$ is even.\\
      $\Z_2 / {Im \, \phi_k ^2} = \{0\}$  if $k$ is odd.
\end{enumerate}

\item      \begin{eqnarray*}
 \phi_k ^3 : \zp &\rightarrow& \zp \\
                  x &\mapsto& kx
\end{eqnarray*}
i) If $(k, p) = 1,$ then $k$ is a unit in $\zp.$ Therefore $kx = 0$ implies $x = 0.$ 
\begin{enumerate}
\item $Im \, \phi_k ^3 =  \zp$ 
\item $ ker \, \phi_k ^3 = \{0\}$ 
\item $\zp / {Im \, \phi_k ^1} = \{0\} $ 
\item $\zp / {ker \, \phi_k ^3} = \zp $ 
\end{enumerate}

ii) Suppose $(k, p) \neq 1,$ then $k = p^m l$ where $p \nmid l.$ 
Then, $kx = o \Rightarrow p^m l x = 0 \Rightarrow p^m x = 0 \Rightarrow x = 0.$
 \begin{enumerate}
\item $Im \, \phi_k ^3 =  p^m \,\zp$ 
\item $ ker \, \phi_k ^3 = \{0\}$ 
\item $\zp / {Im \, \phi_k ^1} = \Z_{p^m} $ 
\item $\zp / {ker \, \phi_k ^3} = \zp $ 
\end{enumerate} 
\end{enumerate}
\newpage
We summarize above three points in the following table. 
Notation in the table, \\
$i) \,d_1 = \dfrac{p - 1}{(p - 1, k)}$,  $ii) \, d_2 = (p - 1, k)$,\\
$iii) \, m$ is such that, $k = p^m l$ and $p \nmid l.$\\
In $\phi_k ^3$, $p$ can be even or odd prime. 
\begin{table}[h]
 \begin{tabular}{|c|c|c|c|c|}
\hline
$\phi_k ^i : X \rightarrow Y$, $x \mapsto kx$ & $Im$ & $ker$ & $X / ker$ & $Y / Im$\\
 \hline
\hline
$\phi_k ^1 : \Z \rightarrow \Z$ & $k \, \Z$ & $\{0\}$ & $\Z$ & $\Z_k$\\
\hline
$\phi_k ^2 : \Z_{p - 1} \rightarrow \Z_{p - 1}$, $p$ odd prime & $\Z_{d_1}$ & $\Z_{d_2}$ & $\Z_{d_1}$ & $\Z_{d_2}$\\
\hline
$\phi_k ^2 : \Z_p \rightarrow \Z_p$, $p = 2$, $k$ even & $\{0\}$ & $\Z_2$ & $\{0\}$ & $\Z_2$\\
\hline
$\phi_k ^2 : \Z_p \rightarrow \Z_p$, $p = 2$, $k$ odd & $\Z_2$ & $\{0\}$ & $\Z_2$ & $\{0\}$\\
\hline
$\phi_k ^3 : \zp \rightarrow \zp$, $(k, p) =1$ & $\zp$ & $\{0\}$ & $\zp$ & $\{0\}$ \\
\hline
$\phi_k ^3 : \zp \rightarrow \zp$, $(k, p) \neq 1$ & $p^m \, \zp$ & $\{0\}$ & $\zp$ & $\Z_{p^m}$ \\
\hline
 \end{tabular}
\end{table}

Combining all the results in the above table , we state following theorem:
\begin{theorem} \label{23}
 Define the map 
\begin{eqnarray*}
 \phi_k : \qp ^* &\rightarrow& \qp ^*\\
               x &\mapsto& x^k
\end{eqnarray*}
Then, for an odd prime $p$,
\begin{enumerate} 
\item $Im \, \phi_k  \cong k \, \Z \times \Z_{d_1} \times \zp $ if $(k, p) = 1.$\\
      $Im \, \phi_k  \cong k \, \Z \times \Z_{d_1} \times p^m \zp $ if $(k, p) \neq 1.$
\item $ker \, \phi_k  \cong \{0\} \times \Z_{d_2} \times \{0\} \cong \Z_{d_2}.$ 
\item $\qp ^* / {ker \, \phi_k} \cong \Z \times \Z_{d_1} \times \zp.$ 
\item $\qp ^* / {Im \, \phi_k} \cong \Z_k \times \Z_{d_2} \times \{0\}$ if $(k, p) = 1.$\\
      $\qp ^* / {Im \, \phi_k} \cong \Z_k \times \Z_{d_2} \times \Z_{p^m}$ if $(k, p) \neq 1.$
\end{enumerate}
 For $p = 2$ and $k$ even,
\begin{enumerate} 
\item $Im \, \phi_k  \cong k \, \Z \times \{0\} \times \zp $ 
\item $ker \, \phi_k  \cong \{0\} \times \Z_{d_2} \times \{0\} \cong \Z_{d_2}.$ 
\item $\qp ^* / {ker \, \phi_k} \cong \Z \times \{0\} \times \zp.$ 
\item $\qp ^* / {Im \, \phi_k} \cong \Z_k \times \Z_2 \times \Z_{2^m}$ if $k = 2^m l$ and $2 \nmid l.$
\end{enumerate}
 For $p = 2$ and $k$ odd,
\begin{enumerate} 
\item $Im \, \phi_k  \cong k \, \Z \times \Z_2 \times \zp.$
\item $ker \, \phi_k  \cong \{0\} \times \{0\} \times \{0\} \cong \{0\}.$ 
\item $\qp ^* / {ker \, \phi_k} \cong \Z \times \Z_2 \times \zp.$ 
\item $\qp ^* / {Im \, \phi_k} \cong \Z_k \times \{0\} \times \{0\}.$ 
\end{enumerate}
\end{theorem}

From the above theorem, let us note the structure of $\qp ^* / {\qp ^*}^k$ explicitly:\\

\begin{theorem} Let  $k = p^m \, l$ where $p \nmid l.$\\
For $p$ an odd prime. Let $d_2 = (p - 1, k)$. 
\begin{enumerate}
\item If $(k, p) = 1$ then $\qp ^* / {\qp ^*}^k \cong \Z_k \times \Z_{d_2}.$
\item If $(k, p) \neq 1$ then $\qp ^* / {\qp ^*}^k \cong \Z_k \times \Z_{d_2} \times \Z_{p^m}$ 
\end{enumerate}
If $p = 2$ then,
\begin{enumerate}
\item for $k$ even,  $\qp ^* / {\qp ^*}^k \cong \Z_k \times \Z_2 \times \Z_{p^m}$
\item for $k$ odd,   $\qp ^* / {\qp ^*}^k \cong \Z_k$ 
\end{enumerate}
 
\end{theorem}

\begin{corollary} Let $k = 2$ in the above theorem. Then\\
$\qp ^* / {\qp ^*}^2 \, \cong \, \Z_2 \times \Z_2$ if $p$ is an odd prime.\\
$\Q_{(2)} ^* / {\Q_{(2)} ^*}^2 \, \cong \, \Z_2 \times \Z_2 \times \Z_2 .$
\end{corollary}

This special case $k = 2$, is an important statement of {\it{quadratic forms}} in number theory, cf. \cite{S}.

\chapter{$p$-ADIC SOLUTION OF $a^x = b$}

\indent \textsc{THIS} chapter is devoted to the main problem of our study and some important examples of our theorems.
The equation $a^x = b$ can be considered in different systems, and accordingly the solution space will also change. More explicitly, for $a, b$ positive real numbers,
solution is given by $x \, \log a = \log b$, hence $x = \log b / \log a$ if $a \neq 1.$\\
\indent Our aim is to find the solution of $a^x \eq b \, (mod \, p^n)$ as $n$ varies in natural numbers and study the distribution of these solutions, 
for fixed integers $a$, $b$ and prime $p$. If $x_n$ denotes solution of this congruence mod $p^n$, 
we show that, as $n$ varies $x_n$ takes the form of $p$-adic integers and the sequence $x_n$ is convergent in the $p$-adic integers.
If $\{x_n\} \rightarrow x_o$ in $p$-adic integers, we say that, $x_o$ is the solution of $a^x = b$ in $\zp$.\\
\indent We will start with $a$ and $b$ natural numbers and try to extend the domain of $a$, $b$ to wider subset of $p$-adic integers.
Throughout, we assume that, $p$ is a prime number. $a$ and $b$ whenever integers are coprime to $p$ and $a \neq 1$. We divide the work of 
finding solutions of $a^x = b$ in 3 sections:\\
\noindent 1. $a$ and $b$ are integers both are coprime to $p$ and $a \neq 1$ \\
2.  $a$ and $b$ are $p$-adic integers both are elements of $U_1$ and $a \neq 1$ \\
3.  $a$ and $b$ are units in $p$-adic integers and $a$ is not congruent to 1 mod $p$.

Since $p$-adic integers $(\zp)$ is defined as projective limit of $\Z_{p^n}$, which means that, $\zp$ is a collective way of thinking all $\Z_{p^n}$'s. Hence, \\
$a^x = b$ in $p$-adic integers if and only if $a^x \eq b \, (mod \, p^n)$ for every $n$.\\

\section{Solutions of $a^x \eq b \, (mod \, p^n)$ as $n \rightarrow \infty$\\ for $a, b \- \Z$ and $(a, p) = (b, p) = 1$} \label{E}

In this section we prove that if the solutions of the congruence $a^{x_n} \eq b \, (mod \, p^n)$ exists for all $n,$ then $x_n$'s take the form of $p$-adic integers.

\begin{proposition} \label{24}
Let $p$ be an odd prime. If $a \neq 1$, $(a, p) = 1$ and order of $a$ in $\Z_p ^*$ is $x_o$ then order of $a$ in $\n$ is of the form $x_n = x_o p^{k_n}$ where
 $k_n \rightarrow \infty$ as $n \rightarrow \infty.$
\end{proposition}

\begin{proof} Let $p$ be an odd prime.\\
\underline{Case 1:} Let $a \eq 1 \, (mod \, p).$ We know that $\n$ is an internal direct product of its Sylow-$q$ subgroups, $q \mid {p - 1}$, and the group generated by $1 + p.$ 
  Since $a \eq 1 \, (mod \, p)$, $a$ must belong to $< 1 + p >.$ If $a$ is of the form $1 + p^k l$ where $p \nmid l$ then order of $a$ in $\n$ is $p^{n - k}.$
Therefore in this case, we get $x_n = p^{n - k}.$ Hence $k_n = n - k$ which goes to infinity as $n \rightarrow \infty.$ \\
\underline{Case 2:}  Let $a \neq 1 \, (mod \, p).$ Write $a = a_1 * a_2$ where $a_1 \- \Z_{p - 1}$ and $a_2 \- <\, 1 + p \,>.$ 
Then order of $a_1$ in $\Z_{p - 1}$ is same as order of $a$ in $\Z_p ^*$ namely $x_o.$  By the earlier case order of $a_2$ is $p^{k_n}$ and 
$k_n$ goes to infinity as $n \rightarrow \infty.$ Since order of $a_1$ is coprime to $p$, we get order of $a = a_1 * a_2$ in $\n$ is $x_o p^{k_n}.$\\
 \end{proof}

\begin{proposition} \label{13} 
 Suppose $a$ is coprime to $p$ and order of $a$ in $\Z_p ^*$ is $x_o$. Then \\
i) there exists a largest natural number $k$ such that, order of $a$ in $\Z_p ^*$, $\Z_{p^2} ^*$, $\cdots, \Z_{p^k} ^* $ is fixed number $x_o.$\\
ii) Order of $a$ in $\n$ is $x_o \, p^{n - k}$ for all $n \geq k$.
\end{proposition}

\begin{proof} Let $p$ be an odd prime.\\
\noindent i)  Let order of $a$ in $\Z_p ^*$ be $x_o.$ From the proposition \ref{24}, order of $a$ in $\n$ cannot be same for all $n.$ Therefore there exists a smallest natural number 
 $k + 1$ such that order of $a$ in $\Z_{p^{k + 1}} ^*$ is \underline{not} $x_o.$ This implies order of $a$ in  $\Z_p ^*$, $\Z_{p^2} ^*$, $\cdots, \Z_{p^k} ^* $ is $x_o.$
This $k$ itself is the largest natural number $k$ such that, order of $a$ in $\Z_p ^*$, $\Z_{p^2} ^*$, $\cdots, \Z_{p^k} ^* $ is fixed number $x_o.$ 

\noindent ii) As order of $a$ in $\Z_p ^*$, $\Z_{p^2} ^*$, $\cdots, \Z_{p^k} ^* $ is $x_o$, we get $a^{x_o} = 1 + p^k l$. \\
Order of $a$ in $\Z_{p^{k + 1}}$ is not $x_o$, therefore $p$ does not divide $l$. Now to prove that, order of $a$ in $\n$ is $x_o \, p^{n - k}$ for all $n \geq k$, we prove:\\
1. $a^{x_o \, p^{n - k}} \eq 1 \, (mod \, p^n)$ and \\
2. $a^{x_o \, p^{n - k - 1}} \neq 1 \, (mod \, p^n)$ 
 \begin{eqnarray*}
    a^{x_o \, p^{n - k}} &=& (1 + p^k l)^{p^{n - k}}\\
                      &=& 1 + p^{n - k} p^k l + \dfrac{p^{n - k} (p^{n - k} - 1)}{2} p^{2 k} l^2 + \cdots + (p^k l)^{p^{n - k}}\\
                      &=& 1 + p^n l + \dfrac{p^{n + k} (p^{n - k} - 1)}{2} l^2 + \cdots + (p^k l)^{p^{n - k}}\\
                      &\eq& 1 \, (mod \, p^n)
 \end{eqnarray*}
Also,
 \begin{eqnarray*}
    a^{x_o \, p^{n - k - 1}} &\eq& 1 + p^{n - 1} l  \, (mod \, p^n)
 \end{eqnarray*}
Since, $p$ does not divide $l$, we get $    a^{x_o \, p^{n - k - 1}} \neq 1 \, (mod \, p^n)$.\\
Therefore, order of $a$ in $\n$ is $p^{n - k}$.

For $p = 2$, $a$ must be an odd integer. Therefore $x_o = 1.$ If $a = 1 + 2^k \, l$ and $2 \nmid l$, we get $a \eq 1 \, (mod \, 2^n)$ for $1 \leq n \leq k.$ 
Therefore order of $a$ mod $2^n$ is 1 for $1 \leq n \leq k,$ For $n \geq k,$ above proof by binomial theorem works, after replacing $p$ by 2.
\end{proof}

\begin{example} (Gauss) Let $p = 29$ and $a = 14$. Then order of $14$ in $\Z_{29} ^*$ and in  $\Z_{29^2} ^*$ is $28$. i.e. $28$ is the smallest positive integer such that,
$$ 14^{28} \eq 1 \, (mod \, 29)$$
$$ 14^{28} \eq 1 \, (mod \, 29^2)$$
But $14^{28} \eq 14298 \, (mod \, 29^3)$. From the above theorem, order of $14 \, (mod \, 29^n)$ is $28 * 29^{n - 2}$ for all $n \geq 2$.
\end{example}

\begin{remark} Existence of largest number $k$ in part (i) of proposition \ref{13} is a special property of integers. Instead if we take $a$ as a Teichm\"{u}ller unit
not equal to 1, having order $x_o$ then order of $a \, mod \, p^n$ is $x_o$ for all $n \geq 1.$ Hence $k$ is $\infty$ in this case.\\
 From proposition \ref{13}, we observe that, as $n$ varies in natural nos. order of $a$ in $\n$ can not be of the form
 $x_o, x_o, x_o p, x_o p, x_o p^2, x_o p^3, x_o p^3, \cdots$
\end{remark}

\begin{theorem} \label{25}
Let $p$ be a prime. Let $a$, $b$ be integers coprime to a prime $p$ and $a \neq 1.$ Let $x_n$ denote the smallest positive solution of the congruence
 $a^x \eq b \, (mod \, p^n)$. If $x_n$ exists for all $n \geq 1$ then the sequence $\{x_n\}$ converges in $p$-adic integers.
\end{theorem}

\begin{proof}
We give a combine proof for $p = 2$ and $p$ an odd prime.
Let $x_n$ and $x_{n + 1}$ be the smallest positive integers such that,
 $$ a^{x_n} \eq b \, (mod \, p^n)$$
 $$ a^{x_{n + 1}} \eq b \, (mod \, p^{n+1})$$
Then, clearly, $x_{n + 1} \geq x_n$. 
$$ a^{x_{n + 1}} - a^{x_n} \eq 0 \, (mod \, p^n)$$
$$ a^{x_n} \, (a^{x_{n + 1} - x_n} - 1) \eq 0 \, (mod \, p^n)$$
Since $a$ and $p$ are coprime, $p^n \nmid a^{x_n}$. Therefore,
$$a^{x_{n + 1} - x_n} - 1 \eq 0 \, (mod \, p^n)$$
$$a^{x_{n + 1} - x_n} \eq 1 \, (mod \, p^n)$$
Therefore, order of $a \, (mod \, p^n)$ divides $x_{n + 1} - x_n.$

We will show that, as $n \rightarrow \infty$ powers of $p$ dividing $x_{n + 1} - x_n$ also goes to infinity. Suppose $k$ is the largest integer such that  order of $a$ is same in
$\Z_p ^*$, $\Z_{p^2} ^*,\cdots, \Z_{p^k} ^*$ say $x_o$. Then $1 < x_o \leq p - 1$. By proposition \ref{13}, order of $a$ in $\n$ is $x_o p^{n - k}$ for all $n \geq k.$
We have $x_o p^{n - k}$ divides $x_{n + 1} - x_n.$ For $n \geq 1$ we have
\begin{eqnarray*}
 x_{k + n} &=& x_{k + n - 1} + l_{k + n - 1} x_o p^{k + n - 1}\\
           &=& x_{k + n - 2} + l_{k + n - 2} x_o p^{k + n - 2 } + l_{k + n - 1} x_o p^{k + n - 1}\\
           &\vdots&\\
           &=& x_k + l_k x_o p + l_{k + 1} x_o p^2 + \cdots + l_{k + n - 1} x_o p^{k + n - 1}
\end{eqnarray*}
Write the $p$-adic expansion of all $l_i$'s then we get, as $n \rightarrow \infty$, powers of $p$ dividing $x_n - x_{n - 1}$ also goes to infinity, for $n \geq k.$
Hence the sequence $\{x_n\}$ converges in $p$-adic integers. 
\end{proof}

 \begin{example} \label{9}
 Consider the crucial example for the above discussion\\
 \framebox(130,20)[c]{$(-3)^x \eq 5 \, (mod \, 2^n)$} 

So, $p = 2$, $a = -3$, $b = 5$. We look for the (2=p)-adic solution of the equation $(-3)^x = 5$.\\
$-3, 5 \in \t.$ We will prove that, $<-3> = <5>$ in $\t$. 
\newpage
order 5 = $2^{n-2}$  \, \,  proved earlier.\\
order $-3$ = $2^{n-2}$\\
Proof: Write $-3$ as $1 - 2^2.$ Since, $5 = 1 + 2^2$, binomial expansion of $(1 - 2^2)^r$ and $(1 + 2^2)^r$ differs only in sign at even places.\\
As \begin{eqnarray*}
                5^{\w} &=& (1 + 2^2)^{\w}\\
                       &\eq& 1 \, (mod \, 2^n)\\
 \Rightarrow (-3)^{\w} &=& (1 - 2^2)^{\w}\\
                       &\eq& 1 \, (mod \, 2^n)
   \end{eqnarray*}
 and
  \begin{eqnarray*}
                5^{\v} &=& (1 + 2^2)^{\v}\\
                       &\eq& 1 + 2^{n - 1} \, (mod \, 2^n)\\
 \Rightarrow (-3)^{\v} &=& 1 - 2^{n-1} \, (mod \, 2^n)\\
                       &\neq& 1   \, (mod \, 2^n)
   \end{eqnarray*}
   
Therefore, order of $-3$ = $\w$.

We now prove that, $<-3> = <5>$ in $\t$. 

Proof: 1. Order of $-3$ = $\w$\\
$(-3)^{\w} \eq 1 \, (mod \, 2^n)$\\
$\Rightarrow ((-3)(-1))^{\w} \eq 1 \, (mod \, 2^n)$\\
and $(-3)^{\v} \eq 1 - 2^{n - 1} \, (mod \, 2^n)$\\
$\Rightarrow 3^{\v} \eq (-1)^{\v} (1 - 2^{n - 1}) = 1 - 2^{n - 1} \neq 1 \, (mod \, 2^n)$.

Therefore, \underline{order of 3 is $\w$.}

2. We see that, \underline{$<3> \neq <5>$ in $\t$.}\\
because, if $3 \in <5>$ then $3 = 5^k$ for some $k$.\\
$5 \eq 1 \, (mod \, 4) \Rightarrow 5^k \eq 1 \, (mod \, 4)$\\
but $3 \eq -1 \, (mod \, 4)$ and $-1 \neq 1 \, (mod \, 4)$\\
Therefore, $3 \notin <5>$.

Therefore, we have $< 3 >$, $< -3 >$ and $< 5 >$ are $3$ subgroups of order $\w$ in $\t$
and $<3> \neq <5>.$

From the structure of $\t$,\\
$\t \cong \Z_2 \times \Z_{\w}$, there are exactly $2$ distinct subgroups of order $2^{n - 2}$ for $n \geq 3.$

Hence we must have, $<-3> = <3>$ or  $<-3> = <5>$, and both cannot be true. 

If $<3> = <-3>$, then $3 \in <-3> \Rightarrow 3 = (-3)^k$ for some $k.$\\
$-3 \eq 1 \, (mod \, 4) \Rightarrow (-3)^k \eq 1 \, (mod \, 4)$\\
But $3 \eq -1 \, (mod \, 4)$. \\
$\Rightarrow 3 \notin <-3>.$

Therefore, $<-3> = <5>$ must be true. $~~~~~~~~~~~~ \square $\\

 We have obtained the actual solution of $(-3)^{x_n} \eq 5 \, (mod \, 2^n)$ \\ 
for $1 \leq n \leq 10 $ using MATHEMATICA.

 \begin{table}[h]
  \begin{tabular}{|c||c|c|c|c|c|c|c|c|c|c|}
\hline
 $n$ & 1 & 2 & 3 & 4 & 5 & 6 & 7 & 8 & 9 & 10 \\ [.5 ex]
 \hline
 $x_n$ & 1 & 1 & 1 & 3 & 3 & 11 & 11 & 11 & 11 & 11 \\[1 ex]
\hline 
  \end{tabular}

 \end{table}

 Hence, 2-adic expansion of \\
$x_1 = x_2 = x_3 = 1 = 1* 2^0$\\
$x_4 = x_5 = 3 = 1 + 1* 2$\\
$x_6 = x_7 = x_8 = x_9 = x_{10} = 1 + 1* 2 + 2^3.$

Using another software SAGE we can directly find p-adic solution to above congruence
$(-3)^x \eq 5 \, (mod \, 2^n)$. Solution is given by\\
$x = 1+ 2+ 2^3+ 2^8+ 2^{10}+ 2^{11}+ 2^{12}+ 2^{13}+ o(2^{18})$
\end{example}

\begin{example} Let $p = 5$, $a = 2$ and $b = 1$. We look for the non-zero solution of $2^{x_n} \eq 1 \, (mod \, 5^n)$.
i.e. we find the order of 2 in the multiplicative group $\Z_{5^n} ^*$. \\
In the following table $x_n$ denotes the order of 2 in $\Z_{5^n} ^*$. These calculations are done with the help of `SAGE'. \\
\begin{table}[h]
\begin{tabular}{|c||c|c|c|c|c|c|c|c|c|c|}
\hline
 $n$ & 1 & 2 & 3 & 4 & 5 & 6 & 7 & 8 & 9 & 10\\
\hline
$x_n$ & $4$ & $4 * 5$ & $4 * 5^2$ & $4 * 5^3$ & $4 *5^4 $& $4 * 5^5$ & $4 * 5^6$ & $4 * 5^7$ & $4 * 5^8$ & $4 * 5^9$\\
\hline
\end{tabular}
\end{table}

In general, order of 2 in $\Z_{5^n} ^*$ is $4 * 5^{n - 1}$.\\
We prove this using binomial theorem: For $n \geq 1$,\\
Step 1: $2^{4 * 5^{n - 1}} \eq 1 \, (mod \, 5^n)$\\
Step 2: $2^{4 * 5^{n - 2}} \neq 1 \, (mod \, 5^n)$
\begin{proof}
 $2^{4 * 5^{n - 1}} = (2^4)^{5^{n - 1}} = 16^{5^{n - 1}}$
\begin{eqnarray*}
 16^{5^{n - 1}} &=& (1 + 3 * 5)^{5^{n - 1}}\\
                &=& 1 + 3 * 5 * 5^{n - 1} + \dfrac{5^{n - 1} (5^{n - 1} - 1)}{2} (3 * 5)^2 + \cdots + (3 * 5)^{5^{n - 1}}\\
                &=& 1 + 3 * 5^n + \dfrac{5^{n - 1} - 1}{2} 3 * 5^n + \cdots + (3 * 5)^{5^{n - 1}}\\
                &\eq& 1 \, (mod \, 5^n)
\end{eqnarray*}
\begin{eqnarray*}
 16^{5^{n - 2}} &=& (1 + 3 * 5)^{5^{n - 2}}\\
                &=& 1 + 3 * 5 * 5^{n - 2} + \dfrac{5^{n - 2} (5^{n - 2} - 1)}{2} (3 * 5)^2 + \cdots + (3 * 5)^{5^{n - 2}}\\
                &\eq& 1 + 3 * 5^{n - 1}  \, (mod \, 5^n)\\
                &\neq& 1 \, (mod \, 5^n)
\end{eqnarray*}

Therefore, $2^{4 * 5^{n - 2}} \neq 1 \, (mod \, 5^n)$ but $2^{4 * 5^{n - 1}} \eq 1 \, (mod \, 5^n)$.\\
Therefore order of 2 in $\Z_{5^n} ^*$ is $4 * 5^{n - 1}$.
\end{proof}
\end{example}
\newpage
\begin{example} Let $p = 5$ and $a = -2$, $b = 3$.\\
We have used `SAGE' to find first few coefficients of the $5$-adic solution of $(-2)^x = 3$.
In the following table, $x_n$ is the smallest positive integer satisfying the congruence $(-2)^{x_n} \eq 3 \, (mod \, 5^n)$.
\begin{table}[h]
\begin{tabular}{|c||c|c|}
\hline
$n$ & $x_n$ & 5-adic expansion of $x_n$\\
\hline
\hline
1 & 1 & 1\\
\hline
2 & 17 & $2 + 3 * 5$\\ 
\hline
3 & 57 & $2 + 5 + 2 * 5^2$\\
\hline
4 & 357 & $2 + 5 + 4 * 5^2 + 2 * 5^3$\\ 
\hline
5 & 1857 & $2 + 5 + 4 * 5^2 + 4 * 5^3 + 2 * 5^4$\\
\hline
6 & 14357 & $2 + 5 + 4 * 5^2 + 4 * 5^3 + 2 * 5^4 + 4 * 5^5$\\
\hline
7 & 14357 & $2 + 5 + 4 * 5^2 + 4 * 5^3 + 2 * 5^4 + 4 * 5^5$\\
\hline
8 & 201857 & $2 + 5 + 4 * 5^2 + 4 * 5^3 + 2 * 5^4 + 4 * 5^5 + 2 * 5^6 + 2 * 5^7$\\
\hline
9 & 1139357 & $2 + 5 + 4 * 5^2 + 4 * 5^3 + 2 * 5^4 + 4 * 5^5 + 2 * 5^6 + 4 * 5^7 + 2 * 5^8$\\
\hline
10 & 5826857 & $2 + 5 + 4 * 5^2 + 4 * 5^3 + 2 * 5^4 + 4 * 5^5 + 2 * 5^6 + 4 * 5^7 + 4 * 5^8 + 2 * 5^9$\\
\hline
\end{tabular}
\end{table}

We observe from the above calculations that, coefficient of $5^n$ becomes eventually constant. 
\end{example}

\section{$a, b \- U_1$ if $p$ odd prime\\ $a, b \- U_2$ if $p = 2$ } \label{F}

Recall $U_n$ denotes the set $1 + p^n \zp$.
Let $a, b$ be elements of $\uone $ and $a \neq 1$ and $p$ be a fixed prime.
As $p \nmid a$ and $p \nmid b$, we get $a \, (mod \, p^n)$, $b \, (mod \, p^n)$ both are units in the ring $\Z_{p^n}$.
We have seen the structure of $\n$: 
 $$\Z_{2^n}^* \cong \Z_2 \times \Z_{2^{n-2}} \, for \, n \geq 2.$$
$$\Z_{p^n}^* \cong \Z_{p-1} \times \Z_{p^{n-1}} \, for \, n \geq 1,$$
\noindent where p is an odd prime. Let $A_n$ and $B_n$ be the subgroups of $\n$ such that, 
$A_n$ is isomorphic to  $\Z_{p - 1} \, (or \, \Z_2)$ and $B_n$ is isomorphic to  $\Z_{p^{n-1}} \, (or \, \Z_{2^{n - 2}})$. 
Then $\n$ is a \underline{canonical internal direct product} of $A_n$ and $B_n$. 
Note that, for $p = 2$, this internal direct product is NOT canonical, since the subgroup of order 2 is not unique.
However, if we allow ourself to represent elements of $\Z_n$ by negative integers then $< -1 >$ is a canonical choice for subgroup isomorphic to $\Z_2.$
Then $\t$ is a canonical internal direct product of $< -1 >$ and $< 3 >. $\\
Taking projective limits of these isomorphisms, we get,
$$ \Z_{(2)} ^* \cong \, \Z_2 \times \Z_{(2)}$$
$$ U = \zp ^* \cong \, \Z_{p - 1} \times \zp$$

We use the method of logarithm to find the solution of $a^x = b$ in $p$-adic integers.
We have already mentioned that,\underline{ if $p$ is an odd prime}, $\log$ of $p$-adic integers is defined from $\uone \rightarrow p \, \zp$.
Hence finding $\log$ of $a$ and $b$ make sense. We will prove the existence of $p$-adic solution later. [Adhoc proof: If $a^x = b$ then $x \log{a} = \log{b}$.]
Suppose $x= \displaystyle\x \- \zp$ denote the solution of $a^x = b$ in $\zp.$ 
Let $a= 1 + \displaystyle\a$ and $b = 1 + \displaystyle\b$ where $a_i, b_i \- \{0, 1, \cdots, p - 1\}$, be the $p$-adic expansions of $a$ and $b$ respectively.\\
$\log$ function is defined as follows:
$$\log : \uone \longrightarrow p \, \zp$$
$ \hspace{13em} 1 + p x \, \mapsto \, px - \dfrac{(px)^2}{2} + \dfrac{(px)^3}{3} - \cdots$\\
Since valuation of $\dfrac{(px)^n}{n}$ goes to infinity as $n \rightarrow \infty$, series on the R.H.S. is convergent in $\zp$.
 
We want to find $x_i \- \{0, 1, \cdots, p - 1\}$ such that,
\begin{center} $ x \, \log{a} = \log{b}$ \end{center} 
\begin{eqnarray*}
(\x) \, \, \log{(1 + \a)} = \, \log{(1 + \b)}\\
((x_o + x_1 p + x_2 p^2 + \cdots) ( \, \a - \dfrac{1}{2} {(\a)^2} + \dfrac{1}{3} (\a)^3 - \cdots)\\
= ( \, \b - \dfrac{1}{2} {(\b)^2} + \dfrac{1}{3} (\b)^3 - \cdots)
\end{eqnarray*}

We can find $x_i$ by comparing coefficients of $p^{i + 1}$ on both sides.\\
For example: 
Comparing coefficients of $p$ we get,

$a_1 x_o = b_1$\\
If $a_1 \neq 0$ then $a_1 ^{- 1}$ exists in $\Z_p$. So $x_o = a_1 ^{-1} b_1$.\\
If $a_1 = 0$ then $b_1$ has to be zero. Therefore $x_o$ can take any value from 0 to $p - 1$.

Comparing coefficients of $p^2$ we get,

$x_o a_2 + x_1 a_1 - \dfrac{1}{2} x_o a_1 ^2 = b_2 - \dfrac{1}{2} b_1 ^2 $\\
If $a_1 \neq 0$ then substituting value of $x_o$,\\
$a_1 ^{-1} b_1 a_2 + x_1 a_1 - \dfrac{1}{2} a_1 b_1 = b_2 - \dfrac{1}{2} b_1 ^2$\\
$x_1 = a_1 ^{-1} \, ( b_2 - \dfrac{1}{2} b_1 ^2 - a_1 ^{-1} b_1 a_2 + \dfrac{1}{2} a_1 b_1)$\\
If $a_1  = 0 \, (\Rightarrow b_1 = 0)$ then we get $x_o a_2 = b_2$ which do not give any condition on $x_1$.\\
Hence if $a_1 = 0 $ then $x_1$  can take any value from 0 to $p - 1$.\\
and so on.

Recall from definition \ref{29}, if $p = 2$ then $\log$ is defined from $1 + 2^2 \, \Z_{(2)}$ to $2^2 \, \Z_{(2)}.$ All the above calculations work similarly as in case of odd prime
$p$.

For our further results we recall some earlier notations: $\n$ is an internal direct product of $A_n$ and $B_n$. \\
If $p$ is an odd prime then $A_n \cong \Z_{p - 1}$ and $B_n \cong \Z_{\p}$ for $n \geq 1$\\
If $p = 2$ then $A_n \cong \Z_2$ and $B_n \cong \Z_{\w}$ for $n \geq 2.$

\begin{proposition} \label{1} 
i) Let $p$ be an odd prime. Let $a \neq 1$ be an element of $U_1$. Then index of group generated by $a \, (mod \, p^n)$ in $B_n (\cong \Z_{\p})$ is eventually constant.\\
ii) Let $p = 2$. Let $a \neq 1$ be an element of $U_2$. Then index of group generated by $a \, (mod \, 2^n)$ in $B_n (\cong \Z_{\w})$ is eventually constant.
\end{proposition}

\begin{proof} We give combine proof for both parts.
 Let $a = 1 + p^r \alpha$ be the $p$-adic expression of $a$ and $p \nmid \alpha$. 
We prove that, $a \, (mod \, p^n)$ as an element of $B_n$ is of order $p^{n - r}$ for all $n \geq r$.\\
i.e. to prove that, $$ a^{\nr} \eq 1 \, (mod \, p^n)$$
and $$a^{\nrone} \neq 1 \, (mod \, p^n)$$
If $n= r$, then $a \, (mod p^n) = 1$, so the above statement is clearly true.\\
For $n \geq r + 1$, we prove the claim by Binomial Theorem:
\begin{eqnarray*}
 a^{\nr} &=& (1 + p^r \alpha)^{\nr}\\
         &=& 1 + \nr p^r \alpha + \dfrac{\nr (\nr - 1)}{2!} p^{2r} {\alpha}^2 + \cdots + \nr (p^r \alpha)^{\nr - 1} + (p^r \alpha)^{\nr}\\
         &=& 1 \, (mod \, p^n)
\end{eqnarray*}
and 
\begin{eqnarray*}
 a^{\nrone} &=& (1 + p^r \alpha)^{\nrone}\\
         &=& 1 + \nrone p^r \alpha + \dfrac{\nrone (\nrone - 1)}{2!} p^{2r} {\alpha}^2 + \cdots + \nrone (p^r \alpha)^{\nrone - 1} + (p^r \alpha)^{\nrone}\\
         &=& 1 + p^{n - 1} \alpha \neq 1 \, (mod \, p^n)
\end{eqnarray*}
Therefore, $< a \, (mod \, p^n) > \, \subseteq B_n$ is of order $\nr$.\\
Therefore index of $< a \, (mod \, p^n)>$ in $B_n$ is equal to $\dfrac{p^{n - 1}}{\nr} = p^{r - 1}$ for all $n \geq r$. 
\end{proof}

\begin{proposition}
 If $a \neq 1$ is a unit in $\zp$ then index of $ < a \, (mod \, p^n) >$ in $A_n \times B_n$ is eventually constant.
\end{proposition}

\begin{proof} Case 1: Let $p$ be an odd prime.
 Every unit in $\zp$ can be uniquely written as a product of a Teichm\"{u}ller unit and an element of $\uone$. Let $a = a_1 a_2$ where $a_1$ is Teichm\"{u}ller unit 
and $a_2 \- \uone$. Since $\zp ^*$ is isomorphic to $\Z_{p - 1} \times (\uone)$, group generated by $a_1$ is a subgroup of $\Z_{p - 1}$.
By the construction of Teichm\"{u}ller unit, for all $n \geq 1$, order of $a_1 \,(mod \, p^n)$ is equal to the order of $a_1$. 
Therefore, index of group generated by $a_1 \,(mod \, p^n)$ in $A_n$ is constant.\\
By the \ref{1}, index of subgroup generated by $a_2 \, (mod \, p^n)$ in $B_n$ is eventually constant.\\
Therefore, index of $a \, (mod \, p^n)$ in $A_n \times B_n$ ($\cong \Z_{p - 1} \times \Z_{\p}$) is eventually constant.\\
Case 2: Let $p = 2.$
 Every unit in $\Z_{(2)}$ can be uniquely written as a product of a Teichm\"{u}ller unit and an element of $U_2$. Let $a = a_1 a_2$ where $a_1$ is Teichm\"{u}ller unit 
and $a_2 \- U_2$. Then $a_1$ can be either 1 or $-1$. If $a_1 = 1$ then $a = a_2$ and by previous proposition we are done. If $a_1 = -1$ then $a_1 \, (mod \, 2^n)$ is of order 2.
Therefore $< a_1 \, (mod \, 2^n) > = A_n.$ Hence index of $< a \, (mod \, 2^n)>$ in $A_n \times B_n$ is same as index of $< a_1 \, (mod \, 2^n) >$ in $B_n$, 
which is eventually constant by previous proposition.
\end{proof}

\begin{proposition} \label{2}
If $p$ is an odd prime, let $a, b \- U_1$, and if $p = 2$ let $a, b \- U_2$, and $a \neq 1.$
 Let $a = 1 + p^r \alpha$. If $b (mod \, p^n)$ belongs to the group generated by  $ a \, (mod \, p^n)$ for some $n_o \geq r$ 
then $b\, (mod \, p^n)$ belongs to group generated by $  a \, (mod \, p^n) $ for all $n \geq n_o$.
\end{proposition}

\begin{proof}
 From the above proposition $[B_n \, : \, < a \, (mod \, p^n) > ]$ and\\ $[B_n \, : \, < b \, (mod \, p^n)>]$ both are eventually constant.\\
 Therefore $[< a \, (mod \, p^n) > \, : \,  < b \, (mod \, p^n) > ]$ is also eventually constant.
\end{proof}

\begin{definition} 
Let $a = 1 + p^r \alpha$ for some $r \geq 1$ if $p$ is odd prime and $r \geq 2$ if $p = 2$, and $\alpha \- \zp$ such that $p$ does not divide $\alpha$.
Define depth of $a$ as $r$. We set depth of 1 as $\infty.$
\end{definition}

\begin{remark}
 Depth of $a$ is same as $p$-adic valuation of $a - 1.$
\end{remark}

\begin{proposition}  \label{3}
 Let $a, b \- U_1$ if $p$ is odd prime and $a, b \- U_2$ if $p = 2$. Let depth of $a$ be equal to $r$. 
Then $a^x = b$ has a solution in $p$-adic integers if and only if $b \, (mod p^n)$ belongs to group generated by
 $a \, (mod \, p^n)$ for some $n \geq r$.
\end{proposition}

\begin{proof} Following proof works for both cases of $p.$
Let $a = 1 + p^r \alpha$, $b = 1 + p^s \beta$ where $p$ does not divide both $\alpha$, $\beta$.
 Suppose  $b \, (mod p^n)$ belongs to group generated by $a \, (mod \, p^n)$ for some $n \geq r$.\\
Then by proposition \ref{2}, we get  $b \, (mod p^n)$ belongs to group generated by $a \, (mod \, p^n)$ for all $n \geq r$.
For $n \geq r $, if $$a^{x_n} \eq b \, (mod \, p^n)$$  
$$a^{x_{n + 1}} \eq b \, (mod \, p^{n + 1})$$
then we have already shown that $x_{n + 1} - x_n$ is divisible by $p^n$ as $n \rightarrow \infty$.
Hence, sequence $\{x_n\}$ is convergent in the $p$-adic integers. Limit of this sequence is the solution of $a^x = b$ in $p$-adic integers.\\
Conversely, Suppose $x_o \- \zp$ is a solution of $a^x = b$. Then, by the definition of $a^x = b$ we have $a^x \eq b \, (mod \, p^n)$ for all $n$. 
i.e. $b \, (mod \, p^n)$ belongs to the group generated by $a \, (mod \, p^n)$ 
\end{proof}

\begin{example} Let $p$ be an odd prime. If $a = 1 + p^2$ and $b = 1 + p$ then $a^x = b \, (mod \, p^n)$ has solution only when $n = 1.$
 Because, $1 + p$ is a generator of $\Z_{p^{n - 1}} \subset \Z_{p^n} ^*$.
Whereas $a = 1 + p^2$ is of order $p^{n - 2}$ in $\n$. Therefore, $< a > $ is a proper subgroup of $< b >$.
 Hence, $a^x = b \, (mod \, p^n)$ cannot have solution for $n > 1$.
\end{example}

 \begin{theorem}  \label{26}
  Let $a, b \- U_1$ if $p$ is odd prime and $a, b \- U_2$ if $p = 2$ and $a \neq 1$. If depth of $a \leq$ depth of $b$ then solution of $a^x = b$ exists as a $p$-adic integer and
 it is given by $x = \log \, b / \log \, a$.
\end{theorem}

\begin{proof} First part of the proof works for both odd and even primes $p.$
 Let $a, b \- U_1$ if $p$ is odd prime and $a, b \- U_2$ if $p = 2$. We can write $$a = 1 + p^r \alpha \, where \, p \nmid \alpha$$
$$b = 1 + p^s \beta \, where \, p \nmid \beta$$
For $1 \leq n \leq r$, we have $a \, (mod \, p^n) = 1 $ hence $o(a) = 1$ in $\n$ \\
For $ n > r$, we have $a \, (mod \, p^n) = 1 + p^r k $  where $ p \nmid k$.\\
Then $o(a) = p^{n - r}$ in $\n$.\\
It has a trivial proof by Binomial theorem:\\
Step 1: $( 1 + p^r k )^{p^{n - r}} \eq 1 \, (mod \, p^n)$\\
Step 2: $( 1 + p^r k )^{p^{n - r - 1}} \neq 1 \, (mod \, p^n)$\\
Similarly, for $1 \leq n \leq s$, order of $b \,(mod \,p^n)$ is 1 and\\
for $n > s$ order of $b \,(mod \,p^n)$ is $p^{n - s}.$

Let $x_n$ and $x_{n + 1}$ be the smallest positive integers such that,
 $$ a^{x_n} \eq b \, (mod \, p^n)$$
 $$ a^{x_{n + 1}} \eq b \, (mod \, p^{n+1})$$
Then, clearly, $x_{n + 1} \geq x_n$. Also, $0 < x_n \leq p^{n - 1}$ and $0 < x_{n + 1} \leq p^{n}$.\\
Therefore, $0 \leq x_{n + 1} - x_n \leq p^n.$
$$ a^{x_{n + 1}} - a^{x_n} \eq 0 \, (mod \, p^n)$$
$$ a^{x_n} \, (a^{x_{n + 1} - x_n} - 1) \eq 0 \, (mod \, p^n)$$
Since $a$ and $p$ are coprime, $p^n \nmid a^{x_n}$. Therefore,
$$a^{x_{n + 1} - x_n} - 1 \eq 0 \, (mod \, p^n)$$
$$a^{x_{n + 1} - x_n} \eq 1 \, (mod \, p^n)$$
Therefore, order of $a \, (mod \, p^n)$ divides $x_{n + 1} - x_n \h3 (*)$

We will show that, as $n \rightarrow \infty$ powers of $p$ dividing $x_{n + 1} - x_n$ also goes to infinity. To prove this we will make 2 cases:\\
$1.\, <a \, (mod \, p^n)> \,=\, B_n $ and $2. \, <a \, (mod \, p^n)> \, \subset \, B_n$.\\

\underline{ Let $p$ be an odd prime.}\\
\underline{Case 1:} Suppose $<a> = B_n \cong \Z_{\p}$. Therefore, order of $a \, (mod \, p^n)$ is $\p$. From $(*)$, we get $\p \mid x_{n + 1} - x_n$.
\begin{eqnarray*}
 x_{n + 1}  &=& x_n + k_{n - 1} \p \\
            &=& x_{n - 1} + k_{n - 2} \q + k_{n - 1} \p\\
	    &=& x_{n - 2} + k_{n - 3} p^{n - 3} + k_{n - 2} \q + k_{n - 1} \p\\
	    &\vdots&\\
	    &=& x_1 + k_1 p + k_2 p^2 + \cdots + k_{n - 2} \q + k_{n - 1} \p
\end{eqnarray*}
Where $k_i \- \{0, 1, \cdots , p - 1\}.$
Hence $x_{n + 1}$ is a polynomial in $p$ with coefficients in $\{0, 1, \cdots , p - 1\}.$

\underline{Case 2:} Suppose $<a \, (mod \, p^n)> \, \subset B_n$. Let $r_n =$ order of $a\, (mod \, p^n)$.
Recall, order of $a \, (mod \, p^n)$ divides $x_{n + 1} - x_n \h3 (*)$\\
Then we get $p^{r_n}$ divides $x_{n + 1} - x_n$.\\
Now, $0 \leq x_{n + 1} - x_n \leq p^n$ implies that $0 \leq \dfrac{x_{n + 1} - x_n}{p^{r_n}} \leq p^{n - r_n}$
\begin{eqnarray*}
 x_{n + 1}  &=& x_n + k_n p^{r_n} \\
            &=& x_{n - 1} + k_{n - 1} p^{r_{n - 1}} + k_n p^{r_n} \\
            &\vdots&\\
            &=& x_1 + k_1 p^{r_1} + k_2 p^{r_2} + \cdots + k_n p^{r_n} 
\end{eqnarray*}
Where $k_i$'s are such that, $0 \leq k_i \leq p^{n - r_n}$. Then write the $p$-adic expansion of each $k_i$.
Hence we can write $x_{n + 1}$ as a polynomial in $p$ with coefficients in $\{0, 1, \cdots , p - 1\}$

Let $p = 2$. Above proof can be mimicked. The only change is, $B_n$ is isomorphic to $\Z_{\w}.$ Therefore order of $a \, (mod \, 2^n)$ is $\w$ in case 1
and $2^{r_n}$, $1 \r_n \leq n - 2$ in case 2. Then further calculations are exactly same as in case of $p$ odd prime. \\
This proves the existence of the $p$-adic solution to $a^x = b.$\\
Further proof works for all primes $p$. Suppose depth of $a$ is $r$ and depth of $b$ is $s$ and $r \leq s.$ As $a, b \- U_1$, we get $\log \, a$ and $\log \, b$ make sense. 
  Let $a = 1 + p^r \alpha$ and   $b = 1 + p^s \beta$. Then 
\begin{eqnarray*}
 \dfrac{\log b}{\log a} &=& \dfrac {p^s \beta - \frac{(p^s \beta)^2}{2} + \cdots + (-1)^{n + 1} \frac{ (p^s \beta)^n }{n} + \cdots}
                                   { p^r \alpha - \frac{(p^r \alpha)^2}{2} + \cdots + (-1)^{n + 1} \frac{ (p^r \alpha)^n }{n} + \cdots }\\
                        &\vspace*{.2cm}&\\
                        &=& \dfrac{ p^{s - r} (p^r \beta - \cdots)} {\alpha - \frac{(p^r \alpha)^2}{2} + \cdots}
\end{eqnarray*}
On the R.H.S. of the above fraction, constant term of the denominator is non zero. Therefore denominator is a unit in $\zp.$ 
Hence $\log b / \log a$ is a $p$-adic integer.
\end{proof}

\begin{theorem} \label{27}
  Let $a, b \- U_1$ if $p$ is odd prime and $a, b \- U_2$ if $p = 2$ and $a \neq 1$. Then $a^x = b$ has unique solution in $p$-adic integers if and only if 
depth of $a \leq$ depth of $b$.
\end{theorem}

\begin{proof} Following proof works for all primes $p.$
  Let $a = 1 + p^r \alpha = 1 + p^r (a_r + a_{r + 1} p + \cdots)$, $a_r \neq 0$ and\\
  $b = 1 + p^s \beta = 1 + p^s (a_s + a_{s + 1} p + \cdots)$, $a_s \neq 0$ and $r \leq s$. \\
claim: There exists $n_o \geq r + 2$ such that $b \, (mod \, p^{n_o})$ belongs to the group generated by $a \, (mod \, p^{n_o})$.\\
Take  $n_o = r + s$. We will show that, this $n_o$ works! \\
We have, constant term of $\alpha = a_r + a_{r + 1} p + \cdots$, is $a_r$ which is non-zero, therefore $\bar{\alpha} = \alpha \,(mod \, p^{n_o})$ 
is a unit in the ring $\Z_{p^{n_o}}$.
We denote by ${\alpha}^{-1}$ the inverse of $\bar{\alpha}$ in $\Z_{p^{n_o}}$. Then,
 \begin{eqnarray*}
  a^{p^{s - r} \beta {\alpha}^{-1}} &=& (1 + p^r \alpha)^{p^{s - r} \beta {\alpha}^{-1}}\\
                                    &=& 1 + p^{s - r} \beta {\alpha}^{-1} p^r \alpha + \dfrac{p^{s - r} \beta {\alpha}^{-1} (p^{s - r} \beta {\alpha}^{-1} - 1)}{2} p^{2r}
                                      {\alpha}^2 + \cdots + (p^r \alpha)^{p^{s - r} \beta {\alpha}^{-1}}\\
 \end{eqnarray*}
Therefore, $b = a^{p^{s - r} \beta {\alpha}^{-1}}$ in the group $\Z_{p^{r + s}} ^*$. Also, $r + s \geq r + 2$. Since $s \geq 2$.
For $s = 1$, $r$ is also 1. Then order of $a$ and $b$ is same in $\Z_{p^n} ^*$ for all $n$. 
Therefore by Proposition \ref{2} and proposition \ref{3}, we get, there exists $p$-adic integer $x$ such that $a^x = b.$\\
Conversely, Suppose there exists $x \- \zp$ such that $a^x = b.$\\ Then $ a^x \eq b \, (mod \, p^n )$ for all $n$. Then clearly, depth of $a \leq$ depth of $b$.
\end{proof}

\begin{theorem}
 $p$-adic solution of $a^x = b$ is a unit in $\zp$ \iff depth of $a =$ depth of $b$.
\end{theorem}

\begin{proof}
 Suppose $p$-adic integer $x$ satisfying $a^x = b$ is a unit. Let $y = x^{-1}$. Then 
\begin{eqnarray*}
 (a^x)^y &=& b^y\\
     a &=& b^y\\
  depth \, of \, b &\leq& \, depth \, of \, a     
\end{eqnarray*}
Also, $a^x = b$ implies that, depth of $a \leq$ depth of $b$.\\
Therefore, depth of $a =$ depth of $b$.\\
Conversely, suppose, depth of $a =$ depth of $b$. Therefore there exists $p$-adic integers $x, y$ such that $a^x = b$ and $a = b^y$. 
Therefore, groups generated by $a$, $b$ in $\n$ are equal for all $n$. Hence the map $a \mapsto a^x$ from $< a > $ to $< a >$ is an automorphism for all $\n$.
Hence $x$ is a unit.
\end{proof}
 
We will see two examples one for $p = 2$ and other for $p$ an odd prime.
\begin{example}
\framebox(130,20)[c]{$(3^2)^x \eq 5^2 \, (mod \, 2^n)$}

\indent Let $p = 2$, $a = 3^2$, $b = 5^2$. Then, $(a, p) = (b, p) = 1.$\\
 Therefore, $3^2$ and $5^2$ both are in $\t.$

From the previous example, $<-3> = <5>$ in $\t$.

Therefore, $-3 = 5^k$ for some $k$. \\
$\Rightarrow (-3)^2 = 3^2 = 5^{2k} = (5^2)^k$\\
$\Rightarrow 3^2 \, \in <5^2>$\\
$\Rightarrow <3^2> \subseteq <5^2>.$

Also, $5 \, \in <-3>$\\
$\Rightarrow 5 = (-3)^k$ for some $k$.\\
$\Rightarrow 5^2 = (-3)^{2k} = 3^{2k} = (3^2)^k$\\
$\Rightarrow 5^2 \, \in <3^2>$\\
$\Rightarrow <5^2> \subseteq <3^2>.$

Order of 5 = $\w$\\
$\Rightarrow$ order of $5^2 = \v.$

Solution of the congruence $(3^2)^{x_n} \eq 5^2 \, (mod \, 2^n)$ is given by\\
\underline{using MATHEMATICA:}
\begin{table}[h]
\begin{tabular}{|c||c|c|c|c|c|c|c|c|c|c|c|c|c|c|c|c|}
 \hline
$n$ & 5 & 6 & 7 & 8 & 9 & 10 & 11 & 12 & 13 & 14 & 15 & 16 & 17 & 18 & 19 & 20 \\ [.5 ex]
\hline
$x_n$ & 3 & 3 & 11 & 11 & 11 & 11 & 11 & 267 & 267 & 1291 & 3339 & 7435 & 15627 & 15627 & 15627 & 15627 \\ [1 ex]
\hline
\end{tabular}
\end{table}

2-adic expansion of $x_n$ is as follows:\\
$3 = 1 + 2$\\
$11 = 1 + 2 + 2^3$\\
$267 = 1 + 2 + 2^3 + 2^8$\\
$1291 = 1 + 2 + 2^3 + 2^8 + 2^{10}$\\
$3339 = 1 + 2 + 2^3 + 2^8 + 2^{10} + 2^{11}$\\
$7435 = 1 + 2 + 2^3 + 2^8 + 2^{10} + 2^{11} + 2^{12}$\\
$15627 = 1 + 2 + 2^3 + 2^8 + 2^{10} + 2^{11} + 2^{12} + 2^{13}.$

\underline{using SAGE:}\\
the 2-adic solution is $x = 1 + 2 + 2^3 + 2^8 + 2^{10} + 2^{11} + 2^{12} + 2^{13} + O(2^{17}).$
\end{example}

\begin{example} \label{10}
\framebox(140,20)[c]{$(1 - p)^x \eq 1 + p \, (mod \, p^n)$} 

 Let $p$ be an odd prime. $a = 1 - p$, $b = 1 + p$.
 Consider the congruence $(1 - p)^x \eq 1 + p \, (mod \, p^n)$.\\
 Then, in group $\n$, order of $1 + p $ and $1 - p$ is  $p^{n-1}$\\
 (proved earlier, on page 7.) (proof for order of $1 - p$ is also similar.)
 
 $\Z_{p^n}^* \cong \Z_{p-1} \times \Z_{p^{n-1}}$\\
 $(p - 1, p^{n-1}) = 1 \Rightarrow \n$ is cyclic group. Therefore, for every divisor $d$ of the order of the group $\n$,
 there exists unique subgroup of order $d$. Here orders of subgroups $<1 + p>$ and $<1 - p>$ are same. hence, they must be equal.\\
\underline{ $<1 - p> = <1 + p>$.} Therefore, for every $x \, \in \{1, 2, \cdots, o(1 - p) = p^{n - 1}\}$,
there exists unique $y \, \in \{1, 2, \cdots, o(1 + p) = p^{n - 1}\}$ such that, $(1 - p)^x \eq (1 + p)^y \, (mod \, p^n)$.

Fix $y = 1.$ Let $x_n$ denote the solution of $(1 - p)^{x_n} \eq 1 + p \, (mod \, p^n)$, $x_n \in \{1, 2, \cdots, p^{n - 1}\}$.
If we express each $x_n$ as $x_o + x_1 p + x_2 p^2 + \cdots$ with $x_i \in \{0, 1, \cdots, p - 1\}.$\\
Then we get a p-adic solution of $(1 - p)^x \eq 1 + p \, (mod \, p^n)$\\
$\Rightarrow x \, ln \, (1 - p) = ln \, (1 + p)$
\begin{eqnarray*}
        \ln(1 - p) &=& -p - \dfrac{p^2}{2} - \dfrac{p^3}{3} - \cdots \\
	           &=& (p + \dfrac{p^2}{2} + \dfrac{p^3}{3} + \cdots) (-1) \\
                -1 &=& (p - 1) \dfrac{1}{1 - p} = (p - 1) (1 + p + p^2 + p^3 + \cdots)\\		   
        \ln(1 - p) &=& (p + \dfrac{p^2}{2} + \dfrac{p^3}{3} + \cdots) ((p -1) + (p -1)p + (p -1)p^2 + \cdots) \\
        \ln(1 + p) &= & p - \dfrac{p^2}{2} + \dfrac{p^3}{3} - \cdots
\end{eqnarray*}
Therefore, we can write,\\
$~~~~~~x \, ln \, (1 - p)= ln \, (1 + p)$\\
$ \Rightarrow (x_o + x_1 p + x_2 p^2 + \cdots) \, (p + \dfrac{p^2}{2} + \dfrac{p^3}{3} + \cdots)$\\
$= ((p -1) + (p -1)p + (p -1)p^2 + \cdots) \, (p - \dfrac{p^2}{2} + \dfrac{p^3}{3} - \cdots).$

We can find $x_k$ by comparing coefficients of $p^{k + 1}$ on both sides. We show here, calculation of $x_o, x_1$ and $x_2$. \\
Coefficient of $p$ gives \\
$x_o = p - 1$\\
Coefficient of $p^2$ gives \\
$\frac{x_o}{2} + x_1 = -\frac{p - 1}{2} + (p - 1)$\\
$\Rightarrow \frac{p - 1}{2} + x_1 = \frac{p - 1}{2}$\\
$\Rightarrow x_1 = 0$.\\
Coefficient of $p^3$ gives \\
$\frac{x_o}{3}+ \frac{x_1}{2} + x_2 = \frac{p - 1}{3} - \frac{p - 1}{2} + (p - 1)$\\
$\Rightarrow \frac{p - 1}{3} + x_2 = \frac{p - 1}{3} + \frac{p - 1}{2}$\\
$\Rightarrow x_2 = \frac{p - 1}{2}$.\\
Hence $x = (p-1) + 0* p + \frac{p - 1}{2} * p^2 + \cdots$.

\underline{Example 3(1):}

\underline{$(-4)^x \eq 6 \, (mod \, 5^n)$}

Here, $p = 5$, $a = 1 - p = -4$, $b = 1 + p = 6.$\\
We have already seen that, if $p$ is a prime then groups generated by $1 - p$ and $1 + p$ are equal.\\
Here, $<-4> = <6> \, \subseteq \Z_{5^n} ^*$ for all $n$. Therefore, for every $n$, there exists 
$x_n \, \in \, \{1, 2, \cdots, 5^{n - 1}\}$ such that, $(-4)^{x_n} \eq 6 \, (mod \, 5^n)$.

Using Mathematica we have calculated $x_n$ for first few $n$.
\begin{table}[h]
 \begin{tabular}{|c||c|c|c|c|c|c|c|c|c|c|}
\hline
$n$ & 1 & 2 & 3 & 4 & 5 & 6 & 7 & 8 & 9 & 10 \\ [.5 ex]
\hline
$x_n$ & 1 & 4 & 4 & 54 & 304 & 929 & 7179 & 22804 & 179054 & $x_9 + 2*5^8$ \\ [1 ex]
\hline
 \end{tabular}
\end{table}

We observe that, each $x_n$ can be written as follows:\\
$x_1 = 1* 5^0$\\
$x_2 = x_3 = 4* 5^0$\\
$x_4 = 4*5^0 + 2* 5^2$\\
$x_5 = 4*5^0 + 2* 5^2 + 2* 5^3$\\
$x_6 = 4*5^0 + 2* 5^2 + 2* 5^3 + 5^4$\\
$x_7 = 4*5^0 + 2* 5^2 + 2* 5^3 + 5^4 + 2*5^5$\\
$x_8 = 4*5^0 + 2* 5^2 + 2* 5^3 + 5^4 + 2*5^5 + 5^6$\\
$x_9 = 4*5^0 + 2* 5^2 + 2* 5^3 + 5^4 + 2*5^5 + 5^6 + 2*5^7$\\
$x_{10} = 4*5^0 + 2* 5^2 + 2* 5^3 + 5^4 + 2*5^5 + 5^6 + 2*5^7 + 2* 5^8.$

Using SAGE we can directly get 5-adic solution as:\\
$x = 4*5^0 + 2* 5^2 + 2* 5^3 + 5^4 + 2*5^5 + 5^6 + 2*5^7 + 2* 5^8 + 3*5^9 + 4*5^{10} + 2*5^{11} + 2*5^{12} + 5^{13} + 4*5^{14} + 3*5^{15} + 4*5^{17} + 3*5^{18} + O(5^{19}).$

Using SAGE p-adic solution can be found with precision $O(20)$, and precision cannot be changed.\\
Hence, MATHEMATICA is more useful if the solution is required with higher precision.
\end{example}

\section{$a$, $b$ units in $\zp$} \label{H}

 This is the last case where $a$ and $b$ are units in $p$-adic integers. Using the theory of \\ Teichm\"{u}ller units and the results in section \ref{E} we will 
find the condition on $a$, $b$, so that the equation $a^x  = b$ has a solution in $p$-adic integers. We record the relevant theorems here which are proved earlier:

Theorem \ref{7} says that, if $p$ is an odd prime then every unit in $p$-adic integers can be uniquely written as product of a Teichm\"{u}ller unit and an element of $U_1.$ \\
If $p = 2$ then every unit in $2$-adic integers can be uniquely written as product of a Teichm\"{u}ller unit and an element of $U_2.$ 

Another main theorem which we proved in section \ref{E} is theorem \ref{27} stated as:\\
\indent {\it{ ``Let $a, b $ be elements of $U_1$ and $a \neq 1$. Then $a^x = b$ has unique solution in $p$-adic integers if and only if depth of $a \leq$ depth of $b$.''}} 

We now state following important theorem in this section whose proof is clear from the above two theorems:

\begin{theorem} \label{28}
Let $p$ be an odd prime. Let $a, b$ be the units in $\zp.$ Write $a = a_1 a_2$, $b = b_1 b_2$ where $a_1, b_1$ are Teichm\"{u}ller units and $a_2, b_2 \- U_1$ then $a^x = b$ has a 
solution in $\zp$ if (i) $b_1$ belong to the group generated by $a_1,$ and (ii) depth of $a_2 \leq$ depth of $b_2.$ \\
Let $p = 2$, Let $a, b$ be the units in $\Z_{(2)}.$ Write $a = a_1 a_2$, $b = b_1 b_2$ where $a_1, b_1$ are Teichm\"{u}ller units and $a_2, b_2 \- U_1$ then $a^x = b$ has a 
solution in $\Z_{(2)}$ if (i) $a_1 = b_1$ and (ii) depth of $a_2 \leq$ depth of $b_2.$ 
\end{theorem}

\section{Special Pairs} \label{I}
Let $a \neq b$ be integers coprime to prime $p$. Therefore $a$, $b$ are units in the ring $\Z_{p^n}$ for all $n$. 
Fix natural number $n$. Consider the group of units $\n$. Then $< a >$ and $< b > $ are subgroups of $\n.$
If $< a > = < b > $ then for every $x \- \{1, 2, \cdots, o(a) = o(b)\}$ there exists unique $y \- \{1, 2, \cdots, o(a) = o(b)\}$
such that $a^x \eq b^y \, (mod \, p^n)$. \\
\indent Fix $y = 1$. Suppose $x_o \- \{1, 2, \cdots, o(a)\}$ is such that $a^{x_o} \eq b \, (mod \, p^n)$.
Clearly, $x_o$ must be coprime to $o(a)$. Consider the map from $< a > \rightarrow < a > $ given by $a \mapsto a^{x_o}$. Since $(x_o, o(a)) = 1$,
the map $\phi_{x_o} : 1 \mapsto x_o$ from $\Z_{o(a)}$ to itself is a group automorphism of $\Z_{o(a)}$, i.e. $\phi_{x_o} \- Aut \, (\Z_{o(a)})$.
We have seen that, the group $ Aut \, (\Z_{n})$ is isomorphic to the multiplicative group $\Z_n ^*$. We have $\phi_{x_o} \- Aut \, (\Z_{o(a)})$ correspondence to 
$x_o \- \Z_{o(a)} ^*.$

\begin{definition}
 Let $a, b$ be integers and $p$ be a fixed prime. We call $(a, b)$ \textbf{special pair mod $p^n$} if 
\begin{enumerate}
\item $a$ and $b$ are coprime to $p$.
\item subgroups of $\n$ generated by $a$ and $b$ are equal.
\item $x_o$ is of maximum order in the group $\Z_{o(a)} ^*$
\end{enumerate}
\end{definition}

\begin{example}
 Example of special pair is example \ref{9}, in which $p = 2$, $(a, b) =\\ (1 - p^2, 1 + p^2) = (-3, 5)$. We recall the $2$-adic solution of $(-3)^x = 5$ from the example \ref{9}.
$x = 1 + 2 + 2^3 + 2^8 + 2^{10} + 2^{11} + 2^{12} + 2^{13} + o(2^{18}).$
With the help of `SAGE' we calculate the multiplicative order of $x, \, mod \, o(a).$ Here $o(a) = 2^{n - 2}.$
Then $\Z_{o(a)} ^* = \Z_{2^{n - 2}} ^* \cong \Z_2 \times \Z_{2^{n - 4}}$ for $n \geq 5.$ Then maximum possible order of element in $\Z_{o(a)} ^*$ is $2^{n - 4}.$
\begin{table}[h]
 \begin{tabular}{|c|c|c|c|}
\hline
$n$ & $\Z_{o(a)} ^*$ & $x \, (mod \, o(a))$ & order of $x \, (mod \, o(a))$\\
\hline
5 & $\Z_{2^3} ^*$ & $1 + 2$ & 2\\
\hline
6 & $\Z_{2^4} ^*$ & $1 + 2 + 2^3$ & $2^2$\\
\hline
7 & $\Z_{2^5} ^*$ & $1 + 2 + 2^3$ & $2^3$\\
\hline
8 & $\Z_{2^6} ^*$ & $1 + 2 + 2^3$ & $2^4$\\
\hline
9 & $\Z_{2^7} ^*$ & $1 + 2 + 2^3$ & $2^5$\\
\hline
10 & $\Z_{2^8} ^*$ & $1 + 2 + 2^3$ & $2^6$\\
\hline
 \end{tabular}
\end{table}

\end{example}

\begin{example}
Another crucial example of special pair is example \ref{10}, where $p$ is an odd prime and $(a, b) = (1 - p, 1 + p)$. Here $\n$ is a cyclic group,
but order of $< a >$ and $< b >$ is $p^{n - 1}$ i.e. it is not the maximum possible order in the group $\n$. Since $\n \cong \Z_{p - 1} \times \Z_{\p}$,
there exists an element of order two viz. $-1$, so the pair $(a, b) = ((-1)(1 - p), (-1)(1 + p))$ is a better special pair in the sense that order of  
$(-1)(1 - p)$ and  $(-1)(1 + p)$ in $\n$ is $2 \p.$ \end{example}

Note that, we allow the representation of $\Z_n$ by negative integers here.

\indent The most important thing about these special pairs is that, they are described only in terms of the prime $p$ and not by the elements of the group $\n$.
Also the index of the groups $< a >$ and $< b >$ is eventually constant and the constant depends only on the prime $p$.

If $(a, b)$ is a special pair then $\phi_{x_o}$ (from condition 4 in the definition) give rise to the permutation of the set $\{1,2, \cdots, o(a) - 1\}$.
We observe that there is some interesting pattern in the permutation and there is some connection between these permutations on consecutive levels $n$, $n + 1$.
Consider the following example:

Example: Let $p = 2$, $(a, b)  = (-3, 5)$. Let $x_n$ denote the solution of $(-3)^{x_n} \eq 5 \, (mod \, 2^n)$. Then we get the permutation of the set $\{1, 2, \cdots, o(a) - 1\}$
by the map $1 \mapsto x_n$. Permutation is given by $1 \mapsto x_n \mapsto x_n ^2 \mapsto \cdots$. Here $x_n ^k$ is calculated modulo $o(a) = 2^{n - 2}$.
\begin{table}[h]
\begin{tabular}{|c|c|c|}
\hline
 $n$ & $\{1, 2, \cdots, 2^{n - 2} - 1\}$ & Permutation\\
\hline
3 & $\{1\}$& $(1)$\\
\hline
4 & $\{1,2, 3\}$ & $(1 \, 3) (2)$\\
\hline
5 & $\{1,2, \cdots, 7\}$ & $(1 \, 3)(2 \, 6) (4) (5 \, 7)$\\
\hline
6 & $\{1, 2, \cdots, 15\}$ & $( 1 \, 11 \, 9 \, 3) (2 \, 6)( 4 \, 12)( 5 \, 7 \, 13 \, 15)(8)(10 \, 14)$\\
\hline
7 & $\{1, 2, \cdots, 31\}$& $ (1 \, 11 \, 25 \, 19 \, 17 \, 27 \, 9 \, 3) \, (2 \, 22 \, 18 \, 6) \, ( 4 \, 12) (5 \, 23 \, 29 \, 31 \, 21 \, 7 \, 13 \, 15) $\\
\hline
\end{tabular}
\end{table}
%
%
%

\chapter{Summary and Conclusions}

In this concluding Chapter, we highlight the main results obtained in the present dissertation. Further, we discuss the possible extensions and the scope for further investigations
in this direction.

\section{Conclusions}

In this dissertation we have studied the congruences of the type $a^x \eq b \, (mod \, p^n )$ where $a, b$ are integers coprime to prime $p.$ We observed that, 
if the solution $x_n$ exists for all $n$ then the sequence $\{x_n\}$ converges in $p$-adic integers. Since $p$-adic integers is a collective way of thinking all $\Z_{p^n}$'s 
we considered the equation $a^x = b$ in $p$-adic integers. We mainly studied when does the solution of $a^x = b$ exists. We solve the case of $a, b \- U_1$ completely. 
The main theorem is: 

\noindent {\it{ ``Let $a, b $ be elements of $U_1$ and $a \neq 1$. Then $a^x = b$ has unique solution in $p$-adic integers if and only if depth of $a \leq$ depth of $b$.''}} 

In case when $a, b$ are units in $p$-adic integers, we used earlier case, $a, b \- U_1$ and the theory of Teichm\"{u}ller units to find the condition on $a, b$, 
when does the solution of $a^x = b$ exists. We arrive at following conclusion: 

\noindent {\it{``Let $p$ be an odd prime. Let $a, b$ be the units in $\zp.$ Write $a = a_1 a_2$, $b = b_1 b_2$ where $a_1, b_1$ are Teichm\"{u}ller units and
 $a_2, b_2 \- U_1$ then $a^x = b$ has a solution in $\zp$ if (i) $b_1$ belong to the group generated by $a_1,$ and (ii) depth of $a_2 \leq$ depth of $b_2.$ \\
Let $p = 2$, Let $a, b$ be the units in $\Z_{(2)}.$ Write $a = a_1 a_2$, $b = b_1 b_2$ where $a_1, b_1$ are Teichm\"{u}ller units and $a_2, b_2 \- U_1$ then $a^x = b$ has a 
solution in $\Z_{(2)}$ if (i) $a_1 = b_1$ and (ii) depth of $a_2 \leq$ depth of $b_2.$ '' }}
 
\section{Scope for further work}

In the case of $a$, $b$ are units in $p$-adic integers, in the theorem \ref{28}, we have obtained one way statement, which gives sufficient condition for the existence 
of the solution of equation $a^x = b.$ For finding the necessary condition or for checking that the sufficient condition is also necessary, we need to find the 
meaning of $a^x$ where $a$ is Teichm\"{u}ller unit and $x$ is a $p$-adic integer. Since $\log$ function is defined only on $1 + p \, \zp$, the 
tools we have used in this problem will not work to define Teichm\"{u}ller unit raised to a $p$-adic integer.

The problem of actual finding $p$-adic solution of $a^x = b$ even in case of integers is important further work. More precisely, we have shown in theorem \ref{25}:\\
\noindent {\it{``Let $a$, $b$ be integers coprime to a prime $p$ and $a \neq 1.$ Let $x_n$ denote the smallest positive solution of the congruence $a^x \eq b \, (mod \, p^n)$.
If $x_n$ exists for all $n \geq 1$ then the sequence $\{x_n\}$ converges in $p$-adic integers.'' }}\\
If sequence $\{x_n\}$ converges to the $p$-adic integer $x_o = \x$ then finding the general rule for the coefficient of $p^n$, is very important. 
Currently available mathematical softwares gives the power series representation of $p$-adic $\log (1 + a)$ upto only first few places. e.g. SAGE gives only upto $p^20.$

We have given the construction of Teichm\"{u}ller units from it's constant term. It will be interesting to find the actual formula for the coefficient of $p^n$ 
whose constant term is given. 

At the end of chapter 4, we have defined special pair in $\n.$ The work of finding all special pairs in $\n$ will be a big job. We can also try to generalize the notion 
 of special pairs to $p$-adic integers. It is our guess that, if $a, b$ are in $U_1$ then $(a, b)$ is a special pair if $(a \, mod \, p^n, b \, mod \, p^n)$ is a special pair in 
$\n$ for all $n.$ But for this we need better understanding of the condition (iii) in the definition of special pairs.  

We have defined $n$-adic integers, we can extend our results if we consider the congruences modulo $n^r$, where $n$ is a fixed natural number and $r$ varies in natural numbers.

%


\chapter{Projective Limit}

\section{ Definitions}

\begin{definition}
 $I$ is \textbf{directed set} if there exists partial order $\leq$ on $I$ such that every two elements in $I$ have an upper bound in $I$. \\
i.e. $\forall \, i, j \- I \Rightarrow \exists \, k \- I \, such \,that \, i \leq k, \, j \leq k $.
\end{definition}

\begin{definition}
 A \textbf{projective system} of sets (groups, rings, topological spaces) indexed by directed set $I$, is a pair $( \{X_i\}_{i \- I}, \{\phi_{ij}\}_{i \leq j} )$ 
where $ \{X_i\}_{i \- I} $ is a family of sets (groups, rings, topological spaces) and for every $i \geq j \- I$, 
$$\phi_{ij} : X_i \rightarrow X_j$$
is morphism of sets (group homomorphism, ring homomorphism,  continuous map), such that $ \phi_{ii} = 1_{X_i} $ and $ \phi_{ik} = \phi_{jk} \circ \phi_{ij}$
for all $i \geq j \geq k$.
$$\displaystyle \begin{xy}
*!C\xybox{
\xymatrix{
{X_i}\ar[rr]^{\phi_{ij}}\ar[dr]_{\phi_{ik}}&&{X_j}\ar[dl]^{\phi_{jk}}\\
&{X_k}&
} }
\end{xy}$$
\end{definition}
 
\begin{definition}
 Given a projective system $( \{X_i\}_{i \- I}, \{\phi_{ij}\}_{i \leq j} )$ the \textbf{ projective (or \\ inverse) limit} of the system is a pair $( X, \{\phi_i\}_{i \- I} )$
where $X$ is a set (group, ring, topological space) and $\phi_i : X \rightarrow X_i$ is a set map (group homomorphism, ring homomorphism,  continuous map) 
such that $\phi_{ij} \circ \phi_{i} = \phi_{j}$ for all $i \geq j$.
$$\displaystyle \begin{xy}
*!C\xybox{
\xymatrix{
{X}\ar[rr]^{\phi_{i}}\ar[dr]_{\phi_{j}}&&{X_i}\ar[dl]^{\phi_{ij}}\\
&{X_j}&
} }
\end{xy}$$
and such that this pair $( X, \{\phi_i\}_{i \- I} )$ is `universal'. i.e. Given any pair $( Y, \{\psi_i\}_{i \- I} )$ of the same type, there exists unique set map
(group homomorphism, ring homomorphism,  continuous map) $\psi : Y \rightarrow X$ such that $\phi_i \circ \psi = \psi_i$ for all $i \- I$.
$$\displaystyle \begin{xy}
*!C\xybox{
\xymatrix{
{Y}\ar[rr]^{\psi}\ar[dr]_{\psi_{i}}&&{X}\ar[dl]^{\phi_{i}}\\
&{X_i}&
} }
\end{xy}$$

As is well know, an object defined by Universal mapping property is unique if it exists. The following construction shows its existence.

\end{definition}
\textbf{
Construction of the projective limit}
$$X = \{(x_i)_{i \- I} \- \prod_{i \- I} X_i \, | \, x_j = \phi_{ij} (x_i), \, \forall \, i \geq j\}$$
We say that, $( X, \{\phi_i\}_{i \- I} )$ is the projective limit of the projective system $( \{X_i\}_{i \- I}, \{\phi_{ij}\}_{i \leq j} )$.
It is denoted as $X = \varprojlim X_i.$


%
\begin{example}
 Consider the family of sets $ X_{\alpha} \subset \Omega $ for $\alpha \- I$. For $\alpha, \beta \- I$ define $\beta \leq \alpha$ \iff $X_{\alpha} \subseteq X_{\beta}$.
Define maps $\phi_{\alpha \beta} : X_{\alpha} \rightarrow X_{\beta}$ as the inclusion map. Then $(\{X_{\alpha}\}_{\alpha \- I} , \phi_{\alpha \beta})$ is a projective system of 
sets. The projective limit of this system is $(X = \bigcap X_{\alpha} , \phi_{\alpha})$ where $\phi_{\alpha}$ is an inclusion map from $X$ to $X_{\alpha}$.
\end{example}

\begin{example}
 The most relevant example is : \\
Let $p$ be a fixed prime. Let $I = \N$ with usual order on $\N$. Consider the family of rings $R_n = \Z_{p^n}$. For $m \leq n$, define map 
$$\phi_{nm} : \Z_{p^n} \rightarrow \Z_{p^m}$$ 
$$~~~~~~~~~~~~ a \, mod \, p^{n} \mapsto a \, mod \, p^{m}$$
This is a projective system of rings. Projective limit of this system is the ring of $p$-adic integers $\zp$.
$$\zp = \{x = (\cdots, x_2, x_1) \, \in \, \prod_{n = 1} ^{\infty} \Z_{p^n}\, \mid \, x_n \, \in \, \Z_{p^n},\phi_{nm} (x_n) = x_m \, \, \forall n, m \}$$
 with map $\phi_n : \zp \rightarrow \Z_{p^n}$ as the projection onto $n^{th}$ co-ordinate.
$$ (\cdots, x_n, \cdots, x_2, x_1) \mapsto x_n$$
\end{example}

\chapter{$\exp$ and $\log$ maps}

\section{On $\Q[[x]]$}

In this section we define $\exp$ and $\log$ maps on the ring of formal power series $\Q[[x]]$. We will show that they are inverses of each other on some specific subsets of 
$\Q[[x]]$. 

Let $R= \Q[[x]] = $ ring of formal power series with rational coefficients. Then $R$ is a local ring with unique maximal ideal generated by $x$. \\
Let $M := M_1 = $ ideal generated by $x$. Let $M_n = $ ideal generated by $x^n$, where $n = 1, 2, \cdots.$\\
Then $U = R - M$ are units in $R$ that is, elements having non-zero constant term in their power series expression. Therefore, $U = \Q^*[[x]].$\\ 
Let $U_n = 1 + x^n R $ for $n = 1, 2, \cdots$. If $u$ is any element of $R$ then $u_o$ denote the constant term of $u$. Then 
\begin{eqnarray*}
U &\rightarrow& \Q ^*\\
u &\mapsto& u_o
\end{eqnarray*}
This is an onto group homomorphism with kernel $U_1$. Therefore $\frac{U}{U_1} \, \cong \, \Q^*$.\\
We also have a map 
\begin{eqnarray*}
 U_n &\rightarrow& \Q\\
1 + x^n \alpha &\mapsto& \alpha_o
\end{eqnarray*}
where $\alpha_o$ is the constant term of $\alpha.$ This is an onto group homomorphism with kernel $U_{n + 1}$. Therefore  $\frac{U_n}{U_{n + 1}} \, \cong \, \Q$ for all $n \geq 1$.\\
\newpage
Now we define $\exp$ and $\log$ map as follows: 
\begin{definition}
\begin{eqnarray*}
 \exp: M &\longrightarrow& U_1\\
         x \alpha  &\mapsto& 1 + x \alpha + \dfrac{(x \alpha )^2}{2!} +  \cdots + \dfrac{(x \alpha )^n}{n!} + \cdots
\end{eqnarray*}
\end{definition}
If $x^k \alpha \- R$ such that $\alpha_o \neq 0$ then in the power series expansion of $\exp (x^k \alpha)$, there are finitely many terms having non-zero coefficient of $x^n$ for
$n = 0, 1, 2, \cdots.$ Therefore $\exp$ map is well defined.

Now we define $\log$ map as follows:
\begin{definition}
\begin{eqnarray*}
 \log: U_1 &\longrightarrow& M\\
 1 + x \alpha &\mapsto& x \alpha - \dfrac{(x \alpha)^2}{2} +  \cdots + (-1)^{n + 1} \dfrac{(x \alpha)^n}{n} - \cdots
\end{eqnarray*}
\end{definition}
If $1 +  x^k \alpha \- R$ such that $\alpha_o \neq 0$ then in the power series expansion of $\log (1 + x^k \alpha)$, there are finitely many terms having non-zero 
coefficient of $x^n$ for $n = 0, 1, 2, \cdots.$ Therefore $\log$ map is well defined.

\begin{proposition}
\begin{eqnarray*}
 \exp: (M, + ) &\longrightarrow& (U_1, *)\\
         x \alpha  &\mapsto& 1 + x \alpha + \dfrac{(x \alpha )^2}{2!} +  \cdots + \dfrac{(x \alpha )^n}{n!} + \cdots
\end{eqnarray*}
is a group homomorphism. 
\end{proposition}
\begin{proof} Let $u, v \- M$.
 \begin{eqnarray*}
  \exp(u + v) &=& \sum_{n = 0} ^{\infty} \dfrac{(u + v)^n}{n!}\\
               &\vspace{0.2cm}&\\
              &=& \sum_{n = 0} ^{\infty} \dfrac{1}{n!} \, \sum_{r = 0} ^{n} \binom{n}{r} \, u^r v^{n - r}\\
               &\vspace{0.2cm}&\\
	      &=& \sum_{n = 0} ^{\infty} \, \sum_{r = 0} ^{n} \, \dfrac{1}{n!} \dfrac{n!}{(n - r)! \,r!} \, u^r v^{n - r}\\
	      &\vspace{0.2cm}&\\
              &=& \sum_{n = 0} ^{\infty} \, \sum_{r = 0} ^{n} \, \dfrac{u^r}{r!} \, \dfrac{v^{n - r}}{(n - r)!}\\
              &\vspace{0.2cm}&\\
	      &=& \sum_{r + s = 0} ^{\infty} \, \sum_{r = 0} ^{r + s} \, \dfrac{u^r}{r!} \, \dfrac{v^s}{s!}\\
             &\vspace{0.2cm}&\\
              &=& 1 + (u + v) + \dl {\dfrac{u^2}{2!} + uv + \dfrac{v^2}{2!} } \dr + \cdots  \\
             &\vspace{0.2cm}&\\
   &+& \dl{ \dfrac{u^{r + s}}{(r + s)!} +  \dfrac{u^{r + s - 1}}{(r + s - 1)!} v + \cdots + \dfrac{u^r}{r!} \, \dfrac{v^s}{s!} + \cdots + \dfrac{v^{r + s}}{(r + s)!} }\dr + \cdots\\
            &\vspace{0.2cm}&\\
            &=& \dl{ \sum_{r = 0} ^{\infty} \, \dfrac{u^r}{r!} }\dr \, \dl{ \sum_{s = 0} ^{\infty} \, \dfrac{v^s}{s!} }\dr\\
            &\vspace{0.2cm}&\\            
            &=& \exp u \, \exp v
 \end{eqnarray*}
\end{proof}

\begin{proposition}
 $\exp$ is an injective map.
\end{proposition}
\begin{proof}
 Suppose $\exp \, u = 1$ That is $1 + u + \dfrac{u^2}{2!} +  \cdots + \dfrac{u^n}{n!} + \cdots = 1$. This is true \iff $u = o$.
\end{proof}

\begin{remark}
 Showing that $\exp$ is surjective and $\log$ is it's inverse is a real tricky part. We will use the method of formal differentiation and integration to prove this part.
I am very much thankful to Prof. M. K. Srinivasan, IIT Bombay, for his kind and quick help in this proof.
\end{remark}

Let $f(x) = \A \- R$ and $g(x) = \B \- M $ then the composition $f \circ g$ is defined, it can be given as 
$$ f \circ g (p) = f(\B) = \sum_{i = 0} ^ {\infty} A_i \, (\B)^i$$
This sum is well defined, since there are only finitely many terms in $A_i, B_j$ which are coefficients of $x^n$ for all $n \geq 1$.\\
Note that, $g \circ f$ cannot be defined, as there will be infinite terms contributing to constant term. 

\textbf{Differentiation:} \\
Define $$D = \dfrac{d}{dx} : R \rightarrow R$$
$$ \dfrac{d}{dx} (\A) = A_1 + 2 A_2 x + \cdots + n A_n x^{n - 1} + \cdots$$
 as the formal differentiation having usual rules. We apply this definition to find derivative of $\exp$ and $\log$.
\begin{proposition}
 $D f(x) = 0 $ \iff $f(x)$ is a constant.
\end{proposition}

\begin{proof}
 Let $f(x) = \A$. Then $Df(x) = A_1 + 2 A_2 x + \cdots + n A_n x^{n - 1} + \cdots.$ \\
$Df(x) = 0$ \iff $A_1 = A_2 = \cdots = 0$.\\
Therefore $Df(x) = 0$ \iff $f(x) = A_o$
\end{proof}

\begin{proposition}
 $D(\exp(x)) = \exp(x)$
\end{proposition}

\begin{proof} \begin{eqnarray*}
  D(\exp(x)) &=& D \dl{  \sum_{i = 0} ^{\infty} \dfrac{x^n}{n!}  }\dr\\
            \vspace{0.2cm}&\\             
            &=& 1 + x + \dfrac{x^2}{2!} + \cdots + \dfrac{x^n}{n!} + \cdots \\
             \vspace{0.2cm}&\\             
            &=& \exp (x)
 \end{eqnarray*} \end{proof}

\begin{proposition}
 $D ( \log (1 + x)) = \dfrac{1}{1 + x}$
\end{proposition}

\begin{proof}
 \begin{eqnarray*}
  D ( \log (1 + x)) &=& D \dl { x - \dfrac{x^2}{2} + \dfrac{x^3}{3} - \cdots + (-1)^{n + 1} \dfrac{x^n}{n} + \cdots  } \dr\\
                     \vspace{0.2cm}&\\            
                     &=& 1 - x + x^2 - \cdots + (-1)^n x^n + \cdots\\
                     \vspace{0.2cm}&\\            
                      &=& \dfrac{1}{1 + x}
 \end{eqnarray*}
\end{proof}

\begin{proposition}
 Partial Chain Rule: $f(x) \- R$, $g(x) \- M $. \\
Then $D(f \circ g (x)) = Df(g(x)) \, Dg(x)$.
\end{proposition}

\begin{proof}
 \begin{eqnarray*}
  f \circ g (x) &=& f(\B)\\
                 &\vspace{0.2cm}&\\                 
                &=& \sum_{i = 0} ^{\infty}  { \dl{ \sum_{j = 1} ^{\infty} B_j x^j }\dr  }^i\\ 
                   &\vspace{0.2cm}&\\            
D(f \circ g (x)) &=& \sum_{i = 0} ^{\infty}  D { \dl{ \sum_{j = 1} ^{\infty} B_j x^j }\dr  }^i \\
                   & \vspace{0.2cm}&\\            
                &=&  \sum_{i = 1} ^{\infty} i \,  { \dl{ \sum_{j = 1} ^{\infty} B_j x^j }\dr  }^{i - 1} \, D { \dl {\B} \dr }  \\
                 &\vspace{0.2cm}&\\            
                &=& Df {\dl {\B} \dr} \, D { \dl {\B} \dr }\\
                 &\vspace{0.2cm}&\\            
                &=& Df(g(x)) Dg(x)
 \end{eqnarray*}
\end{proof}
\textbf{Integration:}
\begin{eqnarray*}
 \int: R &\rightarrow& M\\
       \A &\mapsto& \sum_{i = 0} ^{\infty} A_i \dfrac{x^{i + 1}}{i + 1}
\end{eqnarray*}

\begin{proposition}
 $D(\log \circ \exp (x)) = 1$
\end{proposition}

\begin{proof} Using chain rule,
\begin{eqnarray*}
 D(\log \circ \exp (x)) &=& D \dl \log {  \dl {\sum_{n = 0} ^{\infty} \dfrac{x^n}{n!}  } \dr} \dr \, D ( \exp (x) )\\
                        &\vspace{0.2cm}&\\          
                        &= &  \dfrac{1}{1 + \sum_{n = 1} ^{\infty} \dfrac{x^n}{n!} }  \, \exp (x)\\
                        &=& 1
\end{eqnarray*}
\end{proof}

Therefore $\log \circ \exp (x) = \int 1 = x$. This shows that $\exp$ is surjective. Also for $g(x) \- M$, we get $\log \circ \exp (g(x)) = g(x).$ Since $\exp$ is a group
homomorphism and $\log$ is it's left inverse, from the basic group theory result, we get $\log$ is also right inverse of $\exp$ map. Thus $\exp$ and $\log$ are isomorphisms of groups
$(U_1, *)$ and $(M, +)$. Moreover, it is clear from the definitions of $\exp$ and $\log$ that $\exp$ carries $M_n$ into $U_n$ and $\log$ carries $U_n$ into $M_n$.

\begin{remark}
 The proof of $\log \circ \exp = I$ only by using their power series expressions involves complicated calculations. Showing that coefficient of $x^n$, $n \geq 1$
is zero for the general $n$, contains a non-trivial identity in rational numbers. Above proof reduces all those efforts, hence is a very good proof.
\end{remark}


\begin{thebibliography}{}

\bibitem {Gau} C. F. Gauss,
{\em Disquisitiones Arithmeticae},
English Edition, Springer Verlag, (1985).

 \bibitem {B} D. Burton,
{\em Elementary Number Theory },
 Fourth Edition, The McGraw-Hill Companies, Inc. (1998).

\bibitem {DF} D. Dummit and R. Foote, 
{\em Abstract Algebra},
John Wiley and Sons, Inc. (2002).

\bibitem {Gou} F. Q. Gouv\^{e}a,
{\em $p$-adic Numbers, An Introduction},
Springer Verlag, (1991).

\bibitem{NZM} I. Niven, H. Zukerman, H. Montgomery,
{\em An Introduction to the Theory of Numbers},
Fifth Edition, John Wiley and Sons, Inc. (1991).

\bibitem {S} J. P. Serre,
{\em A Course in Arithmetic},
Springer International Student Edition, (1979).

\end{thebibliography}
\end{document}